\documentclass[final,leqno,pdftex]{siamltex}
\usepackage{comment,amsmath}

\usepackage{epsfig, graphicx}
\usepackage{makeidx}
\usepackage{multicol,multirow}
\usepackage{blkarray}
\usepackage{upgreek}
\usepackage{amsmath, amsfonts, amssymb, mathtools}
\usepackage{bbm, bm, dsfont}
\usepackage{ulem,enumitem}
\usepackage{crop,cancel}
\usepackage{xcolor,color}
\usepackage{todonotes}
\usepackage{hyperref,url}
\usepackage{subfig}
\usepackage[overload]{empheq}
\usepackage{cleveref,scalerel}
\usepackage{xspace,ifthen,mleftright,scalerel}
\usepackage{algorithmic}
\usepackage{xparse}
\usepackage{bbm}
\usepackage{tikz}
\usepackage{tikz-qtree}
\usepackage{tikz-cd}
\usepackage{subfiles}
\usepackage{stmaryrd}
\usepackage{pdflscape}
\usepackage{breqn}
\usepackage{rotating}
\usepackage{array}
\newcolumntype{C}[1]{>{\centering\arraybackslash}m{#1}}  % Vertically-centered & horizontally-centered
\newcolumntype{L}[1]{>{\raggedleft\arraybackslash}m{#1}} % Vertically-centered & horizontally-to-left centered
\newcolumntype{R}[1]{>{\raggedleft\arraybackslash}m{#1}} % Vertically-centered & horizontally-to-right

\usepackage[sort]{cite}
\newcommand{\tr}{{}^{\raisebox{1pt}{\scalebox{0.8}{$\scriptstyle \mathsf{T}$}}}}
\newcommand*{\tran}{^{\raisebox{1pt}{\scalebox{0.8}{$\scriptstyle \mathsf{T}$}}}}

\DeclareMathSizes{10}{10}{7}{4}
\DeclareMathSizes{11}{11}{8}{5}
\newcommand{\slow}{{\scalebox{0.7}{$\textsc{s}$}}}
\newcommand{\fast}{{\scalebox{0.7}{$\textsc{f}$}}}

\renewcommand{\slow}{{\textsc{s}}}
\renewcommand{\fast}{{\textsc{f}}}

\NewDocumentCommand{\component}{m m m}{^{ \left\{ \mkern-1mu #1 \mkern-1mu \right\} #2 \IfBooleanT{#3}{T} }}
\RenewDocumentCommand{\component}{m m m}{^{\left\{ \mkern-1mu \scalebox{0.65}{$#1$} \mkern-1mu \right\} #2 \IfBooleanT{#3}{T} }}
\NewDocumentCommand{\comp}{m O{} s}{\component{#1}{#2}{#3}}
\NewDocumentCommand{\F}{O{} s}{\component{\fast}{#1}{#2}}
\NewDocumentCommand{\FL}{O{\lambda} O{} s}{\component{\fast, \scalebox{0.8}{$#1$}}{#2}{#3}}
\RenewDocumentCommand{\S}{O{} s}{\component{\slow}{#1}{#2}}
\NewDocumentCommand{\SL}{O{\lambda} O{} s}{\component{\slow, \scalebox{0.8}{$#1$}}{#2}{#3}}
\NewDocumentCommand{\FF}{O{} s}{\component{\fast, \fast}{#1}{#2}}
\NewDocumentCommand{\FFL}{O{\lambda} O{} s}{\component{\fast, \fast, \scalebox{0.8}{$#1$}}{#2}{#3}}
\NewDocumentCommand{\FS}{O{} s}{\component{\fast, \slow}{#1}{#2}}
\NewDocumentCommand{\FSL}{O{\lambda} O{} s}{\component{\fast, \slow, \scalebox{0.8}{$#1$}}{#2}{#3}}
\NewDocumentCommand{\SF}{O{} s}{\component{\slow, \fast}{#1}{#2}}
\NewDocumentCommand{\SFL}{O{\lambda} O{} s}{\component{\slow, \fast, \scalebox{0.8}{$#1$}}{#2}{#3}}
\RenewDocumentCommand{\SS}{O{} s}{\component{\slow, \slow}{#1}{#2}}

\newcommand{\nvar}{d}
\newcommand{\nparts}{\mathrm{N}}
\newcommand{\fun}{\mathbf{f}}
\newcommand{\y}{\mathbf{y}}
\newcommand{\yy}{\y}

\newcommand{\yf}[1][]{%
   \ifthenelse{ \equal{#1}{} }
      {{\y\F}}
      {{\y\F_{#1}}}
}

\newcommand{\ys}[1][]{%
   \ifthenelse{ \equal{#1}{} }
      {{\y\S}}
      {{\y\S_{#1}}}
}

\newcommand{\Dt}{h}

\newcommand{\Id}{\mathbf{I}}

\newcommand{\one}{\mathbf{1}} 
\newcommand{\Zero}{\mathbf{0}}
\newcommand{\sfrac}[2]{\mbox{\footnotesize$\displaystyle\frac{#1}{#2}$}} % small frac

\makeatletter
\newcommand{\kron}[1]{{\,\mathrlap{\otimes}{\mathrlap{\hspace{0.275em}\tikz{\path[draw=white,fill=white] (0.025em,0.025em) rectangle (0.22em,0.38em);}}}\scalebox{0.4}{\raisebox{3.5pt}{\hspace{0.7em}#1}} \hspace{0.5em} }}
\makeatother

\def\Re{\mathds{R}}
\providecommand{\diag}{\operatorname{diag}}

\def\d{\mbox{d}}
\def\dt{\mbox{dt}}
\newcommand{\ddt}{\frac{\d}{\dt}}

\newcommand{\tree}{{\mathfrak{t}}}

\newcommand{\BT}{\mathds{T}}
\newcommand{\NT}{{\BT_{\nparts}}}

\newcommand{\Vertiii}{\vert\kern-0.22ex\vert\kern-0.22ex\vert}
\newcommand{\Vertiv}{\vert\kern-0.22ex\vert\kern-0.22ex\vert\kern-0.22ex\vert}

\newcommand{\A}{{\mathbf{A}}}

\renewcommand{\b}{{\mathbf{b}}}
\newcommand{\B}{{\mathbf{B}}}
\renewcommand{\c}{{\mathbf{c}}}

\newcommand{\Ahat}{{\widehat{\mathbf{A}}}}
\newcommand{\Bhat}{{\widehat{\mathbf{B}}}}
\newcommand{\bhat}{{\widehat{\mathbf{b}}}}
\newcommand{\chat}{{\widehat{\mathbf{c}}}}
\newcommand{\dhat}{{\widehat{\mathbf{d}}}}

\newcommand{\ahat}{{\widehat{a}}}

\newcommand{\dnohat}{{{\mathbf{d}}}}

\newcommand{\Ys}[1][]{%
   \ifthenelse{ \equal{#1}{} }
      {{\mathbf{Y}_{\slow}}}
      {{\mathbf{Y}_{\slow,#1}}}
}
\newcommand{\Yf}[1][]{%
   \ifthenelse{ \equal{#1}{} }
      {{\mathbf{Y}_{\fast}}}
      {{\mathbf{Y}_{\fast,#1}}}
}

\providecommand{\q}{\mathbf{q}}
\providecommand{\p}{\mathbf{p}}

\usepackage{pgfplots,makecell}
\usepackage{tikz-qtree}

\newcommand{\sandu}[1]{{\color{blue}\{Sandu: #1\}}}
\newcommand{\guenther}[1]{{\color{magenta}\{Günther: #1\}}}

\newif{\ifreport}
\reportfalse

\specialcomment{leaveout}
 {\begingroup\ttfamily\footnotesize\color{gray!50} }{\color{black}\endgroup}
\excludecomment{leaveout}

\graphicspath{ {Figures/}}

\newcommand{\ones}{\mathbf{1}}

\newtheorem{remark}{Remark}
\newtheorem{example}{Example}
\pgfplotsset{compat=1.16}
\begin{document}

\title{Symplectic GARK methods for partitioned Hamiltonian systems\thanks{Sandu was supported, in part, by the National Science Foundation through awards NSF ACI--1709727 and NSF CCF--1613905, by DOE ASCR through the award  DE--SC0021313, and by the Computational Science Laboratory at Virginia Tech. Günther and Schäfers were supported, in part, bei the Deutsche Forschungsgemeinschaft through Research Unit 5269 on "Future methods for studying confined gluons in QCD".}}

\author{
Michael G\"unther\thanks{Bergische Universit\"at Wuppertal,
        Institute of Mathematical Modelling, Analysis and Computational
        Mathematics (IMACM), Gauss strasse 20, D-42119 Wuppertal, Germany 
        ({\tt guenther@uni-wuppertal.de})} \and
        Adrian Sandu\thanks{Computational Science Laboratory, Department of Computer Science, 2202 Kraft Drive, Virginia Tech, Blacksburg, VA 24060, USA 
        ({\tt sandu@cs.vt.edu})}
        \and 
        Kevin Schäfers\thanks{Bergische Universit\"at Wuppertal,
        Institute of Mathematical Modelling, Analysis and Computational
        Mathematics (IMACM), Gauss strasse 20, D-42119 Wuppertal, Germany 
        ({\tt schaefers@math.uni-wuppertal.de})} \and
        Antonella Zanna\thanks{Matematisk institutt, Universitetet i Bergen, Norway 
        ({\tt Antonella.Zanna@uib.no})}
        }

\maketitle

\begin{abstract}
Generalized Additive Runge-Kutta schemes have shown to be a suitable tool for solving ordinary differential equations with additively partitioned right-hand sides. This work develops symplectic GARK schemes for additively partitioned Hamiltonian systems. In a general setting, we derive conditions for symplecticness, as well as symmetry and time-reversibility. We show how symplectic and symmetric schemes can be constructed based on schemes which are only symplectic, or only symmetric. Special attention is given to the special case of partitioned schemes for Hamiltonians split into multiple potential and kinetic energies. Finally we show how symplectic GARK schemes can leverage different time scales and evaluation costs for different potentials, and provide efficient numerical solutions by using different order for these parts.
\end{abstract}



\pagestyle{myheadings}
\thispagestyle{plain}
\markboth{G\"UNTHER, SANDU, SCHÄFERS AND ZANNA }{SYMPLECTIC GARK METHODS}

%%%%%%%%%%%%%%%%%%%%%%%%%%
\section{Introduction}
\label{dec:introduction}
%%%%%%%%%%%%%%%%%%%%%%%%%%
In many applications, initial value problems of ordinary differential equations are given as {\it additively} partitioned systems of the form: 
\begin{equation} 
\label{eqn:additive-ode}
 \y'= \fun(\y) = \sum_{m=1}^{\nparts} \fun^{\{m\}} (\y), \quad  t \ge t_0, \quad \y(t_0) = \y_0,
\end{equation}
where the right-hand side $\fun : \mathbb{R}^\nvar \rightarrow \mathbb{R}^\nvar$ is split into $\nparts$ different parts 
with respect to, for example, stiffness, nonlinearity, dynamical behavior, and evaluation cost. 

%Additive partitioning also includes the special case of 
%{\it coordinate} partitioning: 
%%
%\begin{equation} 
%\label{eqn:component-ode}
%y' = \begin{bmatrix} y^{\{1\}}  \\ \vdots \\ y^{\{\nparts\}} \end{bmatrix}'   
%= \sum_{m=1}^d  \begin{bmatrix} 0  \\ y^{\{m\}}      \\ 0 \end{bmatrix}'   
%= \sum_{m=1}^d  \begin{bmatrix} 0  \\ \fun^{\{m\}}(\y)  \\ 0 \end{bmatrix}.
%\end{equation}
%%
One step of a GARK method applied to~\eqref{eqn:additive-ode} advances the solution $\y_0$ 
at $t_0$ to the solution $\y_1$ at $t_1=t_0+h$ as follows:
%
%\begin{subequations}
%\label{eqn:GARK-symplectic}
%\begin{eqnarray}
%P_i^{\{q\}} & = & \p_0 + h \sum_{m=1}^{\nparts} \sum_{j=1}^{s^{\{m\}}} a_{i,j}^{\{q,m\}} k_j^{\{m\}}, \\
%Q_i^{\{q\}} & = & \q_0 + h \sum_{m=1}^{\nparts} \sum_{j=1}^{s^{\{m\}}} a_{i,j}^{\{q,m\}} \ell_j^{\{m\}},\\
%\p_1 & = & \p_0 + h \sum_{q=1}^{\nparts} \sum_{i=1}^{s^{\{q\}}} b_{i}^{\{q\}} k_i^{\{q\}}, \\
%\q_1 & = & \q_0 + h \sum_{q=1}^{\nparts} \sum_{i=1}^{s^{\{q\}}} b_{i}^{\{q\}} \ell_i^{\{q\}}, \\
%k_i^{\{m\}} & = & - H^{\{m\}}_{\q}\left(P_i^{\{m\}},Q_i^{\{m\}}\right), 
%\label{eqn:GARK-symplectic.stages}\\
%\ell_i^{\{m\}} & = & H^{\{m\}}_{\p}\left(P_i^{\{m\}},Q_i^{\{m\}}\right). 
%\label{eqn:GARK-symplectic.stage-l}
%\end{eqnarray}
%\end{subequations}
%
\begin{subequations}
\label{eqn:GARK}
\begin{align}
Y_i^{\{q\}}  & = \y_0  + h \sum_{m=1}^{\nparts} \sum_{j=1}^{s^{\{m\}}} a_{i,j}^{\{q,m\}} \fun^{\{m\}}(Y_j^{\{m\}}),\\
\y_1 & = 
\y_0  + h \sum_{q=1}^{\nparts} \sum_{i=1}^{s^{\{q\}}} b_{i}^{\{q\}} \fun^{\{q\}}(Y_i^{\{q\}}). 
\end{align}
\end{subequations} 

The general-structure additive Runge-Kutta (GARK) framework, developed in~\cite{Sandu_2015_GARK}, allows to construct multimethods that apply a different Runge-Kutta scheme, with possibly different time steps \cite{Sandu_2016_MR-GARK}, to discretize each component of \eqref{eqn:additive-ode}. The GARK framework explicitly reveals the structure of the multimethod in form of the numerical discretizations of individual components and coupling terms.  The GARK formalism allowed to construct new schemes such as implicit-implicit  \cite{Sandu_2015_GARK}, multirate methods of high order \cite{Sandu_2016_MR-GARK,Sandu_2019_MR-GARK_High-Order,Sandu_2021_MR-GARK_Implicit}, multirate infinitesimal step schemes \cite{Sandu_2019_MRI-GARK,Sandu_2020_MRI-GARK_Coupled}, and partitioned Rosenbrock methods \cite{Sandu_2021_GARK-ROS,Sandu_2021_MR-GARK_Implicit}. In addition, it was shown in \cite{Sandu_2022_GARK_splittingSchemes} that all classical splitting-based implicit time integration schemes can be understood as GARK methods. 

%\sandu{
Hamiltonian dynamics is fundamental to many fields in science and engineering. Symplectic integrators are special schemes for the numerical solution of Hamiltonian systems that preserve the geometric properties of the flow of the differential equation \cite{Hairer_2006_geometric-book}. Higher order symplectic integrators are typically constructed by ``operator splitting'', and symmetrically alternating fractional steps \cite{Yoshida_1990_splitting,Blanes_2003_composition}. Explicit symplectic schemes for partitioned Hamiltonians are also based on splitting the potential and kinetic energies \cite{SanzSerna_1992_symplectic,armusa97,SanzSerna_1992_symplectic-review}.

In this paper we discuss the GARK numerical solutions of split Hamiltonian systems \eqref{eqn:additive-ode} where each component $\fun^{\{m\}}$ may correspond to a Hamiltonian subsystem or not.  The main contributions of this work are as follows.
(i) Symplecticness, symmetry, and order conditions for GARK schemes applied to partitioned Hamiltonian systems are derived.  
(ii) The GARK formalism allows to consider very general partitions of Hamiltonian systems (e.g., splitting the Hamiltonian or splitting only the potential energy) of type $\fun(\y)=\mathbf{J}   \nabla H(\y)$ with an arbitrary, but skew-symmetric matrix $\mathbf J=- \mathbf J^\top$.
\begin{leaveout}
    Unlike traditional operator splitting approaches, in the GARK framework individual partitions do not have to be themselves Hamiltonian subsystems. Under suitable conditions, the overall method is symplectic even if individual components are not. This even applies to the more general case of  Hamiltonian flows of type $\fun(\y)=\mathbf{J}^{-1}   \nabla H(\y)$ with an arbitrary, but skew-symmetric matrix $\mathbf J=- \mathbf J^\top$.
\end{leaveout}
(iii) The GARK approach allows to integrate each component of the Hamiltonian with a different method, e.g., a high order method for the fast component and a low order method for the slow component. %Alternatively, each component can be integrated with a different time step, e.g., a small one for the fast component and a large one for the slow component. Such possibilities are illustrated in Section \ref{sec-multiratepotential}. \schaefers{would change the order of the last two sentences since MGARK is not illustrated in this paper.}
(iv) We construct symmetric and symplectic partitioned methods starting from partitioned symmetric (but not symplectic) schemes, or starting from  symplectic (but non-symmetric) methods. This is discussed in Section \ref{sec:symmetry-and-reversibility}. We show that explicit symplectic and symmetric partitioned GARK schemes are composition schemes. 
%}

The paper is organized as follows. Section \ref{dec:partitioned} reviews Hamiltonian systems, partitioned forms, and GARK schemes for the integration of partitioned systems. Section \ref{sec:symplectic-GARK}  introduces general symplectic GARK schemes. We derive conditions on the coefficients for symplecticity, which reduce the number of order conditions of GARK schemes drastically, and discuss symmetry and time-reversibility. If the Hamiltonians are split with respect to the potentials or kinetic parts and potentials, respectively,  partitioned versions of symplectic GARK schemes are tailored to exploit this structure. Section~\ref{sec-partitionedsymplecticgark}  introduces these schemes, with a discussion of symplecticity conditions, order conditions, symmetry and time-reversibility, as well as GARK discrete adjoints. Section~\ref{sec-multiratepotential} discusses how symplectic GARK schemes can exploit the multirate potential given by  potentials   of different activity levels. Numerical tests for a coupled oscillator are given. Section~\ref{sec-conclusions}  concludes with a summary.

%%%%%%%%%%%%%%%%%%%%%%%%%%
\section{Partitioned Hamiltonian systems and GARK schemes}
\label{dec:partitioned}
%%%%%%%%%%%%%%%%%%%%%%%%%%
%
A  Hamiltonian system is given by the ODE initial value problem
\begin{equation}
\label{eq-ham}
\y' = \mathbf J \, \nabla H(\y), \quad \y(t_0)=\y_0, \quad \mbox{with } \y = \begin{bmatrix}
\p \\ \q
\end{bmatrix}, \quad \mathbf J = \begin{bmatrix} \Zero_{d_p \times d_p} & -\Id_{d_p \times d_q} \\
\Id_{d_q \times d_p} & \hphantom{-}\Zero_{d_q \times d_q}   \end{bmatrix},
\end{equation}
where $d_p=d_q=d/2$, $q \in \mathbb{R}^{d_q}$ denote the generalized coordinates,  $p\in \mathbb{R}^{d_p}$ the conjugate momenta,  and $H: 
\mathbb{R}^{d_p} \times \mathbb{R}^{d_q} \rightarrow \mathbb{R}$ is a twice continuously differentiable Hamiltonian function.
%\sandu{Don't we always have $d_q = d_p$?}
%
The Hamiltonian flow $\y(t) = \varphi_t(\y_0)$, i.e., the solution to~\eqref{eq-ham}, is characterized by the following properties:
\begin{itemize}
\item The Hamiltonian is an invariant of the flow:
\begin{equation}
\label{inv.ham}
\ddt H(\varphi_t(\y_0))=0.
\end{equation}

\item The Hamiltonian is invariant with respect of changing the sign of momenta, $H(\p,\q) = H(-\p,\q)$. This can be formalized as follows:
\begin{equation}
\label{rhorev}
H = H \circ \rho \quad \textnormal{where} \quad
\rho \coloneqq \begin{bmatrix} -\Id_{d_p \times d_p} & \Zero_{d_p \times d_q} \\
 \hphantom{-} \Zero_{d_q \times d_p} & \Id_{d_q \times d_q}  \end{bmatrix},
\end{equation}
Consequently, the Hamiltonian equation of motion \eqref{eq-ham} are $\rho$-reversible, i.e., $ \rho \circ (\nabla H) = - \nabla (H \circ \rho) $.
%\begin{align*}
% \frac{\partial H(\q,\p)}{\partial \p} 
 %& = - \frac{\partial H(q,-p)}{\partial \p},
% \\
%\frac{\partial H(\q,\p)}{\partial \q}
%    & =       \frac{\partial H(q,-p)}{\partial \q},
%\end{align*}
%the Hamiltonian flow $\varphi_t$ is $\rho$-reversible:
\item The Hamiltonian flow is time-reversible: 
\begin{equation}
\label{eqn:timerev}
\rho \circ \varphi_t \circ \rho \circ  \varphi_t (\y_0) = \y_0
\quad \Leftrightarrow \quad
\rho \circ  \varphi_t = \varphi_{-t}  \circ \rho,  %\qquad \rho = \begin{bmatrix} \Id_{d_q \times d_q} & ~~~\Zero_{d_q \times d_p} \\
%\Zero_{d_p \times d_q} & ~-\Id_{d_p \times d_p}  \end{bmatrix}, 
\end{equation}
%
%which is equivalent to 
%%
%\begin{equation}
% \rho \circ  \varphi_t = \varphi_{-t}  \circ \rho
%\end{equation}
where the second equivalent equation is due to the symmetry $\varphi_t \circ \varphi_{-t} (\y_0)=\y_0$ of the flow.
\item The Hamiltonian flow is symplectic:
\begin{equation}
\label{equation:symplectic}
\left(  \frac{\partial{\varphi_t(\y_0)}}{\partial \y_0 }
\right)^{\mkern-2\thickmuskip \top} \, \mathbf{J^{-1}}\,  
\left(  \frac{\partial{\varphi_t(\y_0)}}{\partial \y_0 }
\right) =\mathbf{J^{-1}},
\end{equation}
and thus volume-preserving 
\begin{equation}
\label{eqn:volpres}
\left| \det \left( \frac{\partial{\varphi_t(\y_0)}}{\partial \y_0 }
\right) \right| =1.
\end{equation}
\end{itemize}
In geometric integration, we demand the mapping $\y_0 \mapsto \Phi_t(\y_0)$ defining the numerical approximation $\Phi_t(\y_0) \approx \varphi_t(\y_0)$ to be time-reversible and symplectic as well:
\begin{subequations}
\label{eqn:numerical-properties}
\begin{eqnarray}
\label{eqn:numerical-reversible}
&& \rho \circ \Phi_t \circ \rho \circ \Phi_t(\y_0) = \y_0,  \\
\label{eqn:numerical-symplectic}
&& \left(  \frac{\partial{\Phi_t(\y_0)}}{\partial \y_0 }
\right)^{\mkern-2\thickmuskip \top}\, \mathbf{J^{-1}}  \,
\left(  \frac{\partial{\Phi_t(\y_0)}}{\partial \y_0 }
\right) = \mathbf{J^{-1}}.
\end{eqnarray}
\end{subequations}
\begin{remark}
\begin{leaveout}
    If we replace $\mathbf{J}$ in~\eqref{eq-ham} by an arbitrary regular skew-symmetric matrix, the invariance of the Hamiltonian~\eqref{inv.ham} \cite{Gonzalez1996TimeIA,McLachlan} and the symplecticeness~\eqref{eqn:numerical-symplectic} of the flow (here the proof of Theorem 2.4 in~\cite{Hairer_2006_geometric-book} directly generalizes to regular skew-symmetric matrices)  still hold.  In this case, however, the unknowns $\q$ and $\p$ might loose their meaning as generalized coordinates and positions of classical mechanics.
\end{leaveout}
If we replace $\mathbf{J}$ in~\eqref{eq-ham} by an arbitrary regular skew-symmetric matrix, the invariance of the Hamiltonian~\eqref{inv.ham} \cite{Gonzalez1996TimeIA,McLachlan} and the symplecticeness~\eqref{equation:symplectic} of the flow (here the proof of Theorem 2.4 in~\cite{Hairer_2006_geometric-book} directly generalizes to regular skew-symmetric matrices)  still hold, as well as volume-preservation~\eqref{eqn:volpres}. In this case, however, the unknowns $\q$ and $\p$ might loose their meaning as generalized coordinates and positions of classical mechanics, and time-reversibility~\eqref{eqn:timerev} loses its significance.
%\guenther{One notes that volume-preservation holds also in the case of singular skew-symmetric matrices}. 
\end{remark}

\subsection{GARK schemes for Hamiltonian systems}

Consider a general splitting of the right-hand side of the type  
\begin{equation}
\label{eqn:partitioned-sympl.system-without-structure}
\renewcommand{\arraystretch}{1.2}
\y'
 =
\sum_{m=1}^{\nparts}
\fun^{\{m\}}(\y).
\end{equation}
One step of a GARK method~\eqref{eqn:GARK} applied to~\eqref{eqn:partitioned-sympl.system-without-structure} advances the solution $(\y_0)$ 
at $t_0$ to the solution $(\y_1)$ at $t_1=t_0+h$ as follows:
%
%\begin{subequations}
%\label{eqn:GARK-symplectic}
%\begin{eqnarray}
%P_i^{\{q\}} & = & \p_0 + h \sum_{m=1}^{\nparts} \sum_{j=1}^{s^{\{m\}}} a_{i,j}^{\{q,m\}} k_j^{\{m\}}, \\
%Q_i^{\{q\}} & = & \q_0 + h \sum_{m=1}^{\nparts} \sum_{j=1}^{s^{\{m\}}} a_{i,j}^{\{q,m\}} \ell_j^{\{m\}},\\
%\p_1 & = & \p_0 + h \sum_{q=1}^{\nparts} \sum_{i=1}^{s^{\{q\}}} b_{i}^{\{q\}} k_i^{\{q\}}, \\
%\q_1 & = & \q_0 + h \sum_{q=1}^{\nparts} \sum_{i=1}^{s^{\{q\}}} b_{i}^{\{q\}} \ell_i^{\{q\}}, \\
%k_i^{\{m\}} & = & - H^{\{m\}}_{\q}\left(P_i^{\{m\}},Q_i^{\{m\}}\right), 
%\label{eqn:GARK-symplectic.stages}\\
%\ell_i^{\{m\}} & = & H^{\{m\}}_{\p}\left(P_i^{\{m\}},Q_i^{\{m\}}\right). 
%\label{eqn:GARK-symplectic.stage-l}
%\end{eqnarray}
%\end{subequations}
%
\begin{subequations}
\label{eqn:GARK2}
\begin{align}
\label{eqn:GARK-symplectic.stage-solutions2}
 Y_i^{\{q\}} 
& = \y_0 + h \sum_{m=1}^{\nparts} \sum_{j=1}^{s^{\{m\}}} a_{i,j}^{\{q,m\}} 
k_j^{\{m\}}, \quad q=1,\dots,\nparts, \\
\label{eqn:GARK-symplectic.stage-functions2}
k_i^{\{m\}}  & \coloneqq \fun^{\{m\}}\left(Y_i^{\{m\}}\right), \\
\label{eqn:GARK-symplectic.solutions2}
\y_1
& = 
\y_0 + h \sum_{q=1}^{\nparts} \sum_{i=1}^{s^{\{q\}}} b_{i}^{\{q\}} 
k_i^{\{q\}}. 
\end{align}
\end{subequations}
The corresponding generalized Butcher tableau is:
\begin{equation}
\label{eqn:general-Butcher-tableau}
\renewcommand{\arraystretch}{1.3}
\begin{array}{c}
\mathbf{A}_{{\textsc{gark}}} \\ 
\Xhline{2\arrayrulewidth}
\mathbf{b}_{{\textsc{gark}}}\tran
\end{array}
~~=~~
\raisebox{17pt}{$
\begin{array}{cccc}
\mathbf{A}^{\{1,1\}} &  \cdots & \mathbf{A}^{\{1,\nparts\}} \\
%\mathbf{A}^{\{2,1\}} & \cdots & \mathbf{A}^{\{2,\nparts\}} \\
\vdots  &\ddots & \vdots \\
\mathbf{A}^{\{\nparts,1\}} &  \cdots & \mathbf{A}^{\{\nparts,\nparts\}} \\ 
\Xhline{2\arrayrulewidth}
\mathbf{b}^{\{1\}}\tr &  \cdots &\mathbf{b}^{\{\nparts\}}\tr
\end{array}$}.
\end{equation}
In contrast to traditional additive methods~\cite{Kennedy_2003} different stage values are used with different components of the right hand side. 
The methods $(\mathbf{A}^{\{q,q\}},\mathbf{b}^{\{q\}})$ can be regarded as stand-alone integration schemes
applied to each individual component $q$. The off-diagonal matrices $\mathbf{A}^{\{q,m\}}$, $m \ne q$, can be 
viewed as a coupling mechanism among components.

We define the abscissae associated with each tableau as $\mathbf{c}^{\{q,m\}} \coloneqq \mathbf{A}^{\{q,m\}}\cdot \one^{\{m\}}$, where $\one^{\{m\}} \in \Re^{s^{\{m\}}}$ is a vector of ones. The method \eqref{eqn:general-Butcher-tableau} is {\it internally consistent} \cite{Sandu_2015_GARK} if all the abscissae along each row of blocks coincide:
\[
\mathbf{c}^{\{q,1\}} = \cdots = \mathbf{c}^{\{q,\nparts\}} \eqqcolon \mathbf{c}^{\{q\}}, \quad q=1,\dots,\nparts.
\]
A particular splitting case is offered by partitioned Hamiltonian systems, which are 
characterized by a Hamiltonian function $H(\y)$ split into $\nparts$ individual Hamiltonians:
\begin{equation}
\label{eqn:partitioned-Hamiltonians}
H(\y)=\sum_{m=1}^{\nparts} H^{\{m\}}(\y).
\end{equation}
Consequently, the equations of motion are a partitioned system \eqref{eqn:partitioned-sympl.system-without-structure} where each component function corresponds to one individual Hamiltonian: 
\begin{equation}
\label{eqn:partitioned-sympl.system}
\renewcommand{\arraystretch}{1.2}
\y'
= \mathbf{J}\, \nabla H(\y) = \sum_{m=1}^{\nparts}
\fun^{\{m\}}(\y), \qquad 
\fun^{\{m\}}(\y) \coloneqq
\mathbf{J}\, \nabla H^{\{m\}}(\y).
\end{equation}   
An efficient numerical integration scheme needs to exploit the different properties of the $\nparts$ individual Hamiltonians, such as slow dynamics with expensive evaluation costs versus fast dynamics with cheap evaluation costs, while preserving time-reversibility and symplecticity. One class of numerical schemes tailored to exploiting different right-hand side component properties are partitioned GARK schemes: one step of a GARK method \eqref{eqn:GARK} applied to~\eqref{eqn:partitioned-sympl.system} advances the solution $(\y_0)$ 
at $t_0$ to the solution $(\y_1)$ at $t_1=t_0+h$ as 
given in~\eqref{eqn:GARK2}, with~\eqref{eqn:GARK-symplectic.stage-functions2} replaced by 
\begin{equation}
    \label{eqn:GARK-symplectic}
    k_i^{\{m\}}  = 
\mathbf{J}\, \nabla H^{\{m\}}(Y_i^{\{m\}} ).
\end{equation}

\begin{leaveout}
\begin{subequations}
\label{eqn:GARK-symplectic2}
\begin{align}
\label{eqn:GARK-symplectic.stage-solutions}
Y_i^{\{q\}} 
 & = 
\y_0
 + h \sum_{m=1}^{\nparts} \sum_{j=1}^{s^{\{m\}}} a_{i,j}^{\{q,m\}} 
k_j^{\{m\}}, \\
\label{eqn:GARK-symplectic.stage-functions}
k_i^{\{m\}}  & = 
\mathbf{J}^{-1}\, \nabla H^{\{m\}}(Y_i^{\{m\}} ), \\
\label{eqn:GARK-symplectic.solutions}
\y_1
& = 
\y_0
 + h \sum_{q=1}^{\nparts} \sum_{i=1}^{s^{\{q\}}} b_{i}^{\{q\}} 
k_i^{\{q\}}. 
\end{align}
\end{subequations}

\sandu{
A non-Hamiltonian splitting of a Hamiltonian system can be obtained via partitioning the $\mathbf{J}$ matrix, in which case the system \eqref{eqn:partitioned-sympl.system-without-structure} reads
\begin{equation}
\label{eqn:partitioned-nonHamiltonians-1}
\quad \fun^{\{m\}}(\y) \coloneqq
\mathbf{J}^{-1}\,\big( \nabla H^{\{m\}}(\y) + \mathbf{g}^{\{m\}}(\y) \big), \quad
\sum_{m=1}^{\nparts} \mathbf{g}^{\{m\}}(\y) \equiv 0.
\end{equation}
The form \eqref{eqn:partitioned-nonHamiltonians-1}  gave some positive results with KdV. Can we show symplecticness based not on B-series, but using the $\mathbf{J}^{-1}$ form directly?

From \eqref{equation:symplectic} we need to show that:
\begin{equation}
\left(  \frac{\partial \y_1}{\partial \y_0 }
\right)^{\mkern-2\thickmuskip \top} \, \mathbf{J}\,  
\left(  \frac{\partial \y_1}{\partial \y_0 }
\right) =\mathbf{J}.
\end{equation}
The Jacobian of the solution is:
\begin{align}
\frac{\partial Y_i^{\{q\}}}{\partial \y_0 } & = 
\Id + h \sum_{m=1}^{\nparts} \sum_{j=1}^{s^{\{m\}}} a_{i,j}^{\{q,m\}} 
\frac{\partial k_j^{\{m\}}}{\partial \y_0 }  , \\
\frac{\partial k_i^{\{q\}}}{\partial \y_0 }   & = 
\mathbf{J}^{-1}\, \nabla \fun^{\{q\}}_i\,\frac{\partial Y_i^{\{q\}}}{\partial \y_0 }, \\
\frac{\partial \y_1}{\partial \y_0 } & = 
\Id + h \sum_{q=1}^{\nparts} \sum_{i=1}^{s^{\{q\}}} b_{i}^{\{q\}} 
\frac{\partial k_i^{\{q\}}}{\partial \y_0 }. 
\end{align}
In matrix notation:
\begin{align}
\frac{\partial Y}{\partial \y_0 } & = 
\one_{s \nparts} \otimes \Id + h \mathbf{A} \kron{n}
\frac{\partial k}{\partial \y_0 }  , \\
\frac{\partial k}{\partial \y_0 }   & = 
\diag_{1 \le q \le \nparts} \left(\diag_{1 \le i \le s} \;\mathbf{J}^{-1}\, \nabla \fun^{\{q\}}_i\right)\,\frac{\partial Y}{\partial \y_0 } \\
& = 
\diag_{1 \le q \le \nparts, 1 \le i \le s} \;\mathbf{J}^{-1}\, \nabla \fun^{\{q\}}_i \,\frac{\partial Y}{\partial \y_0 }, \\
\frac{\partial \y_1}{\partial \y_0 } & = 
\Id + h \, \mathbf{b}^{T} \kron{n}
\frac{\partial k}{\partial \y_0 }. 
\end{align}

Another non-Hamiltonian splitting of a Hamiltonian system can be obtained via partitioning the $\mathbf{J}$ matrix, in which case the system \eqref{eqn:partitioned-sympl.system-without-structure} reads
\begin{equation}
\label{eqn:partitioned-nonHamiltonians-2}
\mathbf{J}^{-1} = \sum_{m=1}^{\nparts} \mathbf{J}^{-1\{m\}},
\quad \fun^{\{m\}}(\y) \coloneqq
\mathbf{J}^{-1\{m\}}\, \nabla H(\y).
\end{equation}
}
\end{leaveout}
%
%Examples of Hamiltonian and non-Hamiltonian splittings are given next.
An example of a system with skew-symmetric, but singular $\mathbf{J}$ is given next.%
\begin{example}[Splitting of Hamiltonian for a two mass oscillator]
\label{example.two.mass.osc}
%%%
\begin{figure}[htb]
\begin{center}
\includegraphics[width=0.75\textwidth]{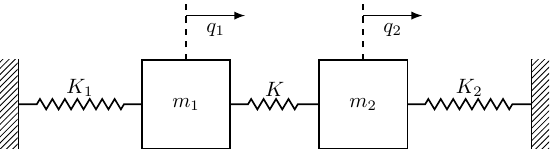}
\end{center}
\caption{\label{fig:twomasses-threesprings} The two masses oscillator.}
\end{figure}
%%%

Consider the following one-dimensional mechanical system consisting of two masses and three linear springs shown in Fig.~\ref{fig:twomasses-threesprings}. It has the Hamiltonian 
\[ 
H=\frac{1}{2} \left( \frac{p_1^2}{m_1} + \frac{p_2^2}{m_2} + K_1 q_1^2 + K (q_1-q_2)^2\  + K_2 q_2^2 \right).
\]
If we split the system into two subsystems consisting of elements $(K_1,m_1,K)$ and $(m_2,K_2)$, the  $H$ is split into the two Hamiltonians $H_1$ and $H_2$ of the subystems with
\[
    H_1(p_1,q_1,q_1-q) = \frac{1}{2} \left( \frac{p_1^2}{m_1} + K_1 q_1^2 + K (q_1-q)^2\right), \quad H_2(p_2,q_2) = \frac{1}{2} \left( \frac{p_2^2}{m_2} + K_2 q_2^2 \right),
    \]
respectively, where $q$ is a port variable that defines the coupling parameter from the first to the second system; for the coupling configuration above $q=q_2$. Now the dynamics can be defined by two coupled port-Hamiltonian systems (see ~\cite{paper_BGJR}), which yields in condensed form with $x=(x_1,x_2)^\top$,  $x_1=(p_1,q_1,q_1-q)^\top$ and $x_2=(p_2,q_2)^\top$: 
    \begin{eqnarray}
    \label{ex.phs}
    \dot x = \mathbf{J} \cdot \nabla H(x), \quad  \mathbf{J}= \begin{bmatrix}
    \mathbf{J}_1 & \mathbf{B} \\ -\mathbf{B}^\top & \mathbf{J}_2
    \end{bmatrix}, 
    \end{eqnarray}
    with $\nabla$ denoting the derivative with respect to $x$ and 
    \[ 
    \mathbf{J}_1 = \begin{bmatrix} 0 & -1 & -1 \\ 1 & 0 & 0 \\ 1 & 0 & 0
    \end{bmatrix}, \quad
    \mathbf{J}_2 = \begin{bmatrix} 0 & -1  \\ 1 & 0
    \end{bmatrix}, \quad
    \mathbf{B}=\begin{bmatrix}
    0 & 0 \\ 0 & 0 \\ -1 & 0
    \end{bmatrix}.
    \]
    A Hamiltonian  splitting  is given by
    \[
    f^{\{1\}}(x_1,x_2) =\begin{bmatrix}
    \mathbf{J}_1 & \mathbf{B} \\ -\mathbf{B}^\top & \mathbf{J}_2
    \end{bmatrix} \nabla H_1(x_1), \quad f^{\{2\}}(x_1,x_2)=\begin{bmatrix}
    \mathbf{J}_1 & \mathbf{B} \\ -\mathbf{B}^\top & \mathbf{J}_2
    \end{bmatrix} \nabla H_2(x_2).
    \]
    One may also split $\mathbf{J}$ instead of the Hamiltonian $H$ and obtain either a non-Hamiltonian (component-wise) splitting by 
    \begin{eqnarray*}
    f^{\{1\}}(x_1,x_2) &=& \begin{bmatrix}
    \mathbf{J}_1 & \mathbf{B} \\ 0 & 0
    \end{bmatrix}
    \cdot \nabla H(x), \\
   % \nabla_x \big(H_1(x_1)+H_2(x_2)\big), \\
    \quad
    f^{\{2\}}(x_1,x_2) &=& \begin{bmatrix}
    0 & 0  \\ -\mathbf{B}^\top & \mathbf{J}_2
    \end{bmatrix}
    \cdot
    \nabla H(x), 
    %\nabla_x \big(H_1(x_1)+H_2(x_2)\big).
    \end{eqnarray*}
  or a (non-)Hamiltonian splitting with respect to subsystems and coupling parts
   \begin{eqnarray*}
    f^{\{1\}}(x_1,x_2) &=& \begin{bmatrix}
    \mathbf{J}_1 & 0\\ 0 & \mathbf{J}_2
    \end{bmatrix}
    \cdot \nabla H(x), \\
   % \nabla_x \big(H_1(x_1)+H_2(x_2)\big), \\
    \quad
    f^{\{2\}}(x_1,x_2) &=& \begin{bmatrix}
    0 & \mathbf{B}  \\ -\mathbf{B}^\top & 0
    \end{bmatrix}
    \cdot
    \nabla H(x).
    %\nabla_x \big(H_1(x_1)+H_2(x_2)\big).
    \end{eqnarray*}
    \end{example}
\color{black}
{
\begin{remark}
    Note that the skew-symmetric matrix $\mathbf{J}$ in~\eqref{ex.phs} is singular. Nevertheless, it fulfills a symplectic structure as we can see as follows: as the Hamiltonian is quadratic, \eqref{ex.phs} can be written as $\dot x = \mathbf{J} Qx$ with $Q$ positive-definite, %semidefinit, 
    which can be transformed into
    $\dot w = \widetilde{\mathbf{J}} w$ with $w:=Q^{1/2} x$ and $\widetilde{\mathbf{J}}= Q^{1/2} \mathbf{J} Q^{1/2}$. 
    Consider now the variational equations 
     $\dot \Phi= \mathbf{J} Q \Phi$, $\Phi(0)=\Id_{5 \times 5}$
    and $ \dot \Psi= \widetilde{\mathbf{J}}\Psi$, $\Psi(0)=\Id_{5 \times 5}$ of the original and tranformed system, resp. 
    Defining $\mathbf{J}^+$ as the Drazin inverse of $\mathbf{J}$ and $\widetilde{\mathbf{J}}^+$ as the Drazin inverse of $\widetilde{\mathbf{J}}$, we get on the one hand
    \begin{eqnarray*}
        \ddt \left( \Psi^\top \widetilde{\mathbf{J}}^+ \Psi\right) & = &  \left( \dot \Psi^\top \widetilde{\mathbf{J}}^+ \Psi + \Psi^\top \widetilde{\mathbf{J}}^+ \dot \Psi\right) \\
        & = &  (\widetilde{\mathbf{J}}\Psi)^\top \widetilde{\mathbf{J}}^+ \Psi + \Psi^\top \widetilde{\mathbf{J}}^+ (\widetilde{\mathbf{J}}\Psi)  \\
        & = & \Psi^\top \left( \widetilde{\mathbf{J}}^\top \widetilde{\mathbf{J}}^+ + \widetilde{\mathbf{J}}^+ \widetilde{\mathbf{J}} \right) \Psi \\
        & = & \Psi^\top \left( - \widetilde{\mathbf{J}} \widetilde{\mathbf{J}}^+ + \widetilde{\mathbf{J}}^+  \widetilde{\mathbf{J}} \right) \Psi \\ &&= 0,
    \end{eqnarray*}
    i.e., $ \Psi^\top \widetilde{\mathbf{J}}^+ \Psi = \widetilde{\mathbf{J}}^+$ is a quadratic invariant. On the other hand we have with 
    $\Psi= Q^{1/2} \Phi$
   \begin{eqnarray*}
        \Psi^\top \widetilde{\mathbf{J}}^+ \Psi &=& (Q^{1/2} \Phi)^\top (Q^{1/2} \mathbf{J} Q^{1/2})^+ (Q^{1/2} \Phi) \\
        & = & \Phi^\top Q^{1/2} Q^{-1/2} \mathbf{J}^+ Q^{-1/2} Q^{1/2} \Phi \\
       & = & \Phi^\top \mathbf{J}^+ \Phi,
    \end{eqnarray*}
   which shows that $ \Phi^\top \mathbf{J}^+ \Phi = \mathbf{J}^+$ is a quadratic invariant, too. Summing up, for a singular skew-symmetric matrix $\mathbf{J}$ the symplectic structure given by~\eqref{equation:symplectic} holds, if one replaces $\mathbf{J}^{-1}$ by the Drazin inverse of $\mathbf{J}$.
\end{remark}
}

%%%%%%%%%%%%%%%%%%%%%%%%%%
%%%%%%%%%%%%%%%%%%%%%%%%%%%%%%

\section{Symplectic GARK schemes -- the general case}
\label{sec:symplectic-GARK}
%%%%%%%%%%%%%%%%%%%%%%%%%%%%%%
In this section we consider the general case of GARK schemes~\eqref{eqn:GARK2} applied to a Hamiltonian system~\eqref{eq-ham} based on a general splitting~\eqref{eqn:partitioned-sympl.system-without-structure}.

Several matrices are defined from the coefficients of \eqref{eqn:GARK} for $m,\ell = 1,\dots,\nparts$:
\begin{subequations}
\label{eqn:symplecticness-matrix}
\begin{eqnarray}
\label{eqn:B-matrix}
\mathbf{B}^{\{m\}} &\coloneqq& \mbox{diag}\bigl(\mathbf{b}^{\{m\}}\bigr) \in \Re^{s^{\{m\}} \times s^{\{m\}}}, \\
\label{eqn:P-matrix-blocks}
\mathbf{P}^{\{m,\ell\}} &\coloneqq& \mathbf{A}^{\{\ell,m\}}\tr\,  \mathbf{B}^{\{\ell\}} + \mathbf{B}^{\{m\}}\, \mathbf{A}^{\{m,\ell\}} -  \mathbf{b}^{\{m\}} \,\mathbf{b}^{\{\ell\}}\tr \in \Re^{s^{\{m\}} \times s^{\{\ell\}}}, \\
\label{eqn:P-matrix}
\mathbf{P} &\coloneqq& \big[ \mathbf{P}^{\{m,\ell\}} \big]_{1 \le \ell,m \le \nparts}  \in \Re^{s \times s},
\quad \textnormal{where} \quad s \coloneqq \sum_{m=1}^{\nparts} s^{\{m\}}.
\end{eqnarray}
\end{subequations}
The matrix $\mathbf{P} \in \Re^{s \times s}$ \eqref{eqn:P-matrix} is symmetric since  $\mathbf{P}^{\{\ell,m\}} = \mathbf{P}^{\{m,\ell\}}\tr$. It was shown in \cite{Sandu_2015_GARK} that the GARK method is algebraically stable iff the matrix $\mathbf{P}$ is non-negative definite.
Using the Butcher tableau \eqref{eqn:general-Butcher-tableau} the matrix \eqref{eqn:P-matrix} is constructed as:
\begin{equation}
\label{eqn:P-matrix-gark}
\B_{\textsc{gark}} \coloneqq \mbox{diag}\left(\b_{\textsc{gark}}\right), \quad 
\mathbf{P} = \mathbf{A}_{\textsc{gark}}\tran\,  \B_{\textsc{gark}} + \B_{\textsc{gark}}\, \mathbf{A}_{\textsc{gark}} 
-  \b_{\textsc{gark}} \,\b_{\textsc{gark}}\tran.
\end{equation}

We have the following property that generalizes the characterization of symplectic Runge Kutta schemes \cite{Hairer_book_I}.

\begin{theorem}[Symplectic GARK schemes]
\label{theorem.symplectic}
Consider a GARK scheme~\eqref{eqn:GARK-symplectic} applied to an Hamiltonian splitting~\eqref{eqn:partitioned-sympl.system},
and its matrix $\mathbf{P}$ defined by \eqref{eqn:P-matrix}.
The GARK scheme is symplectic  if and only if $\mathbf{P}=\mathbf{0}_{s \times s}$, which is equivalent to:
\begin{equation}
\label{eqn:GARK-symplectic-condition}
\begin{split}
\mathbf{P}^{\{m,\ell\}}  &= \mathbf{A}^{\{\ell,m\}}\tr\,  \B^{\{\ell\}} + \B^{\{m\}}\, \mathbf{A}^{\{m,\ell\}} -  \b^{\{m\}} \,\b^{\{\ell\}}\tr
= \Zero_{s^{\{m\}} \times s^{\{\ell\}}}, \\
& \forall\; \ell,m=1,\ldots,\nparts. 
\end{split}
\end{equation}
\end{theorem}

%\sandu{Does the proof below need the fact that each component is symplectic? Or does it work for any type of partitioning, as long as the sum of components is symplectic? That would be a major difference from other composition schemes.}

\begin{proof}
The proof is based on symplectic NB-series introduced in~\cite{armusa97}, and is similar to the proof for Runge-Kutta, partitioned Runge-Kutta 
methods~\cite{saca93}  and ARK schemes~\cite{armusa97}. 
\ifreport
An alternative proof using differential forms is given in \Cref{app:proof-of-symplecticness}.
\fi

N-trees~\cite{armusa97} are a generalization of P-trees from the case of component partitioning  to  the  general  case  of  right-hand  side  partitioning~\eqref{eqn:additive-ode}.   The  set  $\NT$  of N-trees consists of all Butcher trees with colored vertices; each vertex is assigned one of $\nparts$ different  colors  corresponding  to  the $\nparts$ components  of  the  partition.   Similar to regular Butcher trees each vertex is also assigned a label.  The order $\rho(u)$ is the number of nodes of $u \in \NT$. 

The  empty  N-tree  is  denoted  by $\emptyset$.   The  N-tree  with  a  single  vertex  of  color $m$ is denoted by $\tau_{\{m\}}$.  The N-tree $u \in \NT$ with $\rho(u)>1$ and a root of color $m$ can be represented  as $u=  [u_1,...,u_r]_{\{m\}}$,  where $\{u_1,...,u_r\}$ are  the  non-empty  subtrees (N-trees) arising from removing the root of $u$. 
The elementary differential associated with the N-tree  $u$ and evaluated at $\yy$ is:
\begin{equation}\label{eq:elem-diffs}
F(u)(\yy)  \coloneqq 
\begin{cases} \yy, & u = \emptyset, \\
 \fun^{\{m\}}_{\yy^r}\big( F(u_1)(\yy), \dots, F(u_r)(\yy) \big), & u=[u_1,\ldots,u_r]_{\{m\}}.
 \end{cases}
\end{equation}
An NB-series is a formal power expansion:
\begin{equation}
\label{eqn:NB-series-definition}
\textnormal{NB}(\mathfrak{a},\yy(t)) \coloneqq  \sum_{\tree \in \NT} \mathfrak{a}(\tree)\, \frac{\Dt^{\rho(\tree))}}{\sigma(\tree)}\, F(\tree)(\yy(t))\,,
\end{equation}
where $\mathfrak{a} : \NT \to \Re$ is a mapping that assigns a real number to each N-tree; with some abuse of nomenclature we call the mappings NB-series as well. It can be shown that the stage vectors and the solution of the GARK scheme \eqref{eqn:GARK-symplectic} can be written as NB-series:
\begin{equation}
\begin{split}\label{eq:NB-series}
 k_i^{\{m\}}  &= \textnormal{NB}(\mathfrak{g}_i^{\{m\}},\y_0), \quad
\mathfrak{g}^{\{m\}} = [\mathfrak{g}_1^{\{m\}} \dots \mathfrak{g}_{s^{\{m\}}}^{\{m\}}]\tr;\\
\y_1 &= \textnormal{NB}(\mathfrak{a},\y_0).
\end{split}
\end{equation}

The Butcher product $u \bullet v$ of the NT-trees $u,v$ is defined as follows:
\begin{equation}
\label{eqn:Butcher-product}
\begin{array}{rl}
u &= [u_1,\dots,u_r]_{\{m\}}, \\
v &= [v_1,\dots,v_p]_{\{n\}},
\end{array}
\quad 
u \bullet v \coloneqq \begin{cases}
u, & v = \emptyset, \\
[v]_{\{m\}}, & u = \tau_{\{m\}}, \\
[u_1,\dots,u_r,v]_{\{m\}}, & \textnormal{otherwise}.
\end{cases}
\end{equation}
Consider the NB-series associated with a partitioning where each component is Hamiltonian. In Araujo et al \cite{armusa97} it is shown that the NB-series $\mathfrak{a}$ is symplectic (for the special case of $\mathbf{J}$ given by~\eqref{eq-ham}) iff for each pair $u,v \in \NT \backslash \{\emptyset\}$ it holds that
\begin{align}
\label{cond.sympl.ara}
\mathfrak{a}(u \bullet v) + \mathfrak{a}(v \bullet u) &= \mathfrak{a}(u)\,\mathfrak{a}(v).
\end{align}
This result also holds for an arbitrary regular skew-symmetric matrix $\mathbf{J}$, as the argumentation in~\cite{armusa97} is based only on the skew-symmetry of $\mathbf{J}$, and not on the special structure of $\mathbf{J}$ in the Hamiltonian dynamics case~\eqref{eq-ham}.

Consider the non-empty NT-trees $u$ and $v$ in \eqref{eqn:Butcher-product} and define:
\begin{equation}
\label{eqn:U-and-V}
    U \coloneqq  \mathfrak{g}^{\{m\}} (u_1) \times \cdots \times \mathfrak{g}^{\{m\}} (u_r), \quad
    V \coloneqq  \mathfrak{g}^{\{n\}} (v_1) \times \cdots \times \mathfrak{g}^{\{n\}}(v_p),
\end{equation}
where $\times$ denotes the element-by-element product of vectors. From \eqref{eqn:GARK-symplectic}
we have the following expressions for the corresponding NB-series:
\begin{align*}
    \mathfrak{a}(u) & = \b^{\{m\}}\tr\, U = \sum_{i=1}^{s^{\{m\}}} b_i^{\{m\}} U_i, \quad
    \mathfrak{a}(v)  = \b^{\{n\}}\tr\, V = \sum_{j=1}^{s^{\{n\}}} b_j^{\{n\}} V_j, \\
    \mathfrak{a}(u \bullet v) & =
    \b^{\{m\}}\tr\, \left( U \times \A^{\{m,n\}} V \right) 
    =
    \sum_{i=1}^{s^{\{m\}}} \sum_{j=1}^{s^{\{n\}}} b_i^{\{m\}}\, U_i \, \A_{i,j}^{\{m,n\}}\, V_j, \\
    \mathfrak{a}(v \bullet u) & =
    \b^{\{n\}}\tr \left( V \times \A^{\{n,m\}} U \right) 
    =
    \sum_{j=1}^{s^{\{n\}}} \sum_{i=1}^{s^{\{m\}}} b_j^{\{n\}}\, V_j\, \A_{j,i}^{\{n,m\}}\, U_i,
\end{align*}
After reordering the coefficients,  the symplecticness condition~\eqref{cond.sympl.ara} reads:
\begin{align}
\label{eq.derivation.sympl}
  \sum_{i=1}^{s^{\{m\}}} \sum_{j=1}^{s^{\{n\}}} \left( \b_i^{\{m\}} \, \A_{i,j}^{\{m,n\}} + 
  (\A^{\{n,m\}}\tr)_{i,j} \, \b_j^{\{n\}} -  \b_i^{\{m\}}\,\b_j^{\{n\}} \right)\,  U_i\, V_j & = 0,
\end{align}
which is equivalent to $\mathbf{P}^{\{m,n\}} = \mathbf{0}$. 
\end{proof}

\begin{corollary}A GARK scheme \eqref{eqn:GARK2} based on a general splitting~\eqref{eqn:partitioned-sympl.system-without-structure} is symplectic iff~\eqref{eqn:GARK-symplectic-condition} and
\begin{equation}
    \label{add.sypl.comd}
    \b^{\{\mu\}}=\b^{\{\sigma\}}=\b, \quad s^{\{\mu\}}=s^{\{\sigma\}}, \quad \forall \, \mu,\sigma=1,\ldots,N
\end{equation}
hold.
\end{corollary}

\begin{proof} 
This follows directly from Araujo et al \cite{armusa97}:   for a general splitting, the NB-series $\mathfrak{a}$ is symplectic iff in addition to~\eqref{cond.sympl.ara}  the  following condition holds: 
$\mathfrak{a}(u)=\mathfrak{a}(v)$ for each pair of nonempty 
N-trees $u$, $v$ that differ only in the color of their roots.
\end{proof}

\begin{remark}
\label{rem:single-RK}
\begin{leaveout}
As for additive Runge-Kutta schemes, symplectic GARK schemes applied to a general splitting can always be rewritten as a single symplectic Runge-Kutta method for the non-decomposed system (see also Example~\ref{example-GARK-IMIM}). \color{red}{MG: Is this really true?} 
\sandu{This seems to contradict our Example \ref{example-GARK-IMIM}?
From \eqref{eqn:GARK-symplectic-condition} we have:
\begin{equation*}
\begin{split}
\mathbf{A}^{\{\ell,m\}}\tr\,  \B + \B\, \mathbf{A}^{\{m,\ell\}} &= \b \,\b\tr, \quad \forall\; \ell,m=1,\ldots,\nparts. \\
a_{j,i}^{\{\ell,m\}}\,\b_j + \b_i\, a_{i,j}^{\{m,\ell\}}&= \b_i \,\b_j.
\end{split}
\end{equation*}
Does this really imply that $\mathbf{A}^{\{m,\ell\}} =\A$ $\forall\; \ell,m$? Perhaps not in general. %
Consider decoupled schemes with $b_i \ne 0$.
For $\ell=m$ we have that 
$\mathbf{A}^{\{m,m\}}=\A$ equal to the column format of implicit symplectic schemes:
\[
\A = \textnormal{tril}(\b\,\one^T,-1) + \frac{1}{2}\,\textnormal{diag}(\b),
\]
and
\[
\begin{split}
& a_{j,i}^{\{\ell,m\}}\,a_{i,j}^{\{m,\ell\}} =0, \quad \ell \ne m \quad \Rightarrow \quad \\
& a_{j,i}^{\{\ell,m\}} = 0 ~~\&~~ a_{i,j}^{\{m,\ell\}}= \b_j \quad \textnormal{or} \quad
a_{i,j}^{\{m,\ell\}} = 0 ~~\&~~ a_{j,i}^{\{\ell,m\}} = \b_i. 
\end{split}
\]
This can be a composition of steps: 
\[
\begin{array}{ccc | ccc}
\frac{b_1}{2} & 0 & 0 & 0 & 0 & 0 \\
b_1 & \frac{b_2}{2} & 0& 0 & 0 & 0 \\
b_1 & b_2 & \frac{b_3}{2}& 0 & 0 & 0 \\
\hline
b_1 & b_2 & b_3 & \frac{b_1}{2} & 0 & 0 \\
b_1  & b_2 & b_3 & b_1 & \frac{b_2}{2} & 0 \\
b_1  & b_2 & b_3 & b_1 & b_2 & \frac{b_3}{2}
\end{array}
\]
Another legitimate structure that is decoupled, but is not a traditional composition since the stages are computed in mixed order:
\[
\begin{array}{ccc | ccc}
\frac{b_1}{2} & 0 & 0 & 0 & 0 & 0 \\
b_1 & \frac{b_2}{2} & 0& b_1 & 0 & 0 \\
b_1 & b_2 & \frac{b_3}{2}& b_1 & b_2 & 0 \\
\hline
b_1 & 0  & 0 & \frac{b_1}{2} & 0 & 0 \\
b_1 & b_2 & 0 & b_1 & \frac{b_2}{2} & 0 \\
b_1 & b_2 & b_3 & b_1 & b_2 & \frac{b_3}{2}
\end{array}
\]
%
%How about the following coupled, second order Gauss-based method?
%%
%\[
%\begin{array}{c| cc | cc}
%\frac{1}{2} - \frac {\sqrt{3}}{6} & \frac{1}{4}  & \frac{1}{4} - \frac {\sqrt{3}}{6} & a_{1,1} & a_{1,2}  \\
%\frac{1}{2} + \frac {\sqrt{3}}{6} & \frac{1}{4} + \frac {\sqrt{3}}{6}  &  \frac{1}{4}  & a_{2,1} & a_{2,2} \\
%\hline
% & \frac{1}{2}-a_{1,1}  &  \frac{1}{2}-a_{2,1} & \frac{1}{4}  & \frac{1}{4} - \frac {\sqrt{3}}{6}  \\
% & \frac{1}{2}-a_{1,2} & \frac{1}{2}-a_{2,2} & \frac{1}{4} + \frac {\sqrt{3}}{6}  &  \frac{1}{4}  \\
%\hline
%& \frac{1}{2} &  \frac{1}{2} &  \frac{1}{2} &  \frac{1}{2} 
%\end{array}
%\]
%%
%or
%
%\[
%\begin{array}{cc | cc}
%\frac {b_1}{2} & 0 & a_{1,1} & a_{1,2}  \\
%b_1 &\frac {b_2}{2}  &  a_{2,1} & a_{2,2} \\
%\hline
%b_1  - a_{1,1}  &  b_2  - \frac{b_2}{b_1}\,a_{2,1} & \frac {b_1}{2} & 0   \\
% b_1  - \frac{b_1}{b_2}\,a_{1,2} &  b_2  - a_{2,2} & b_1 &\frac {b_2}{2}   \\
%\hline
%b_1 &  b_2 &  b_1 &  b_2 
%\end{array}
%\]
% a_{j,i}^{\{\ell,m\}} = \b_i  - \frac{b_i}{b_j}\, a_{i,j}^{\{m,\ell\}} .
}

\guenther{In the first example $\A^{\{1,2\}}=0$, and thus the order limited by one. The same for the second example: here is a contradiction in the order 2 condition. Is this then a good example? }
\sandu{Remark \ref{rem:single-RK} seems to contradict our Example~\ref{example-GARK-IMIM}. We show there that the method cannot be formulated as a single scheme applied to the non-partitioned system. We can refer to that and delete the remark?}
\end{leaveout}
Additive Runge-Kutta schemes applied to a general splitting can always be rewritten as a single symplectic Runge-Kutta method for the non-decomposed system, see~\cite{armusa97}. However this is not the case for symplectic GARK schemes, as shown in Example~\ref{example-GARK-IMIM}.
\end{remark}

\subsection{Order conditions}

As shown in~\cite{Sandu_2015_GARK}, 
the order conditions for a GARK method
\eqref{eqn:GARK} are obtained from the order conditions of
ordinary Runge--Kutta methods.  The usual labeling of the Runge-Kutta
coefficients (subscripts $i,j,k,\ldots$) is accompanied by a
corresponding labeling of the different partitions
(superscripts $m,s,t,\ldots$).
%Let $\one^{{\{s\}}}$ be a vector of ones of dimension
%$s^{{\{s\}}}$, and
%$\mathbf{c}^{\{m,s\}} \coloneqq \mathbf{A}^{\{m,s\}} \cdot
%\one^{{\{s\}}}$. 
The conditions for orders one to four are as follows: 
\begin{subequations} 
\label{eqn:GARK-order-conditions-1-to-4}
\begin{align}
\label{eqn:GARK-order-conditions-1}
\mathbf{b}^{\{m\}}\tr \cdot \one^{{\{m\}}} = 1,  & \quad \forall\; m, &  (\textnormal{order}~ 1) \\
%\end{eqnarray}
%%
%\begin{eqnarray}
\label{eqn:GARK-order-conditions-2}
 \mathbf{b}^{\{m\}}\tr \cdot \mathbf{c}^{\{m,\ell\}}  =  \sfrac{1}{2}, &\quad  \forall\; m, \ell,  &  (\textnormal{order}~ 2) \\
%\end{eqnarray}
%%
%\begin{eqnarray}
\label{eqn:GARK-order-conditions-3a}
\mathbf{b}^{\{m\}}\tr \cdot \left(  \mathbf{c}^{\{m,\ell\}} \times \mathbf{c}^{\{m,s\}} \right)
 =  \sfrac{1}{3}, & \quad \forall\; m,\; \forall\; \ell\leq s, & (\textnormal{order}~ 3) \\
\label{eqn:GARK-order-conditions-3b}
\mathbf{b}^{\{m\}}\tr \cdot \mathbf{A}^{\{m,\ell\}}  \cdot \mathbf{c}^{\{\ell, s\}} 
 =  \sfrac{1}{6}, & \quad \forall\; m,\ell,s,  & (\textnormal{order}~ 3) \\
%\end{eqnarray}
%%
%\begin{eqnarray}
\label{eqn:GARK-order-conditions-4a}
\mathbf{b}^{\{m\}}\tr \cdot \left( \mathbf{c}^{\{m,\ell\}} \times \mathbf{c}^{\{m,s\}} \times \mathbf{c}^{\{m,t\}} \right)
 =  \sfrac{1}{4}, & \quad \forall\; m,\; \forall\; \ell\leq s\leq t, & (\textnormal{order}~ 4) \\
\label{eqn:GARK-order-conditions-4b}
\mathbf{b}^{\{m\}}\tr \times \left(\mathbf{c}^{\{m,\ell\}} \times \mathbf{A}^{\{m,s\}}  \cdot \mathbf{c}^{\{s,t\}} \right) 
 =  \sfrac{1}{8}, & \quad \forall\; m,\ell,s,t, & (\textnormal{order}~ 4) \\
\label{eqn:GARK-order-conditions-4c}
\mathbf{b}^{\{m\}}\tr \cdot \mathbf{A}^{\{m,\ell\}}  \cdot
\left( \mathbf{c}^{\{\ell,s\}} \times \mathbf{c}^{\{\ell,t\}} \right)
 =  \sfrac{1}{12},  &\quad \forall\; m,\ell,\; \forall\;s \leq t, &  (\textnormal{order}~ 4) \\
\label{eqn:GARK-order-conditions-4d}
\mathbf{b}^{\{m\}}\tr \cdot  \mathbf{A}^{\{m,\ell\}}  \cdot \mathbf{A}^{\{\ell,s\}}  \cdot \mathbf{c}^{\{s,t\}}
 =  \sfrac{1}{24}, & \quad \forall\; m,\ell,s,t. & (\textnormal{order}~ 4)
\end{align}
\end{subequations}
Here, the standard matrix and vector
multiplication is denoted by dot (e.g.,
$\mathbf{b}\tr \cdot \mathbf{c}$ is a dot product), whereas the cross
denotes component-wise multiplication (e.g.,
$\mathbf{b}\times \mathbf{c}$ is a vector of element-wise products).
For internally consistent schemes these order conditions simplify considerably. Moreover,
for symplectic GARK schemes many of the order conditions \eqref{eqn:GARK-order-conditions-1-to-4} are redundant.

\begin{remark}[Redundancy of order conditions for symplectic GARK schemes]
\label{rem:order-redundancy}
Assume that the symplectic GARK method has a solution with NB-series coefficients $\mathfrak{a}$, and that it satisfies all conditions up to order $k$:
\begin{equation*}
\mathfrak{a}(u) =  \sfrac{1}{\gamma(u)}  \quad \forall u:~\rho(u) \le k.
\end{equation*}
Symplecticness equation~\eqref{cond.sympl.ara}  implies 
\begin{align}
\label{rel.oc}
\mathfrak{a}(u \bullet v) + \mathfrak{a}(v \bullet u) &= \sfrac{1}{\gamma(u)\,\gamma(v)},  \quad \forall\; u,v:~\rho(u),\rho(v) \le k.
\end{align}
Since $\rho(u\bullet v) = \rho(v\bullet u) = k + 1$, equation \eqref{rel.oc} involves two order $k+1$ conditions; if one is satisfied, then the other is satisfied as well. Specifically, assuming that $\mathfrak{a}(u \bullet v) = 1/\gamma(u \bullet v)$ we have
\begin{equation*}
\begin{split}
& \mathfrak{a}(u \bullet v) = \sfrac{1}{\gamma(u \bullet v)} = \sfrac{\rho(u)}{\rho(u)+\rho(v)}\,\sfrac{1}{\gamma(u)\,\gamma(v)} 
\quad \Rightarrow \\
& \mathfrak{a}(v \bullet u) = \sfrac{1}{\gamma(u)\,\gamma(v)}-\mathfrak{a}(u \bullet v)
= \sfrac{\rho(v)}{\rho(u)+\rho(v)}\,\sfrac{1}{\gamma(u)\,\gamma(v)} = \sfrac{1}{\gamma(v \bullet u)}.
\end{split}
\end{equation*}
For $v=\tau_{\{\ell\}}$ and $u=[u_1,\dots,u_r]_{\{m\}}$, $\rho(u) = k < p$, we have
\begin{equation*}
u \bullet v = [u_1,\dots,u_r,\tau_{\{\ell\}}]_{\{m\}}, \quad 
v \bullet u = [[u_1,\dots,u_r]_{\{m\}}]_{\{\ell\}}, 
\end{equation*}
with $\rho(u\bullet v) = \rho(v\bullet u) = k + 1 \le p$. Equation~\eqref{rel.oc} yields:
\begin{align*}
    \mathfrak{a}([u_1,\dots,u_r,\tau_{\{\ell\}}]_{\{m\}})+\mathfrak{a}([[u_1,\dots,u_r]_{\{m\}}]_{\{\ell\}}) = \mathfrak{a}(\tau_{\{\ell\}}) \,\mathfrak{a}([u_1,\dots,u_r]_{\{m\}}),
\end{align*}
which implies the order $k+1$  relation
\begin{align}
\label{eq.symp.cond.bs}
\b^{\{m\}}\tr\, (U \times \c^{\{m,\ell\}}) +  
\b^{\{\ell\}}\tr\,\A^{\{\ell,m\}}\,U = \frac{1}{\gamma([u_1,\dots,u_r]_{\{m\}})},
\end{align}
where $U$ is defined in \eqref{eqn:U-and-V}.
\end{remark}

This redundancy of order conditions discussed in Remark  \ref{rem:order-redundancy} yields the following reduction.

\begin{theorem}[Reduced order conditions for symplectic GARK schemes]
\label{th.sympl.order.cond}
For symplectic GARK schemes the number of order conditions is reduced due to redundancy:
\begin{itemize}
    \item order 2: the number of $\nparts^2$ order conditions reduces to only $\nparts(\nparts-1)/2$ order conditions;
    \item order 3: the number of $(3\nparts^3 + \nparts^2)/2$ order conditions reduces to only $(\nparts^3 + \nparts^2)/2$ order conditions, if the scheme is of at least order 2;
    \item order 4: the number of $(8\nparts^4 + 3\nparts^3 + \nparts^2)/3$ order conditions reduces to only $(4 \nparts^4 + 3\nparts^3 - \nparts^2)/6$ order conditions, if the scheme is of  at least order 3.
\end{itemize}
\end{theorem}
\begin{proof}
We assume that the symplectic GARK scheme has at least order one.

\paragraph{Order two} Using the redundancy relation~\eqref{eq.symp.cond.bs} with $u=\tau_{\{m\}}$, we get 
\begin{align}
    \label{eq.symplgarko2}
       \b^{\{m\}}\tr\, \c^{\{m,\ell\}} +   \b^{\{\ell\}}\tr\,\c^{\{\ell,m\}} -1 =0,
    \end{align}
which yields for $\ell=m$ the order two conditions
\begin{align*}
    \b^{\{m\} }\tr \c^{\{m,m\}} & = \sfrac{1}{2}.
\end{align*}
Assuming that the order two condition
\begin{align*}
    \b^{\{m\} }\tr \c^{\{m,\ell\}} & = \sfrac{1}{2}
\end{align*}
holds for $\ell<m$, \eqref{eq.symplgarko2} yields the order two condition for $\ell>m$.  This condition is automatically fulfilled for internally consistent schemes.

\paragraph{\it Order three} Using \eqref{eq.symp.cond.bs} with $v=\tau_{\{\ell\}}$ and $u=[\tau_{\{s\}}]_{\{m\}}$, we get $U = \c^{\{m,s\}}$, and:
\begin{align*}
\b^{\{m\}}\tr\, (\c^{\{m,s\}} \times \c^{\{m,\ell\}}) +  
\b^{\{\ell\}}\tr\,\A^{\{\ell,m\}}\,\c^{\{m,s\}} = \sfrac{1}{2}.
\end{align*}
Thus the order three condition \eqref{eqn:GARK-order-conditions-3a} (for a set of partitions $m,s,\ell$) yields the corresponding order condition  \eqref{eqn:GARK-order-conditions-3b}, and vice versa:
\begin{align}
\label{eq.symplgarko3}
    \b^{\{\ell\}}\tr (\c^{\{\ell,m\}} \times \c^{\{\ell,s\}}) & =  \sfrac{1}{3}
    \quad \Leftrightarrow \quad 
    \b^{\{m\}}\tr \A^{\{m,\ell\}} \c^{\{\ell,s\}}  = \sfrac{1}{6}.
\end{align}
Since \eqref{eqn:GARK-order-conditions-3a} consists of $(\nparts^3 + \nparts^2)/2$ order conditions, the total number of order three conditions becomes $(\nparts^3 + \nparts^2)/2$ for symplectic GARK schemes.

{\it Order four.}
Using the redundancy relation~\eqref{eq.symp.cond.bs} with $u=[[\tau_{\{\ell\}}]_{\{s\}}]_{\{m\}}$, we get $U = \A^{\{m,s\}}\c^{\{m,s\}}$, and the following relation:
\begin{align}
\label{eq.symplgarko4.1}
\b^{\{m\}}\tr\, (  \A^{\{m,s\}}\,\c^{\{s,\ell\}} \times \c^{\{m,t\}}) +  
\b^{\{t\}}\tr\,\A^{\{t,m\}}\, \A^{\{m,s\}}\,\c^{\{s,\ell\}} - \sfrac{1}{6} = 0.
\end{align}
If an order condition \eqref{eqn:GARK-order-conditions-4b} is satisfied (for a set of partitions $m,s,\ell,t$) then so is the corresponding order condition \eqref{eqn:GARK-order-conditions-4d}, and vice-versa.

Using the redundancy relation~\eqref{eq.symp.cond.bs} with $u=[\tau_{\{s\}},\tau_{\{t\}}]_{\{\ell\}}$, we get $U = \c^{\{\ell,s\}} \times \c^{\{\ell,t\}}$, and the following relation:
\begin{align}
\label{eq.symplgarko4.2}
\b^{\{\ell\}}\tr\, (  \c^{\{\ell,s\}} \times \c^{\{\ell,t\}} \times \c^{\{\ell,m\}}) 
+  
\b^{\{m\}}\tr\,\A^{\{m,\ell\}}\,
(\c^{\{\ell,s\}} \times \c^{\{\ell,t\}}) - \sfrac{1}{3} =0.
\end{align}
If an order condition \eqref{eqn:GARK-order-conditions-4a} is satisfied then so is the corresponding order condition \eqref{eqn:GARK-order-conditions-4c}, and vice-versa. \eqref{eqn:GARK-order-conditions-4a} consists of $(\nparts^4 + 3\nparts^3 + 2\nparts^2)/6$ order conditions, while the equivalent order condition \eqref{eqn:GARK-order-conditions-4c} consists of $(\nparts^4 + \nparts^3)/2$ equations. Since $(\nparts^4 + 3\nparts^3 + 2\nparts^2)/6 \leq (\nparts^4 + \nparts^3)/2 \; \forall\; \nparts \in \mathbb{N}$, \eqref{eqn:GARK-order-conditions-4a} and \eqref{eqn:GARK-order-conditions-4c} reduce to $(\nparts^4 + 3\nparts^3 + 2\nparts^2)/6$ conditions in case of symplectic GARK schemes.

For $v=[\tau_{\{s\}}]_{\{\ell\}}$ and  $u=[\tau_{\{t\}}]_{\{m\}}$ we have
\begin{align*}
    u \bullet v & =[\tau_{\{t\}},[\tau_{\{s\}}]_{\{\ell\}}]_{\{m\}}, &
    v \bullet u & = [\tau_{\{s\}},[\tau_{\{t\}}]_{\{m\}}]_{\{\ell\}},
\end{align*}
and thus~\eqref{rel.oc} leads to
\begin{align*}
    \mathfrak{a}([\tau_{\{t\}},[\tau_{\{s\}}]_{\{\ell\}}]_{\{m\}})+\mathfrak{a}([\tau_{\{s\}},[\tau_{\{t\}}]_{\{m\}}]_{\{\ell\}}) = \mathfrak{a}([\tau_{\{t\}}]_{\{m\}}) \,
     \mathfrak{a}([\tau_{\{s\}}]_{\{\ell\}}),
\end{align*}
which implies the following relation between order four conditions \eqref{eqn:GARK-order-conditions-4b}:
\begin{equation}
\label{eq.symp.cond.bs2}
\begin{split}
\b^{\{m\}}\tr\, (\c^{\{m,t\}} \times \A^{\{m,\ell\}}\,\c^{\{\ell,s\}}) & +   
\b^{\{\ell\}}\tr\,(\c^{\{\ell,s\}} \times \A^{\{\ell,m\}}\,\c^{\{m,t\}})  \\
& = \frac{1}{\gamma([\tau_{\{t\}}]_{\{m\}})}
\frac{1}{\gamma([\tau_{\{s\}}]_{\{\ell\}})}=\sfrac{1}{4}.
\end{split}
\end{equation}
Setting $t=\ell$ and $m=s$, the second redundancy relation~\eqref{eq.symp.cond.bs2} yields the $\nparts^2$ order four conditions
\begin{align*}
\b^{\{\ell\}}\tr\,(\c^{\{\ell,t\}} \times \A^{\{\ell,\ell\}}\,\c^{\{\ell,t\}}) - \sfrac{1}{8}= 0
%    (\b^{\{\ell\}} \times\c^{\{\ell,t\}})\tr  \A^{\{\ell,\ell\}}  \c^{\{\ell,m\}}) - \frac{1}{8}= 0
    \end{align*}
as part of~\eqref{eqn:GARK-order-conditions-4b}. 
If in addition 
\begin{align*}
    % \label{eq.symplgarko4.3}
  \b^{\{\ell\}}\tr\,(\c^{\{\ell,m\}} \times \A^{\{\ell,t\}}\,\c^{\{t,s\}}) - \sfrac{1}{8}= 0  
%    
%    (\b^{\{\ell\}} \times \c^{\{\ell,m\}})\tr  \A^{\{\ell,t\}}  \c^{\{t,s\}} - \frac{1}{8}  = 0 
\end{align*}
holds for $\ell+m\le t+s$, then the overall 
$\nparts^2(\nparts^2-1)/2$ conditions are equivalent to~\eqref{eqn:GARK-order-conditions-4b}.
Assuming now that these $\nparts^2(\nparts^2-1)/2$ conditions hold, \eqref{eq.symplgarko4.1} yields the $\nparts^4$ 
order conditions~\eqref{eqn:GARK-order-conditions-4d}.
\end{proof}

\begin{remark}
Theorem~\eqref{th.sympl.order.cond} contains all reductions implied by symplecticness, as follows:
\begin{itemize}
    \item for order two, there is only one symplecticness condition:   $a(u \bullet v)$ defines an order two condition, if  both $u$ and $v$ contain one node.
    \item for order three, there is only one symplecticness condition:   $a(u \bullet v)$ defines an order three condition, if  $u$ and $v$ contain one and node and two nodes, resp.
    \item for order four, there are only three symplecticness conditions:   $a(u \bullet v)$ defines an order four condition, if  $u$ and $v$ contain one and three nodes and  two and two nodes, resp., which gives $2+1=3$ conditions. 
\end{itemize}
Overall, symplecticness yields 5 conditions reducing the number of order conditions, which have all been discussed in Theorem~\ref{th.sympl.order.cond}. 
\end{remark}

\begin{corollary}[Reduced number of order conditions for internally consistent symplectic GARK schemes]
\label{th.sympl.order.cond.intcons}
If the symplectic GARK scheme is internally consistent, then
\begin{itemize}
    \item order 2: the order two conditions~\eqref{eqn:GARK-order-conditions-2} are automatically fulfilled;
    \item order 3: if the scheme has at least order two, only the order conditions~\eqref{eqn:GARK-order-conditions-3a} ($\nparts$  equations)
%    \begin{align*}
 %   \b^{\{\ell\}}\tr (\c^{\{\ell\}} \times %\c^{\{\ell\}}) & =  \sfrac{1}{3}
%\end{align*}
have to be fulfilled; 
    \item order 4: if the scheme has at least order three, then only the order four conditions \eqref{eqn:GARK-order-conditions-4a} ($\nparts$ equations) and ~\eqref{eqn:GARK-order-conditions-4b} for $m<s$ ($\nparts(\nparts-1)/2$ equations) have to be fulfilled.
\end{itemize}
\end{corollary}

\begin{proof}
The proposition for internally consistent schemes follows directly from the results above in theorem~\ref{th.sympl.order.cond}. 

{\it Order two.\/} From \eqref{eq.symplgarko2} we get the $\nparts$ order two conditions~\eqref{eqn:GARK-order-conditions-2}.
%\begin{align*}
%    \b^{\{m\} }\tr \c^{\{m\}} & = \frac{1}{2}.
%\end{align*}

{\it Order three.\/} 
Here the order conditions~\eqref{eq.symplgarko3}
%\begin{align*}
 %   \b^{\{\ell\}}\tr (\c^{\{\ell\}} \times \c^{\{\ell\}}) & =  \sfrac{1}{3}.
%\end{align*}
reduce to the $\nparts$ order conditions~\eqref{eqn:GARK-order-conditions-3a}.

{\it Order four.\/} 
If the $\nparts(\nparts-1)/2$ order conditions
\begin{align*}
    \b^{\{\ell\}}\tr\, \left( \c^{\{\ell\}} \times \A^{\{\ell,t\}}  \c^{\{t\}} \right) & = \sfrac{1}{8}
\end{align*}
for $\ell<t$ hold, then the $\nparts^2$ order conditions
~\eqref{eqn:GARK-order-conditions-4b} are fulfilled, and as before the $\nparts^3$ order conditions~\eqref{eqn:GARK-order-conditions-4d}.
Assuming now that the $\nparts$ order conditions of~\eqref{eqn:GARK-order-conditions-4a} are fulfilled, the $\nparts^2$ order conditions~\eqref{eqn:GARK-order-conditions-4c} 
hold. 
\end{proof}

%%%%%%%%%%%%%%%%%%%%%%%%%%%
\subsection{Symmetry and time-reversibility}
%%%%%%%%%%%%%%%%%%%%%%%%%%%

\par
\begin{remark}[Time-reversed GARK method]
Let $\mathcal{P}^{\{m\}} \in \Re^{s^{\{m\}} \times s^{\{m\}}}$  be the permutation matrix that reverses the order of the entries of a vector. In matrix notation the time-reversed GARK method is:
\begin{subequations}
\label{eqn:GARK-time-reversed-1}
\begin{eqnarray}
\label{eqn:GARK-time-reversed-1.a}
\underline{\b}^{\{m\}} & \coloneqq & \mathcal{P}^{\{m\}}\,\b^{\{m\}} \quad \Leftrightarrow \quad\underline{b}_j^{\{m\}}  =  b_{s^{\{m\}}+1-j}^{\{m\}}, ~~ \forall\, j; \\
%\quad \underline{\hat{\b}}^{\{m\}}  \coloneqq \mathcal{P}\,\hat{\b}^{\{m\}}, \\
\label{eqn:GARK-time-reversed-1.b}
\underline{\A}^{\{\ell,m\}} & \coloneqq & \one^{\{\ell\}}\,\b^{\{m\}}\tr - \mathcal{P}^{\{\ell\}}\,\A^{\{\ell,m\}}\,\mathcal{P}^{\{m\}} \\
\nonumber
&& \Leftrightarrow \quad
\underline{a}_{i,j}^{\{\ell,m\}}  =   b_{j}^{\{m\}} - a_{s^{\{\ell\}}+1-i,s^{\{m\}}+1-j}^{\{\ell,m\}}, ~~ \forall\, i,j.
%i=1,\ldots,s^{\{m\}}, \, j=1,\ldots,s^{\{q\}}.
 %\\ %\quad
%\hat{\underline{\A}}^{\{\ell,m\}} \coloneqq \one\,\hat{\b}^{\{m\}}\tr - \mathcal{P}\,\hat{\A}^{\{\ell,m\}}\,\mathcal{P}. 
\end{eqnarray}
The general Butcher tableau \eqref{eqn:general-Butcher-tableau} of the time-reversed GARK method is:
\begin{equation}
\label{eqn:GARK-time-reversed-tableau}
\begin{split}
& \underline{\b}_{\textsc{gark}}   = \mathcal{P}\,\b_{\textsc{gark}} , \quad
\underline{\A}_{\textsc{gark}} = \one_{s \times 1}\,\b_{\textsc{gark}}\tran - \mathcal{P}\,\A_{\textsc{gark}}\,\mathcal{P}, \\
& \textnormal{where} \quad  \mathcal{P} \coloneqq \underset{m = 1,\dots,\nparts}{\textnormal{blkdiag}} \big( \mathcal{P}^{\{m\}} \bigr) \in \Re^{s \times s}.
\end{split}
\end{equation}
\end{subequations}
\end{remark}

\begin{definition}[Symmetric GARK schemes]
The GARK scheme  \eqref{eqn:GARK2} is symmetric if it is invariant with respect to time reversion \eqref{eqn:GARK-time-reversed-1}:
\begin{subequations}
\label{eqn:GARK-symmetry}
\begin{align}
\label{eqn:GARK-symmetry.a}
\b^{\{m\}} & =  \underline{\b}^{\{m\}} & \Leftrightarrow \quad&
b_j^{\{m\}}  =  b_{s^{\{m\}}+1-j}^{\{m\}}, ~~ \forall j; \\
\label{eqn:GARK-symmetry.b}
\A^{\{\ell,m\}} & =  \underline{\A}^{\{\ell,m\}} & \Leftrightarrow \quad &
a_{i,j}^{\{\ell,m\}}  =   b_{j}^{\{m\}} - a_{s^{\{\ell\}}+1-i,s^{\{m\}}+1-j}^{\{\ell,m\}},  ~~ \forall i,j.
\end{align}
\end{subequations}
Using~\eqref{eqn:GARK-time-reversed-tableau}, the symmetry condition \eqref{eqn:GARK-symmetry} can be written compactly as
%\begin{equation}
\begin{align}
\label{cond.symm.comp}
{\b}_{\textsc{gark}}   & = \mathcal{P}\,\b_{\textsc{gark}}, \qquad
{\A}_{\textsc{gark}} = \one_{s \times 1}\,\b_{\textsc{gark}}\tran - \mathcal{P}\,\A_{\textsc{gark}}\,\mathcal{P}.
\end{align}
%\end{equation}
\end{definition}
Note that~\eqref{eqn:GARK-symmetry.b} implies \eqref{eqn:GARK-symmetry.a}.  From \eqref{eqn:GARK-time-reversed-1.b} and \eqref{eqn:GARK-symmetry.b}
\begin{equation*}
\mathcal{P}^{\{\ell\}}\,\A^{\{\ell,m\}}\,\mathcal{P}^{\{m\}} + \A^{\{\ell,m\}} = \mathcal{P}^{\{\ell\}}\,\A^{\{\ell,m\}}\,\mathcal{P}^{\{m\}} + \underline{\A}^{\{\ell,m\}}  = \one^{\{\ell\}}\,\b^{\{m\}}\tr,
\end{equation*}
and multiplying this equation from left with $\mathcal{P}^{\{\ell\}}$ and from the right with $\mathcal{P}^{\{m\}}$ yields
\begin{eqnarray*}
%\A^{\{\ell,m\}} + \mathcal{P}\,\A^{\{\ell,m\}}\,\mathcal{P} & = & \one\,\b^{\{m\}}\tr, \Rightarrow \\
\mathcal{P}^{\{\ell\}}\,\A^{\{\ell,m\}}\,\mathcal{P}^{\{m\}} + \A^{\{\ell,m\}}   & = & \one^{\{\ell\}}\,\b^{\{m\}}\tr\,\mathcal{P}^{\{m\}},
\end{eqnarray*}
which implies $\mathcal{P}^{\{m\}}\,\b^{\{m\}} = \b^{\{m\}}$.

\begin{remark}[Symplecticness of time-reversed GARK methods]
\label{remark.sympl.tr}
Using \eqref{eqn:P-matrix-gark} and \eqref{eqn:GARK-time-reversed-tableau} we have:
\begin{equation}
\begin{split}
\underline{\mathbf{P}} &= (\underline{\A}_{\textsc{gark}})\tr\,  \underline{\B}_{\textsc{gark}} + \underline{\B}_{\textsc{gark}}\, \underline{\A}_{\textsc{gark}}  -  \underline{\b}_{\textsc{gark}} \,\underline{\b}_{\textsc{gark}}\tran \\
&=  \b_{\textsc{gark}}\, \underline{\b}_{\textsc{gark}}\tran  +  \underline{\mathbf{b}}_{\textsc{gark}} \,\underline{\b}_{\textsc{gark}}\tran - 2\,\underline{\mathbf{b}}_{\textsc{gark}} \,\underline{\b}_{\textsc{gark}}\tran -\mathcal{P} \, \mathbf{P}\,\mathcal{P}.
\end{split}
\end{equation}
Consider a symplectic GARK method~\eqref{eqn:GARK-symplectic-condition} with $ \mathbf{P} = \Zero$.  The symplecticness condition $\underline{\mathbf{P}} = \Zero$ for the time-reversed scheme \eqref{eqn:GARK-time-reversed-1} reads: 
\begin{equation*}
\b_{\textsc{gark}}\, \underline{\b}_{\textsc{gark}}\tran  +  \underline{\mathbf{b}}_{\textsc{gark}} \,\underline{\b}_{\textsc{gark}}\tran = 2\,\underline{\mathbf{b}}_{\textsc{gark}} \,\underline{\b}_{\textsc{gark}}\tran.
\end{equation*}
Multiply this equation from the right by a vector of ones $\one_{s \times 1}$, and divide both sides by $\nparts$. We conclude that a {%\color{blue} 
necessary and sufficient} condition for symplecticness of the time-reversed scheme is that all weight vectors are palindomic:
\begin{equation*}
\b_{\textsc{gark}} = \underline{\b}_{\textsc{gark}} \quad \Leftrightarrow \quad 
\b^{\{m\}} = \mathcal{P}^{\{m\}}\,\b^{\{m\}} = \underline{\b}^{\{m\}}, ~~ m = 1,\dots,\nparts.
\end{equation*}
\end{remark}

With the help of time-reversed symplectic GARK methods one can derive symmetric and symplectic GARK methods:
\begin{theorem}
\label{theorem.construction.symm.sympl.}
Consider a GARK scheme $\left(\mathbf{A}_{\textsc{gark}},\b_{\textsc{gark}}\right)$ that is symplectic and has palindromic weights, $\b^{\{m\}}  ={\mathcal P}^{\{m\}} \b^{\{m\}}=\underline{\b}^{\{m\}}$ for all $m$. The GARK scheme defined by applying one step with the GARK scheme, followed by one step with its time-reversed GARK scheme, is defined by the Butcher tableau \eqref{eqn:general-Butcher-tableau}
\begin{align}
\label{eq.augmented}
\renewcommand{\arraystretch}{1.3}
    \begin{array}{cc}
        \mathbf{A}_{\textsc{gark}} & \Zero_{s \times s} \\ \one_{s \times 1}\, \underline{\b}_{\textsc{gark}}\tran & \underline{\mathbf{A}}_{\textsc{gark}} \\
 \Xhline{2\arrayrulewidth}
        \b_{\textsc{gark}}\tran & \underline{\b}_{\textsc{gark}}\tran
    \end{array}
%    \begin{pmatrix}
%        \hat \A^{\{\ell,m\}} & 0 \\ \one \b^{\{m\}}\tr & \underline{\hat \A}^{\{\ell,m\}}
%    \end{pmatrix},
%    \begin{pmatrix}
%        \hat \b^{\{m\}} \\ {\hat \b}^{\{m\}} 
%    \end{pmatrix}
\end{align}
and is both symmetric and symplectic.
\end{theorem}

\begin{proof}
According to remark~\ref{remark.sympl.tr}, the time-reversed scheme is symplectic, too, and so is~\eqref{eq.augmented} as composition of two symplectic schemes. Symmetry is given by the fact that we have a composition of a scheme with its time-reversed scheme.
\end{proof}

\begin{remark}
Consider a GARK method (with possibly some weights equal to zero). Multiplying the symmetry equation \eqref{eqn:GARK-time-reversed-tableau} by $\b_{\textsc{gark}}\tran$ from the left leads to:
\begin{equation*}
\begin{split}
%\mathcal{P}\,\A_{\textsc{gark}}\,\mathcal{P}  &= \one_{s \times 1}\,\b_{\textsc{gark}}\tran - \mathbf{A}_{\textsc{gark}}, \\
\b_{\textsc{gark}}\tran\,\mathcal{P}\,\A_{\textsc{gark}}\,\mathcal{P}  &= \b_{\textsc{gark}}\,\b_{\textsc{gark}}\tran - \b_{\textsc{gark}}\tran\, \mathbf{A}_{\textsc{gark}}.
\end{split}
\end{equation*}
and from the symplecticness equation \eqref{eqn:P-matrix-gark} we have:
\begin{equation*}
\begin{split}
\b_{\textsc{gark}} \,\underline{\b}_{\textsc{gark}}\tran -  \b_{\textsc{gark}}\, \mathbf{A}_{\textsc{gark}}  = \mathbf{A}_{\textsc{gark}}\tran\,  \b_{\textsc{gark}}.
\end{split}
\end{equation*}
Therefore, for symmetric and symplectic methods it holds that:
\begin{equation}
\label{cond.merge}
\begin{split}
\b_{\textsc{gark}}\,\mathcal{P}\,\A_{\textsc{gark}}\,\mathcal{P}  =  \mathbf{A}_{\textsc{gark}}\tran\,  \b_{\textsc{gark}}
\quad \Leftrightarrow \quad
b^{\{\ell\}}_i\, a^{\{\ell,m\}}_{s+1-i,s+1-j} = b^{\{m\}}_j\, a^{\{m,\ell\}}_{j,i}.
\end{split}
\end{equation}

%Symplecticness condition (for $\b^{\{\ell\}}_i \ne 0$) plus symmetry requires in the general case that:
%%
%\begin{eqnarray}
%0 &= & \b^{\{\ell\}-1}\, \mathbf{A}^{\{m,\ell\}}\tr\,  \b^{\{m\}} +\mathbf{A}^{\{\ell,m\}} -  \one \,\b^{\{m\}}\tr \nonumber \\
%&=&  \b^{\{\ell\}-1}\, \mathbf{A}^{\{m,\ell\}}\tr\,  \b^{\{m\}} - \mathcal{P}\,\A^{\{\ell,m\}}\,\mathcal{P} \nonumber \\
%%\mathcal{P}\,\mathbf{A}^{\{m,\ell\}}\tr\,  \b^{\{m\}} &=&  \b^{\{\ell\}}\,\A^{\{\ell,m\}}\,\mathcal{P} \nonumber \\
%\mathbf{A}^{\{m,\ell\}}\tr\,  \b^{\{m\}} &=& 
%%(\mathcal{P}\, \b^{\{\ell\}}\,\mathcal{P})\,
%\b^{\{\ell\}}\,
%(\mathcal{P}\,\A^{\{\ell,m\}}\,\mathcal{P})
%\label{cond.merge}\\
%b^{\{m\}}_j\, a^{\{m,\ell\}}_{j,i} &= &b^{\{\ell\}}_i\, a^{\{\ell,m\}}_{s+1-i,s+1-j}. \nonumber
%\end{eqnarray}

This condition together with symmetry implies symplecticness, and vice versa, together with symplecticness it implies symmetry.
\begin{equation*}
\{symmetric\} \cap \{symplectic\} \quad\Leftrightarrow\quad
symmetric  \cap  \eqref{cond.merge} \quad\Leftrightarrow\quad
\eqref{cond.merge}  \cap  symplectic.
\end{equation*}
\end{remark}

\ifreport
\begin{remark}
The condition~\eqref{cond.merge} is only necessary for symplecticness and symmetry, but not sufficient for either symplecticness or symmetry, as the following counterexample for $\nparts =1$ and two stages shows.  Condition~\eqref{cond.merge} reads in this case 
%%
%\begin{equation*}
%\begin{pmatrix}
%   0 & 1 \\ 1 & 0
%\end{pmatrix}
%\begin{pmatrix}
%    a_{1,1} & a_{2,1} \\
%    a_{1,2} & a_{2,2}
%\end{pmatrix}
%\begin{pmatrix}
%    b_1 & 0 \\ 0 & b_2 
%\end{pmatrix} =
%\begin{pmatrix}
%    b_1 & 0 \\ 0 & b_2 
%\end{pmatrix}
%\begin{pmatrix}
%    a_{1,1} & a_{1,2} \\
%    a_{2,1} & a_{2,2}
%\end{pmatrix}
%\begin{pmatrix}
%    0 & 1 \\ 1 & 0
%\end{pmatrix},
%\end{equation*}
%which demands the following conditions on the coefficients:
%\begin{align}
%    a_{2,2} & = \frac{b_1}{b_2} a_{1,1}
%    \label{cond.symmsympl}
%\end{align}
%
%{\color{red}
%
\begin{align}
\label{cond.symmsympl}
a_{2,2} = a_{1,1} \quad \textnormal{and} \quad b_1 = b_2.
\end{align}
%}
Symmetry requires that
\begin{equation*}
\begin{pmatrix}
    a_{2,2} & a_{2,1} \\
    a_{1,2} & a_{1,1}
\end{pmatrix}
+
\begin{pmatrix}
    a_{1,1} & a_{1,2} \\
    a_{2,1} & a_{2,2}
\end{pmatrix}
=
\begin{pmatrix}
    b_1 & b_2 \\ b_1 & b_2
\end{pmatrix},
\end{equation*}
holds, which yields %(for convergent schemes, i.e., $b_1=b_2=1/2$) 
the following symmetry conditions:
\begin{align}
\label{cond.symm}
\begin{split}
b_1 & = b_2, \\
a_{1,1}+a_{2,2} & = b_1, \\
a_{1,2}+a_{2,1} & = b_1.
\end{split}
\end{align}
Symplecticness is achieved if the following relations are fulfilled:
\begin{equation*}
\begin{pmatrix}
    \frac{1}{b_1} & 0 \\ 0 & \frac{1}{b_2}
\end{pmatrix}
\begin{pmatrix}
    a_{1,1} & a_{2,1} \\
    a_{1,2} & a_{2,2}
\end{pmatrix}
\begin{pmatrix}
    b_1 & 0 \\ 0 & b_2 
\end{pmatrix} +
\begin{pmatrix}
    a_{1,1} & a_{1,2} \\
    a_{2,1} & a_{2,2}
\end{pmatrix} -
\begin{pmatrix}
    b_1 & b_2\\ b_1 & b_2
\end{pmatrix}=0,
\end{equation*}
which yields %(for convergent schemes, i.e., $b_1=b_2=1/2$) 
the symplectic conditions
\begin{align}
\begin{split}
\label{cond.sympl}
    a_{1,1} & = \sfrac{b_1}{2}, \\
    a_{2,2} & = \sfrac{b_2}{2}, \\
    a_{1,2} b_1+ a_{2,1} b_2  & = b_1 b_2.
    \end{split}
\end{align}
Clearly, condition~\eqref{cond.symmsympl} neither implies~\eqref{cond.symm} nor~\eqref{cond.sympl}.
\end{remark}
\fi

When dealing with Hamiltonian systems, time-reversibility of a scheme $\Phi_h$ \eqref{eqn:numerical-reversible} is a desirable property.
In the following we show that the symmetry of a GARK scheme ~\eqref{eqn:numerical-reversible} ensures its
time-reversibilty, if each component is $\rho$-reversible. %a symplectic subsystem in its own right, i.e., the GARK scheme is of type \eqref{eqn:GARK-symplectic}.
\begin{theorem}[Symmetric GARK schemes are time-reversible]
\label{theorem.symmetry.time-reversibility}
A symmetric GARK scheme \eqref{eqn:GARK2}
%{eqn:GARK-symplectic} 
is $\rho$-reversible, provided that all individual components are $\rho$-reversible: 
\begin{leaveout}%
\[
\rho \circ \fun^{\{m\}}(\p,\q)= -\fun^{\{m\}}\big(\rho \circ(\p,\q)\big), \quad m = 1,\dots,\nparts.
\]
\end{leaveout}
{\[
\rho \circ \fun^{\{m\}}(\y)= -\fun^{\{m\}}\big(\rho \circ \y \big), \quad m = 1,\dots,\nparts.
\]}
\end{theorem}
\begin{leaveout}
\begin{proof}
Apply the GARK step~\eqref{eqn:GARK-symplectic} to the initial values $\rho(\q_0,\p_0)$ = $(\q_0,-\p_0)$ with step size $-h$ to obtain:
\begin{align*}
P_i^{\{q\}} & =  (-\p_0) + (-h) \sum_{m=1}^{\nparts} \sum_{j=1}^{s^{\{m\}}} a_{i,j}^{\{q,m\}} k_j^{\{m\}}, &
Q_i^{\{q\}} & =  \q_0 + (-h) \sum_{m=1}^{\nparts} \sum_{j=1}^{s^{\{m\}}} a_{i,j}^{\{q,m\}} \ell_j^{\{m\}},\\
\p_1 & =  (-\p_0) + (-h) \sum_{q=1}^{\nparts} \sum_{i=1}^{s^{\{q\}}} b_{i}^{\{q\}} k_i^{\{q\}}, &
\q_1 & =  \q_0 + (-h) \sum_{q=1}^{\nparts} \sum_{i=1}^{s^{\{q\}}} b_{i}^{\{q\}} \ell_i^{\{q\}}. 
\end{align*}
As the partitions of the Hamiltonian are $\rho$-reversible, %and therefore time-reversible
\begin{eqnarray*}
\begin{pmatrix}
k_i^{\{m\}} \\ \ell_i^{\{m\}} 
\end{pmatrix}  & = & 
\fun^{\{m\}}(P_i^{\{m\}},Q_i^{\{m\}}) = 
- \rho \circ \fun^{\{m\}}(\rho(P_i^{\{m\}},Q_i^{\{m\}})) \\
& = &
\begin{pmatrix} \Id_{d_q \times d_q} & \hphantom{-}\Zero_{d_q \times d_p} \\
 \Zero_{d_p \times d_q} & -\Id_{d_p \times d_p}  \end{pmatrix}
 \fun^{\{m\}}(-P_i^{\{m\}},Q_i^{\{m\}}),
%k_i^{\{m\}} & = &  - H^{\{m\}}_{\q}(P_i^{\{m\}},Q_i^{\{m\}}) = - H^{\{m\}}_{\q}(-P_i^{\{m\}},Q_i^{\{m\}}), 
%\label{eqn:GARK-symplectic.35}\\
%\ell_i^{\{m\}} & = &  H^{\{m\}}_{\p}(P_i^{\{m\}},Q_i^{\{m\}}) = -H^{\{m\}}_{\p}(-P_i^{\{m\}},Q_i^{\{m\}}). \label{eqn:GARK-symplectic.36}
\end{eqnarray*}
renaming the variables $P_i^{\{m\}} \coloneqq -P_i^{\{m\}}$ shifts the increments  $\ell_i^{\{m\}}$ to $-\ell_i^{\{m\}}$ and leaves the increments $k_i^{\{m\}}$ unchanged, which leads to the scheme:
\begin{align*}
P_i^{\{q\}} & =  \p_0 + h \sum_{m=1}^{\nparts} \sum_{j=1}^{s^{\{m\}}} a_{i,j}^{\{q,m\}} k_j^{\{m\}}, &
Q_i^{\{q\}} & =  \q_0 + h \sum_{m=1}^{\nparts} \sum_{j=1}^{s^{\{m\}}} a_{i,j}^{\{q,m\}} \ell_j^{\{m\}},\\
\p_1 & =  -\left(\p_0 + h \sum_{q=1}^{\nparts} \sum_{i=1}^{s^{\{q\}}} b_{i}^{\{q\}} k_i^{\{q\}} \right), &
\q_1 & =  \q_0 + h \sum_{q=1}^{\nparts} \sum_{i=1}^{s^{\{q\}}} b_{i}^{\{q\}} \ell_i^{\{q\}},
\end{align*}
and immediately shows that~\eqref{eqn:numerical-reversible} holds for the GARK scheme~\eqref{eqn:GARK-symplectic}.
%\qed
\end{proof}
\end{leaveout}
\begin{proof}
Apply the GARK step~\eqref{eqn:GARK-symplectic} to the initial values $\rho(\y_0)$  with step size $-h$ to obtain:
\begin{eqnarray*}
Y_i^{\{q\}} & =&  \rho(\y_0)  + (-h) \sum_{m=1}^{\nparts} \sum_{j=1}^{s^{\{m\}}} a_{i,j}^{\{q,m\}} \fun^{\{m\}}(Y_j^{\{m\}}), \\
\y_1 & = & \rho(\y_0)  + (-h) \sum_{q=1}^{\nparts} \sum_{i=1}^{s^{\{q\}}} b_{i}^{\{q\}} \fun^{\{q\}}(Y_i^{\{q\}}).
\end{eqnarray*}
As the partitions of the Hamiltonian are $\rho$-reversible, %and therefore time-reversible
\begin{eqnarray*}
\rho \circ \fun^{\{m\}}(Y_i^{\{m\}}) &= &
-\fun^{\{m\}}(\rho(Y_i^{\{m\}})) 
\end{eqnarray*}
renaming the internal variables $\widetilde Y_i^{\{m\}} \coloneqq \rho^{-1} (Y_i^{\{m\}})$ leads to the scheme:
\begin{eqnarray*}
\rho(\widetilde Y_i^{\{q\}})& =&  \rho(\y_0)  + (-h) \sum_{m=1}^{\nparts} \sum_{j=1}^{s^{\{m\}}} a_{i,j}^{\{q,m\}} \fun^{\{m\}}(\rho(\widetilde Y_j^{\{m\}})) \\
& = & \rho(\y_0)  + h \sum_{m=1}^{\nparts} \sum_{j=1}^{s^{\{m\}}} a_{i,j}^{\{q,m\}} \rho \circ \fun^{\{m\}}(\widetilde Y_j^{\{m\}}) , \\
\y_1 & = & \rho(\y_0)  + h \sum_{q=1}^{\nparts} \sum_{i=1}^{s^{\{q\}}} b_{i}^{\{q\}} \rho \circ \fun^{\{q\}}(\widetilde Y_i^{\{q\}}),
\end{eqnarray*}
which yields for $\rho$ linear and regular (see, for example, the mapping given in~\eqref{rhorev}),
\begin{eqnarray*}
\widetilde Y_i^{\{q\}}& =&  \y_0  + h \sum_{m=1}^{\nparts} \sum_{j=1}^{s^{\{m\}}} a_{i,j}^{\{q,m\}} \fun^{\{m\}}(\widetilde Y_j^{\{m\}}) , \\
\y_1 & = & \rho \left( \y_0  + h \sum_{q=1}^{\nparts} \sum_{i=1}^{s^{\{q\}}} b_{i}^{\{q\}} \rho \circ \fun^{\{q\}}(\widetilde Y_i^{\{q\}}) \right),
\end{eqnarray*}
which immediately shows that~\eqref{eqn:numerical-reversible} holds for the GARK scheme~\eqref{eqn:GARK2}.%\qed
\end{proof}

\begin{remark}
For a Hamiltonian splitting \eqref{eqn:partitioned-sympl.system} with $\y=(\p,\q)^\top$, we have 
\begin{eqnarray*}
\rho \circ \fun^{\{m\}} (\p,\q) = 
\begin{pmatrix}
H^{\{m\}}_{\q}(\p,\q) \\
H^{\{m\}}_{\p}(\p,\q)
\end{pmatrix} = 
\begin{pmatrix}
H^{\{m\}}_{\q}(-\p,\q) \\
-H^{\{m\}}_{\p}(-\p,\q)
\end{pmatrix} = 
- \fun^{\{m\}} (\rho \circ (\p,\q)),
\end{eqnarray*}
i.e., all components define $\rho$-reversible flows. Hence the  GARK method~\eqref{eqn:GARK2} is time-reversible.
\end{remark}
\begin{comment}
The proof can be generalized to the setting of $\rho$-reversibility
with $\rho(\q,\p)=(\rho_1(q),\rho_2(p))$.
%
\begin{corollary}
Symmetric GARK schemes are $\rho$-reversible for $\rho(p,q)=(\rho_1 \cdot p,\rho_2 \cdot  q)$ with linear invertible transformations $\rho_1,\rho_2$.
\end{corollary}
%
\begin{proof}
The proof follows the line of the proof of theorem~\ref{theorem.symmetry.time-reversibility} , replacing $k_i^{\{m\}}$ and $l_i^{\{m\}}$ by
$-\rho_1 \cdot k_i^{\{m\}}$ and $-\rho_2 \cdot \ell_i^{\{m\}}$, respectively. \sandu{$+\rho_2 \cdot \ell_i^{\{m\}}$?}
\end{proof}
\end{comment}

We finish this section with two examples of symmetric and/or symplectic schemes for $\nparts=2$ partitions.
\begin{example}[A symplectic implicit-implicit scheme]
Consider the GARK scheme~\eqref{eqn:GARK2} defined by the generalized Butcher tableau \eqref{eqn:general-Butcher-tableau} %\sandu{oder?}
\begin{equation*}
\renewcommand{\arraystretch}{1.3}
%\begin{array}{c|c|c|c}
%\mathbf{c}^{\{1,1\}} & \mathbf{A}^{\{1,1\}} &\mathbf{A}^{\{1,2\}} & 
%\mathbf{c}^{\{1,2\}} \\
%\hline 
%\mathbf{c}^{\{2,1\}} & \mathbf{A}^{\{2,1\}} & \mathbf{A}^{\{2,2\}}  & \mathbf{c}^{\{2,2\}} \\
%\hline
%& \b^{\{1\}}\tr & \b^{\{2\}}\tr
%\end{array}
~~\raisebox{8pt}{$\begin{array}{c|c}
\mathbf{A}^{\{1,1\}} &\mathbf{A}^{\{1,2\}} \\
\hline 
\mathbf{A}^{\{2,1\}} & \mathbf{A}^{\{2,2\}}   \\
\Xhline{2\arrayrulewidth}
 \b^{\{1\}}\tr & \b^{\{2\}}\tr
\end{array}
$}
~~=~~
%\begin{array}{c|cc|cc|c} 
%\frac{1}{8} ~&~ \frac{1}{8} ~&~ 0 ~&~ 0 ~&~ 0 ~&~0 \\ 
%\frac{5}{8} ~&~ \frac{1}{4} ~&~ \frac{3}{8} ~&~ \frac{2}{3} ~&~ 0 ~&~  \frac{2}{3} \\ 
%\hline
%\frac{1}{4}  ~&~ \frac{1}{4} ~&~ 0 ~&~ \frac{1}{3} ~&~ 0 ~&~ \frac{1}{3}\\ 
%1 ~&~ \frac{1}{4} ~&~ \frac{3}{4}  ~&~ \frac{2}{3} ~&~ \frac{1}{6} ~&~ \frac{5}{6}\\
%\hline
%~&~ \frac{1}{4} ~&~ \frac{3}{4} ~&~ \frac{2}{3} ~&~ \frac{1}{3} ~&~  
%\end{array}~.
\raisebox{24pt}{$
\begin{array}{cc|cc} 
\frac{1}{8} ~&~\scriptstyle 0 ~&~\scriptstyle  0 ~&~\scriptstyle  0  \\ 
 \frac{1}{4} ~&~ \frac{3}{8} ~&~ \frac{2}{3} ~&~\scriptstyle  0 \\ 
\hline
\frac{1}{4} ~&~\scriptstyle  0 ~&~ \frac{1}{3} ~&~\scriptstyle  0 \\ 
 \frac{1}{4} ~&~ \frac{3}{4}  ~&~ \frac{2}{3} ~&~ \frac{1}{6} \\
\Xhline{2\arrayrulewidth}
 \frac{1}{4} ~&~ \frac{3}{4} ~&~ \frac{2}{3} ~&~ \frac{1}{3}  
\end{array}$}\;.
\end{equation*}
\renewcommand{\arraystretch}{1}

This scheme is symplectic for a Hamiltonian splitting according to Theorem~\ref{theorem.symplectic}, and of second order. However, 
%it is straightforward to show that 
%a second order DIRK scheme of type~\eqref{eqn:imimRK} with two stages cannot be symmetric and symplectic at the same time.
the scheme is neither internally consistent nor symmetric.
\end{example}

    \begin{example}[A symplectic and symmetric  implicit-implicit GARK]
\label{example-GARK-IMIM}

An example of a symmetric and symplectic GARK method of order two for a general splitting \eqref{eqn:additive-ode} (the component subsystems are not necessarily Hamiltonian, i.e., \eqref{eqn:partitioned-sympl.system} may not hold), based on the Verlet scheme in the coupling parts, is given by the following Butcher tableau~\eqref{eqn:general-Butcher-tableau}
\begin{equation}
\label{eqn:GARK-IMIM2}
\renewcommand{\arraystretch}{1.3}
\raisebox{8pt}{$\begin{array}{c|c}
\mathbf{A}^{\{1,1\}} &\mathbf{A}^{\{1,2\}} \\
\hline 
\mathbf{A}^{\{2,1\}} & \mathbf{A}^{\{2,2\}}   \\
\Xhline{2\arrayrulewidth}
 \b^{\{1\}}\tr & \b^{\{2\}}\tr
\end{array}$}
~~=~~
\raisebox{23pt}{$
\begin{array}{cc|cc} 
\frac{1}{4} &\scriptstyle \alpha  & a_{1,1} & a_{1,2} \\ 
\frac{1}{2} - {\scriptstyle \alpha} & \frac{1}{4} & a_{2,1} & a_{2,2} \\ 
\hline
\frac{1}{2} - a_{1,1}  &  \frac{1}{2} - a_{2,1}  & \frac{1}{4} &\scriptstyle \beta \\ 
\frac{1}{2} - a_{1,2} &  \frac{1}{2}  - a_{2,2} & \frac{1}{2} - {\scriptstyle \beta} & \frac{1}{4} \\
\Xhline{2\arrayrulewidth}
\frac{1}{2} & \frac{1}{2} & \frac{1}{2} & \frac{1}{2} 
\end{array}$}\;,
\end{equation}
%\begin{equation}
%\label{eqn:GARK-IMIM}
%\renewcommand{\arraystretch}{1.3}
%\raisebox{8pt}{$\begin{array}{c|c}
%\mathbf{A}^{\{1,1\}} &\mathbf{A}^{\{1,2\}} \\
%\hline 
%\mathbf{A}^{\{2,1\}} & \mathbf{A}^{\{2,2\}}   \\
%\Xhline{2\arrayrulewidth}
% \b^{\{1\}}\tr & \b^{\{2\}}\tr
%\end{array}$}
%~~=~~
%\raisebox{23pt}{$
%\begin{array}{cc|cc} %
%\frac{1}{4} &\scriptstyle \alpha  & \scriptstyle 0 & \scriptstyle0 %\\ 
%\frac{1}{2} - {\scriptstyle \alpha} & \frac{1}{4} & \frac{1}{2} & \frac{1}{2} \\ 
%\hline
% \frac{1}{2} & \scriptstyle 0 & \frac{1}{4} &\scriptstyle \beta \\ 
% \frac{1}{2} & \scriptstyle 0 & \frac{1}{2} - {\scriptstyle \beta} & \frac{1}{4} \\
%\Xhline{2\arrayrulewidth}
%\frac{1}{2} & \frac{1}{2} & \frac{1}{2} & \frac{1}{2} 
%\end{array}$}\;,
%\end{equation}
%
where $\alpha$, $\beta$, $a_{1,1},\ldots,a_{2,2}$ are free real parameters with $a_{1,1}+a_{2,2}=1/2$ and $a_{1,2}+a_{2,1}=1/2$.
\end{example}
    \begin{remark}
This GARK scheme has the following properties: 
\begin{itemize}
    \item it is decoupled for $\alpha=\beta=0$ provided that $\mathbf{A}^{\{1,2\}} \times \mathbf{A}^{\{2,1\}\top} = \mathbf{0}$ holds, when it becomes a DIRK-DIRK scheme; 
    \item it is internally consistent, if  the four conditions $1/4 + \alpha = a_{1,1} + a_{1,2}$, $3/4 - \alpha = a_{2,1} + a_{2,2}$, $1/4 + \beta = 1-a_{1,1}-a_{2,1}$, $3/4 - \beta = 1 - a_{1,2} - a_{2,2}$ are satisfied; 
    \item in general, if $\mathbf{A}^{\{1,2\}} \times \mathbf{A}^{\{2,1\}\top} = \mathbf{0}$ does not hold, it is NOT a composition scheme; 
    \item One notes that this scheme, when applied to a separable Hamiltonian system (see Section \ref{sec-partitionedsymplecticgark}) in a coordinate partitioning way, is equivalent to the original Verlet scheme, if $a_{1,1}=a_{2,1}=1/2$ and $a_{1,2} =a_{2,2}=0$ or vice versa holds. 
    \end{itemize}
    \end{remark}
    \begin{remark}
     Whereas in the case of additive Runge-Kutta schemes  the component sums $\fun^{\{1\}}\left(P_i,Q_i\right) + \fun^{\{2\}}\left(P_i,Q_i\right)$ = $\fun\left(P_i,Q_i\right) $ equal the total right-hand side, we have in the GARK case different arguments and the components do not add to the total right-hand side in general. Consequently, the symplectic GARK is not equivalent to a single RK scheme applied to the non-partitioned system.  This is also the case in Example~\ref{example-GARK-IMIM}. For the choice $\alpha=\beta=a_{1,1}=a_{1,2}=0$ and $a_{2,1}=a_{2,2}=1/2$ the GARK scheme~\eqref{eqn:GARK2} reads 
     with $\fun_i^{\{m\}}  \coloneqq \fun^{\{m\}}\left(Y_i^{\{m\}}\right)$: %, \fun_i \coloneqq \sum_m \fun_i^{\{m\}}$
$$
\begin{array}{lcllcl}
Y_1^{\{1\}}  &=& \y_0 +  \sfrac{h}{4} { \fun_1^{\{1\}}  }, & \quad 
Y_1^{\{2\}} 
&=& Y_1^{\{1\}} + \sfrac{h}{4} \fun_1^{\{1\}}    + \sfrac{h}{4} { \fun_1^{\{2\}} }, %\neq  Y_1^{\{1\}} + \sfrac{h}{4} \fun_1, 
 \\[1ex]
Y_2^{\{2\}}  
%=  \y_0 + \sfrac{h}{2}  \fun_1^{\{1\}} + \sfrac{h}{2} \fun_1^{\{2\}}  +  \sfrac{h}{4} \fun_2^{\{2\}}
  &=& Y_1^{\{2\}} + \sfrac{h}{4}  \fun_1^{\{2\}}    + \sfrac{h}{4} \fun_2^{\{2\}}, & \quad 
%  \\
%
Y_2^{\{1\}} 
%& = \y_0 + \sfrac{h}{2} \fun_1^{\{1\}}   
%+ \sfrac{h}{2}  \fun_1^{\{2\}}   
%+ \sfrac{h}{2}  \fun_2^{\{2\}}  + \sfrac{h}{4} { \fun_2^{\{1\}} }  \\
&=& Y_2^{\{2\}}  + \sfrac{h}{4}  \fun_2^{\{1\}}  + \sfrac{h}{4} { \fun_2^{\{2\}} }, %\ne  Y_2^{\{2\}}  + \sfrac{h}{4}  \fun_2,  
%\\
%%
%\y_1  &=& \y_0 + \sfrac{h}{2} ( \fun_1^{\{1\}} + \fun_1^{\{2\}} )
% + \sfrac{h}{2} ( \fun_2^{\{1\}} + \fun_2^{\{2\}} ) %\neq \sfrac{h}{2} ( \fun_1+ \fun_2 ).
\end{array}
$$
$$
\y_1  = \y_0 + \sfrac{h}{2} ( \fun_1^{\{1\}} + \fun_1^{\{2\}} )
 + \sfrac{h}{2} ( \fun_2^{\{1\}} + \fun_2^{\{2\}} ). %\neq \sfrac{h}{2} ( \fun_1+ \fun_2 ).
$$
Here $\fun_1^{\{1\}}+\fun_1^{\{2\}}=\fun^{\{1\}}\left(Y_1^{\{1\}}\right)+\fun^{\{2\}}\left(Y_1^{\{2\}}\right)$, as well as  
$\fun_2^{\{1\}}+\fun_2^{\{2\}}=\fun^{\{1\}}\left(Y_2^{\{1\}}\right)+\fun^{\{2\}}\left(Y_2^{\{2\}}\right)$ are evaluated at different stage values and are thus not equivalent to a single evaluation of the function $\fun(\cdot)$ at the same stage value as in the ARK case. 
\end{remark}

\subsection{Backward error analysis}
Performing a backward error analysis, one can show \cite{Hairer_2006_geometric-book} that the modified equation of symplectic numerical integration schemes, applied to the Hamiltonian system \eqref{eq-ham}, is also Hamiltonian. Consequently, the scheme preserves a nearby \textit{shadow Hamiltonian}. Consider the GARK scheme \eqref{eqn:GARK2} whose numerical solution $\y_1 = \Phi_h(\y_0)$, written as NB-series, is given by \eqref{eq:NB-series}. In terms of backward error analysis, the numerical solution can be regarded as the exact solution to the modified system 
\begin{equation}\label{eq:modified_system}
    \dot{\tilde{\y}} = \sum\limits_{\tree \in \NT} \mathfrak{b}(\tree) \frac{h^{\rho(\tree) - 1}}{\sigma(\tree)} F(\tree)(\tilde{\y}), 
\end{equation}
with elementary differentials $F(\tree)(\y)$ defined recursively via \eqref{eq:elem-diffs}. Defining the set of all splittings
\begin{equation*}
    \mathrm{SP}(\tree) := \{ \theta \in \mathrm{OST}(\tree)\ \vert \ \tree \setminus \theta \text{ consists of a single element}\},
\end{equation*}
with $\mathrm{OST}(\tree)$ being the set of ordered subtrees, the real coefficients $\mathfrak{b}(\tree)$ are recursively defined by $\mathfrak{b}(\emptyset) = 0,\ \mathfrak{b}(\tau_{\{m\}}) = 1$ and 
\begin{equation*}
    \mathfrak{b}(\tree) = \mathfrak{a}(\tree) - \sum\limits_{j=2}^{\rho(\tree)} \frac{1}{j!} \partial_{\mathfrak{b}}^{j-1} \mathfrak{b}(\tree) \quad \text{for} \quad \tree \in \NT.
\end{equation*}
Here, $\partial_{\mathfrak{b}}^{j-1}$ denotes the $(j-1)$-th iterate of the Lie derivative 
\begin{equation*}
    \partial_{\mathfrak{b}} \mathfrak{c}(\tree) = \sum\limits_{\theta \in \mathrm{SP}(\tree)} \mathfrak{c}(\theta) \mathfrak{b}(\tree \setminus \theta).
\end{equation*}
For a given smooth Hamiltonian function $H : \mathbb{R}^d \to \mathbb{R}$ and for $\tree \in \NT$, the elementary Hamiltonian $\mathcal{H}(\tree) : \mathbb{R}^d \to \mathbb{R}$ is given by 
\begin{equation*}
    \begin{split}
        \mathcal{H}(\tau_{\{m\}})(\y) &= H^{\{m\}}(\y), \\
        \mathcal{H}(\tree)(\y) &= H^{\{m\}(r)}(\y)(F(\tree_1)(\y),\ldots,F(\tree_r)(\y)), \quad \tree = [\tree_1,\ldots,\tree_r]_{\{m\}}.
    \end{split}
\end{equation*}
Similar to investigations for B-series and P-series \cite{Hairer_2006_geometric-book}, we select representatives from the equivalence class $u \circ v \sim v \circ u$, resulting in the set 
\begin{equation*}
    \mathbb{T}_{\mathrm{N}}^{*} = \{\tau_{\{1\}},\ldots,\tau_{\{\mathrm{N}\}} \} \cup \bigg\{ \tree \in \NT \bigg\vert \parbox{7cm}{$\tree$ cannot be written as $\tree = u \circ v$ with $u < v$, also not if the color of the root is changed.}\bigg\}.
\end{equation*}
Then, the modified system \eqref{eq:modified_system} is Hamiltonian with shadow 
\begin{equation}\label{eq:shadow_Hamiltonian}
    \tilde{H}(\y) = \sum\limits_{k=1}^{\infty} h^{k-1} H_k(\y),\quad \text{with\ \ } H_k(\y) = \sum\limits_{\tree \in \mathbb{T}_{\mathrm{N}}^{*},\rho(\tree) = k} \frac{\mathfrak{b}(\tree)}{\sigma(\tree)} \mathcal{H}(\tree)(\y).  
\end{equation}

\begin{example} 
Consider the symplectic and symmetric implicit-implicit GARK scheme given by the Butcher tableau \eqref{eqn:GARK-IMIM2} with $a_{1,2}=a_{2,2}=0$ and $a_{1,1} = a_{2,1} = \tfrac{1}{2}$. The scheme preserves the shadow Hamiltonian
\begin{equation*}
\begin{split}
    \tilde{H} &= H + h^2 \Big( \left(\tfrac{\alpha^2}{2} - \tfrac{\alpha}{4} - \tfrac{1}{96}\right) \mathcal{H}([\tau_{\{1\}},\tau_{\{1\}}]_{\{1\}}) - \tfrac{1}{12} \mathcal{H}([\tau_{\{1\}},\tau_{\{2\}}]_{\{1\}}) \\
     &\qquad - \tfrac{1}{24} \mathcal{H}([\tau_{\{2\}},\tau_{\{2\}}]_{\{1\}}) + \left(\tfrac{\beta^2}{2} - \tfrac{\beta}{4} - \tfrac{1}{96}\right) \mathcal{H}([\tau_{\{2\}},\tau_{\{2\}}]_{\{2\}})\\
    &\qquad +  (\tfrac{1}{24} - \tfrac{\beta}{2}) \mathcal{H}([\tau_{\{2\}},\tau_{\{1\}}]_{\{2\}}) + \tfrac{1}{12} \mathcal{H}([\tau_{\{1\}},\tau_{\{1\}}]_{\{2\}}) \Big) + \mathcal{O}(h^4).
\end{split}
\end{equation*}
\end{example}

\begin{leaveout}
\sandu{
In the general scenario~\eqref{eqn:GARK2}, e.g., in the case of component-wise splitting discussed in Example~\ref{example.two.mass.osc}, we get 
\begin{subequations}
\begin{align}
Y^{\{1\}}_1 
 & = 
\y_0
 +  \sfrac{h}{4} 
\fun^{\{1\}}_1
  + h \alpha 
\fun^{\{1\}}_2
 \\
Y_1^{\{2\}} 
 & = 
\y_0
+ \sfrac{h}{2}  \fun^{\{1\}}_1  
 +  \sfrac{h}{4}  \fun^{\{2\}}_1
 + h \beta  \fun^{\{2\}}_2
 \\
Y_2^{\{1\}} 
 & = 
\y_0
 + h (\sfrac{1}{2}-\alpha)  \fun^{\{1\}}_1  
+ \sfrac{h}{4}  \fun^{\{1\}}_2   + \sfrac{h}{2}  \fun^{\{2\}}_1  
+ \sfrac{h}{2}  \fun^{\{2\}}_2   \\
Y_2^{\{2\}} 
 & = 
\y_0
 + \sfrac{h}{2} \fun_1^{\{1\}}
+ h (\sfrac{1}{2}-\beta) \fun_1^{\{2\}}
  +  \sfrac{h}{4} \fun_2^{\{2\}} \\
\fun^{\{m\}}_i & = \fun^{\{m\}}\left(Y_i^{\{m\}}\right), \\
\y_1& = 
\y_0
 + \sfrac{h}{2} 
\fun_1^{\{1\}} + \fun_1^{\{2\}}  
 + \sfrac{h}{2} 
 \fun_2^{\{1\}} + \fun_2^{\{2\}} . 
\end{align}
\end{subequations}
Whereas in the ARK case the component sums $\fun^{\{1\}}\left(P_i,Q_i\right) + \fun^{\{2\}}\left(P_i,Q_i\right)$ = $\fun\left(P_i,Q_i\right) $ equal the total right-hand side, we have in the GARK case different arguments and the components do not add to the total right-hand side in general. \sandu{Consequently, the symplectic GARK is not equivalent to a single RK scheme applied to the non-partitioned system.}
%I

%
A generic representation with $\alpha=\beta=0$ shows the decoupled nature:
\begin{align*}
\fun_i^{\{m\}} & \coloneqq \fun^{\{m\}}\left(Y_i^{\{m\}}\right), \quad
\fun_i \coloneqq \sum_m \fun_i^{\{m\}}, \\
Y_1^{\{1\}} & = \y_0 +  \sfrac{h}{4} { \fun_1^{\{1\}}  },
 \\
Y_1^{\{2\}} & %= \y_0 + \sfrac{h}{2}  \fun_1^{\{1\}}    +  \sfrac{h}{4} { \fun_1^{\{2\}} }
= Y_1^{\{1\}} + \sfrac{h}{4} \fun_1^{\{1\}}    + \sfrac{h}{4} { \fun_1^{\{2\}} } \neq  Y_1^{\{1\}} + \sfrac{h}{4} \fun_1, 
 \\
Y_2^{\{2\}} & 
%=  \y_0 + \sfrac{h}{2}  \fun_1^{\{1\}} + \sfrac{h}{2} \fun_1^{\{2\}}  +  \sfrac{h}{4} \fun_2^{\{2\}}
  = Y_1^{\{2\}} + \sfrac{h}{4}  \fun_1^{\{2\}}    + \sfrac{h}{4} \fun_2^{\{2\}}, \\
%  \\
%
Y_2^{\{1\}} 
%& = \y_0 + \sfrac{h}{2} \fun_1^{\{1\}}   
%+ \sfrac{h}{2}  \fun_1^{\{2\}}   
%+ \sfrac{h}{2}  \fun_2^{\{2\}}  + \sfrac{h}{4} { \fun_2^{\{1\}} }  \\
&= Y_2^{\{2\}}  + \sfrac{h}{4}  \fun_2^{\{1\}}  + \sfrac{h}{4} { \fun_2^{\{2\}} } \ne  Y_2^{\{2\}}  + \sfrac{h}{4}  \fun_2,  \\
\y_1 & = \y_0 + \sfrac{h}{2} ( \fun_1^{\{1\}} + \fun_1^{\{2\}} )
 + \sfrac{h}{2} ( \fun_2^{\{1\}} + \fun_2^{\{2\}} ) \neq \sfrac{h}{2} ( \fun_1+ \fun_2 ).
\end{align*}
}
\end{leaveout}

\section{Partitioned GARK schemes  for separable Hamiltonian systems}
\label{sec-partitionedsymplecticgark}
%%%%%%%%%%%%%%%%%%%%%%%%%

In this section we consider schemes for separable Hamiltonians $H(\q,\p)=T(\p)+V(\q)$.
We discuss two types of partitioned Hamiltonians, first when both the potential part $V(\q)$ and kinetic part $T(\p)$ are split, and second when only the potential is split.

%%%%%%%%%%%%%%%%%
\subsection{Partitioned symplectic GARK schemes for kinetic and potential splitting}

We consider systems where both the potential and the kinetic parts are split:
    \begin{subequations}
       \label{eqn:split-TV}
    \begin{equation}
       \label{eqn:split-TV-Hamiltonian-sum}
    H(\p,\q)= \sum_{m=1}^{\nparts} \left( T^{\{m\}}(\p) + V^{\{m\}}(\q) \right),
    \end{equation}
     we consider the $2\nparts-way$ partitioned Hamiltonian \eqref{eqn:partitioned-Hamiltonians}
     \begin{align}
       \label{eqn:split-TV-Hamiltonian}
        H(\p,\q)= &\sum_{m=1}^{2\nparts} H^{\{m\}}(\p,\q) \quad \mbox{with} \quad
        \begin{cases}
        H^{\{m\}}(\p,\q)  = T^{\{m\}}(\p), & m=1,\ldots,\nparts, \\[3pt]
        H^{\{m+\nparts\}}(\p,\q)  = V^{\{m\}}(\q), & m=1,\ldots,\nparts.
        \end{cases}
%        T(\p) &=\sum_{i=1}^{\nparts} T^{\{i\}}(\p), \\
%        V(\q) &=\sum_{i=1}^{\nparts} V^{\{i\}}(\q).
    \end{align}   
    \end{subequations}
The GARK scheme \eqref{eqn:GARK-symplectic} applied to a system with splitting~\eqref{eqn:split-TV} reads: 
\begin{equation}
\label{eq:sympl.part.gark.2}
\begin{split}
\widetilde P_i^{\{q\}} & =  \p_0 + h \sum_{m=1}^{\nparts} \sum_{j=1}^{s^{\{\nparts + m\}}} \widetilde a_{i,j}^{\{q,\nparts+m\}} \widetilde k_j^{\{\nparts+m\}}, \\
\widetilde Q_i^{\{q\}} & =  \q_0 + h \sum_{m=1}^{\nparts} \sum_{j=1}^{s^{\{m\}}} \widetilde  a_{i,j}^{\{q,m\}} \widetilde \ell_j^{\{m\}},\\
\p_1 & =  \p_0 + h \sum_{q=1}^{\nparts} \sum_{i=1}^{s^{\{\nparts + q\}}} \widetilde b_{i}^{\{\nparts+q\}} \widetilde k_i^{\{\nparts+q\}}, \\
\q_1 & =  \q_0 + h \sum_{q=1}^{\nparts} \sum_{i=1}^{s^{\{q\}}} \widetilde{b}_{i}^{\{q\}} \widetilde \ell_i^{\{q\}}, \\
\widetilde k_i^{\{\nparts+m\}} & %=  - H^{\{\nparts+m\}}_{\q}\left(\widetilde P_i^{\{\nparts+m\}},\widetilde Q_i^{\{\nparts+m\}}\right)
= - V_{\q}^{\{m\}}(\widetilde Q_i^{\{\nparts+m\}}), \\
\widetilde \ell_i^{\{m\}} & %=  H^{\{m\}}_{\p}\left(\widetilde P_i^{\{m\}},\widetilde Q_i^{\{m\}}\right)
= T_{\p}^{\{m\}}(\widetilde P_i^{\{m\}}).
\end{split}
\end{equation}
The stage vectors $\widetilde Q_i^{\{\ell\}}$ and $\widetilde P_i^{\{\nparts+\ell\}}$ are not needed for any $\ell=1,\ldots,\nparts$. Using the notation $P_i^{\{m\}}\coloneqq\widetilde P_i^{\{m\}}$,  $Q_i^{\{m\}}\coloneqq\widetilde Q_i^{\{\nparts+m\}}$, $\hat s^{\{q\}}\coloneqq s^{\{\nparts+q\}}$, and
\begin{equation*}
\Ahat^{\{\ell,m\}}\coloneqq\widetilde{\mathbf{A}}^{\{\ell,\nparts+m\}}, \quad
 {\mathbf{A}}^{\{\ell,m\}}\coloneqq\widetilde{\mathbf{A}}^{\{\nparts+\ell,m\}}, \quad
 \bhat^{\{m\}}\coloneqq \widetilde{\b}^{\{\nparts+m\}}, \quad
 {\b}^{\{m\}}\coloneqq \widetilde{\b}^{\{m\}},
\end{equation*}
for $m=1,\ldots,\nparts$, the partitioned GARK scheme \eqref{eq:sympl.part.gark.2} reads:
\begin{subequations}
\label{eq:sympl.part.gark}
\begin{align}
P_i^{\{q\}} & =  \p_0 + h \sum_{m=1}^{\nparts} \sum_{j=1}^{\hat s^{\{m\}}} \ahat_{i,j}^{\{q,m\}} k_j^{\{m\}},
%, \qquad i=1,\ldots,s^{\{q\}},\\
& Q_i^{\{q\}} & =  \q_0 + h \sum_{m=1}^{\nparts} \sum_{j=1}^{s^{\{m\}}}  a_{i,j}^{\{q,m\}} \ell_j^{\{m\}}, \\
%\qquad i=1,\ldots,\hat s^{\{q\}},\\
\p_1 & =  \p_0 + h \sum_{q=1}^{\nparts} \sum_{i=1}^{\hat s^{\{q\}}} \widehat{b}_{i}^{\{q\}} k_i^{\{q\}}, &
\q_1 & =  \q_0 + h \sum_{q=1}^{\nparts} \sum_{i=1}^{s^{\{q\}}} {b}_{i}^{\{q\}} \ell_i^{\{q\}}, \\
k_i^{\{m\}} & =  -V_{\q}^{\{m\}}(Q_i^{\{m\}}), &
\ell_i^{\{m\}} & =  T_{\p}^{\{m\}}(P_i^{\{m\}}).
\end{align}
\end{subequations}
This scheme can also be analyzed in the framework presented in this section.
The corresponding generalized Butcher tableau \eqref{eqn:general-Butcher-tableau} is:
\begin{equation}
\label{eqn:Partitioned-Butcher-tableau}
\renewcommand{\arraystretch}{1.3}
\begin{array}{c|c}
\Zero & \Ahat \\ \hline
\mathbf{A} & \Zero
\\
\Xhline{2\arrayrulewidth}
\b \tr & \bhat\tr
\end{array}
~~\raisebox{-8pt}{$=$}~~
\raisebox{35pt}{$
\begin{array}{ccc|ccc}
\Zero &  \ldots & \Zero  &
\Ahat^{\{1,1\}} &  \ldots & \Ahat^{\{1,\nparts\}} 
\\
\vdots & \ddots & \vdots & \vdots & \ddots & \vdots 
\\
\Zero &  \ldots & \Zero  &
\Ahat^{\{\nparts,1\}} &  \ldots & \Ahat^{\{\nparts,\nparts\}} \\
\hline
\mathbf{A}^{\{1,1\}} &  \ldots & \mathbf{A}^{\{1,\nparts\}} 
& \Zero &  \ldots & \Zero 
\\
\vdots & \ddots & \vdots & \vdots & \ddots & \vdots
\\
\mathbf{A}^{\{\nparts,1\}} &  \ldots & \mathbf{A}^{\{\nparts,\nparts\}} 
& \Zero &  \ldots & \Zero \\
\Xhline{2\arrayrulewidth}
\b ^{\{1\}}\tr &  \ldots &\b ^{\{\nparts\}}\tr &
\bhat^{\{1\}}\tr &  \ldots &\bhat^{\{\nparts\}}\tr 
\end{array}$}. 
\end{equation}
%

%where the empty positions can be filled with arbitrary coefficients, and
%%
%\begin{align*} 
%\mathbf{A} &\coloneqq 
%\begin{bmatrix}
%\mathbf{A}^{\{1,1\}} &  \ldots & \mathbf{A}^{\{1,\nparts\}} 
%\\
%\vdots & \ddots & \vdots 
%\\
%\mathbf{A}^{\{\nparts,1\}} &  \ldots & \mathbf{A}^{\{\nparts,\nparts\}} 
%\end{bmatrix}, 
%&\mathbf{\hat A} &\coloneqq 
%\begin{bmatrix}
%\mathbf{\hat A}^{\{1,1\}} &  \ldots & \mathbf{\hat A}^{\{1,\nparts\}} 
%\\
%\vdots & \ddots & \vdots 
%\\
%\mathbf{\hat A}^{\{\nparts,1\}} &  \ldots & \mathbf{\hat A}^{\{\nparts,\nparts\}} 
%\end{bmatrix}, 
%\\
%%
%\b \tr &\coloneqq 
%\begin{bmatrix}
%\b ^{\{1\}}\tr &  \ldots &\b ^{\{\nparts\}} \tr
%\end{bmatrix}, 
%& \mathbf{\hat  b}\tr &\coloneqq 
%\begin{bmatrix}
%\mathbf{\hat  b}^{\{1\}}\tr &  \ldots &\bhat^{\{\nparts\}}\tr 
%\end{bmatrix}.
%\end{align*}
%
%An argument similar to the one used in the proof of Corollary~\ref{cor.part.ham} \sandu{Now this Corollary is later, we need to move the main argument  here} shows that %
%The necessary and sufficient condition \eqref{eqn:GARK-symplectic-condition} for symplecticness of the scheme \eqref{eqn:Partitioned-Butcher-tableau} is:
%%
%\begin{align}
%\label{eqn:symplecticness}
%\Ahat^{\{\ell,m\}}\tr\,  \B^{\{\ell\}} + \Bhat^{\{m\}}\, \mathbf{A}^{\{m,\ell\}} -  \bhat^{\{m\}} \,\b^{\{\ell\}}\tr = \Zero,
%\quad \ell,m = 1,\dots,\nparts,
%\end{align}
%%
%which is shown in eqn:GARK-symplectic.sta
\begin{corollary}[Symplecticity]
\label{cor.part.ham}
The necessary and sufficient conditions for the symplecticity of the GARK scheme~\eqref{eq:sympl.part.gark.2}
are $\mathbf{P}^{\{\ell,m\}}=0$ for $\ell \in \{1,\ldots,N\}$ and $m \in \{N+1,\ldots, 2\nparts\}$. Using notation \eqref{eqn:Partitioned-Butcher-tableau}, the necessary and sufficient condition \eqref{eqn:GARK-symplectic-condition} for symplecticness is:
\begin{align}
\label{eqn:symplecticness}
\Ahat^{\{\ell,m\}}\tr\,  \B^{\{\ell\}} + \Bhat^{\{m\}}\, \mathbf{A}^{\{m,\ell\}} -  \bhat^{\{m\}} \,\b^{\{\ell\}}\tr = \Zero,
\quad \ell,m = 1,\dots,\nparts.
\end{align}
\end{corollary}
\begin{proof}
From \eqref{eq:sympl.part.gark.2}  we see that 
$\widetilde k_i^{\{\ell\}}=0$ and $\widetilde \ell_i^{\{\nparts+\ell\}}=0$ for $\ell=1,\ldots,\nparts$.
 Consequently, using the partition
\[
\mathfrak{g}_i^{\{m\}} = \begin{bmatrix}
\mathfrak{g}_{k,i}^{\{m\}} \\
\mathfrak{g}_{l,i}^{\{m\}}
\end{bmatrix}, \quad 1 \le m \le  2\nparts,
\]
in the proof of Theorem~\ref{theorem.symplectic} gives:
%the NB series of $k_i^{\{s\}}$ and $\ell_i^{\{s\}}$, $1 \le s \le  2\nparts$,  can be written as
\[
\begin{bmatrix}
k_i^{\{t\}} \\ \ell_i^{\{t\}}
\end{bmatrix} 
=
\begin{bmatrix}
0 \\ \textnormal{NB}(\mathfrak{g}_{l,i}^{\{m\}},[\q_0, \p_0])
\end{bmatrix};
\qquad 
%\]
%and 
%\[
\begin{bmatrix}
k_i^{\{N+t\}} \\ \ell_i^{\{N+t\}}
\end{bmatrix} 
=
\begin{bmatrix}
\textnormal{NB}(\mathfrak{g}_{k,i}^{\{m\}},[\q_0, \p_0]) \\
0
\end{bmatrix}, ~~ 1 \le t \le \nparts. 
\]
Hence we have for the non-empty NT trees $u=[u_1,\dots,u_r]_{\{m\}}$ and $v=[v_1,\dots,v_p]_{\{n\}}$
\begin{align*}
    1 \le s \le \nparts: \qquad & U \coloneqq  
    \begin{bmatrix}
    0 \\ \mathfrak{g}_\ell^{\{m\}} (u_1)
    \end{bmatrix} 
    \times \cdots \times
    \begin{bmatrix}
    0 \\ \mathfrak{g}_\ell^{\{m\}} (u_r)
    \end{bmatrix},  \quad
    V \coloneqq  
    \begin{bmatrix}
    0 \\ \mathfrak{g}_\ell^{\{n\}} (v_1)
    \end{bmatrix} 
    \times \cdots \times
    \begin{bmatrix}
    0 \\ \mathfrak{g}_\ell^{\{n\}} (v_r)
    \end{bmatrix}, \\
    \nparts + 1 \le s \le  2\nparts: \qquad & 
   U \coloneqq  \begin{bmatrix}
    \mathfrak{g}_k^{\{m\}} (u_1) \\ 0
    \end{bmatrix} \times \cdots \times
    \begin{bmatrix}
    \mathfrak{g}_k^{\{m\}} (u_r) \\ 0
    \end{bmatrix},  \quad
    V \coloneqq  \begin{bmatrix}
    \mathfrak{g}_k^{\{n\}} (v_1) \\ 0
    \end{bmatrix} \times \cdots \times
    \begin{bmatrix}
    \mathfrak{g}_k^{\{n\}} (v_r) \\ 0
    \end{bmatrix},   
\end{align*}
%for $m>1$.
and thus $U_i V_j \neq 0$ for $m \le \nparts, n>N$  or  $m>N, n \le \nparts$.  Hence
\eqref{eq.derivation.sympl} implies
$\mathbf{P}^{\{\ell,m\}}=0$ for $\ell \in \{1,\ldots,N\}$ and $m \in \{N+1,\ldots, 2\nparts\}$.
\end{proof}
\begin{remark}[Dimensions] 
\label{rem:dimensions}
Equation \eqref{eqn:symplecticness} implies $\hat{s}^{\{m\}} = s^{\{m\}}$ for all $m$.
\end{remark}

Consider the $\Bhat^{\{\ell\}}$ weights to be degrees of freedom, with $\Bhat^{\{\ell\}}$ regular. The symplecticness equation \eqref{eqn:symplecticness} can be solved for $\Ahat^{\{\ell,m\}}$ to obtain:
\begin{equation}
\label{part.GARK.sympl.cond}
\begin{split}
& \Ahat^{\{\ell,m\}} = \one^{\{\ell\}}\, \bhat^{\{m\}}\tr
-\B^{\{\ell\}-1}\, \mathbf{A}^{\{m,\ell\}}\tr\, \Bhat^{\{m\}}, \\
& \Rightarrow \quad \chat^{\{\ell,m\}} = \one^{\{\ell\}} - \B^{\{\ell\}-1}\,\mathbf{d}^{\{\ell,m\}}
 \quad \textnormal{with} \quad \mathbf{d}^{\{\ell,m\}} \coloneqq \mathbf{A}^{\{m,\ell\}}\tr\, \bhat^{\{m\}}, \\
&  \Leftrightarrow \quad \mathbf{d}^{\{\ell,m\}}= \b^{\{\ell\}} \times( \one^{\{\ell\}} - \chat^{\{\ell,m\}}).
\end{split}
\end{equation}
We note that if the weights are equal, $\b =\bhat$, then $\Ahat^{\{\ell,m\}}$ is fixed via \eqref{part.GARK.sympl.cond} by the choice of the base methods $(\mathbf{A}^{\{m,\ell\}},\b ^{\{\ell\}})$:
\begin{equation}
\label{part.GARK.sympl.cond2}
\begin{split}
& \Ahat^{\{\ell,m\}} =  \one^{\{\ell\}}\, \b^{\{m\}}\tr
-\mathbf{\B}^{\{\ell\}-1}\, \mathbf{A}^{\{m,\ell\}}\tr\, \B^{\{m\}}.
%\quad \chat^{\{\ell,m\}} = \one^{\{\ell\}} -\B^{\{\ell\}-1}\,\mathbf{d}^{\{\ell,m\}}, \\
%& \mathbf{d}^{\{\ell,m\}} = \mathbf{A}^{\{m,\ell\}}\tr\, \b^{\{m\}} = \b^{\{\ell\}} \times( \one^{\{\ell\}} - \chat^{\{\ell,m\}}).
\end{split}
\end{equation}

\begin{remark}
The Runge-Kutta $D(1)$ simplifying assumption extends to GARK $D^{\{m,\ell\}}(1)$ simplifying assumption \cite{Tanner_2018_PhD} , which for our method reads:
\begin{equation}
\label{part.GARK.sympl.extended-D1}
\begin{split}
\mathbf{A}^{\{m,\ell\}}\tr\, \bhat^{\{m\}} &= \b ^{\{\ell\}} \times \bigl( \one - \chat^{\{\ell, m\}} \bigr), \\
\Ahat^{\{m,\ell\}}\tr\, \b^{\{m\}} &= \bhat^{\{\ell\}} \times \bigl( \one - \mathbf{c}^{\{\ell, m\}} \bigr).
\end{split}
\end{equation}
If the symplecticness condition \eqref{eqn:symplecticness} holds then both \eqref{part.GARK.sympl.extended-D1} equations are fulfilled. This can be seen by multiplying \eqref{eqn:symplecticness} with a vector of ones from the left and from the right.
\end{remark}

%\begin{remark}
%For $T_2=\ldots=T_N=0$ the  Hamiltonian partitioning \eqref{eqn:split-TV} reduces to the potential splitting case~\eqref{eqn:split-V}, and method \eqref{eqn:Partitioned-Butcher-tableau} reduces to ~\eqref{eqn:Partitioned-Butcher-tableau-PS}.
%\end{remark}
%\color{black}

%%%%%%%%%%%%%%%%%%%%%%%
\subsection{GARK discrete adjoints}

GARK discrete adjoints were developed in \cite{Sandu_2021_GARK-adjoint}. As in the case of standard Runge-Kutta methods \cite{Sandu_2006_dadjRK}, if all the weights are nonzero, $b_i^{\{q\}}\neq 0$, one can reformulate the discrete GARK adjoint as another GARK method
to advance (in reverse time) the adjoint variables $\lambda_{n}$:
\begin{subequations}
\label{eq:GARK_costate}
\begin{align}
\label{eq:costate_step_GARK} 
\lambda_{n} &= \lambda_{n+1} + h \sum_{q=1}^{\nparts} \sum_{j=1}^{s^{\{q\}}} \bar{b}_j^{\{q\}} \, \ell_{n,j}^{\{q\}}, \\
\label{eq:costate_stages_GARK} 
\Lambda_{n,i}^{\{q\}} &= \lambda_{n+1} + h \sum_{m=1}^{\nparts} \sum_{j=1}^{s^{\{m\}}} \bar{a}_{i,j}^{\{q,m\}} \, \ell_{n,j}^{\{m\}}, 
\quad i = s^{\{q\}}, \dots, 1, \\
\label{eq:costate_stages_ell} 
\ell_{n,i}^{\{q\}} &= \fun^{\{q\}}_{\yy}\tr(Y_{n,i}^{\{q\}}) \cdot \Lambda_{n,i}^{\{q\}},\\
\label{eq:GARK_costate_coefficients}
\bar{b}_i^{\{q\}}&=b_i^{\{q\}}, \qquad \bar{a}_{i,j}^{\{q,m\}}=\frac{b_j^{\{m\}}\,a_{j,i}^{\{m,q\}}}{b_i^{\{q\}}}.
\end{align}
\end{subequations}
Reverting the time $h \to -h$ the method \eqref{eq:GARK_costate} reads:
\begin{equation}
\label{eq:GARK_costate.fwd}
\begin{split}
\Lambda_{n,i}^{\{q\}} 
%&= \lambda_{n} + h \sum_{m=1}^{\nparts} \sum_{j=1}^{s^{\{m\}}} \bar{b}_j^{\{m\}} \, \ell_{n,j}^{\{m\}} - h \sum_{m=1}^{\nparts} \sum_{j=1}^{s^{\{m\}}} \bar{a}_{i,j}^{\{q,m\}} \, \ell_{n,j}^{\{m\}} \\
&= \lambda_{n} + h \sum_{m=1}^{\nparts} \sum_{j=1}^{s^{\{m\}}} \ahat_{i,j}^{\{q,m\}} \, \ell_{n,j}^{\{m\}}, \\
\lambda_{n+1} &= \lambda_{n} + h \sum_{q=1}^{\nparts} \sum_{j=1}^{s^{\{q\}}} \widehat{b}_j^{\{q\}} \, \ell_{n,j}^{\{q\}}, \\
\ahat_{i,j}^{\{q,m\}} &= \bar{b}_j^{\{m\}} - \bar{a}_{i,j}^{\{m,q\}} 
= b_j^{\{m\}} - \frac{b_j^{\{m\}}\,a_{j,i}^{\{m,q\}}}{b_i^{\{q\}}},
\quad \widehat{b}_j^{\{q\}} = \bar{b}_j^{\{q\}} = b_j^{\{q\}},
\end{split}
\end{equation}
which is called the formal discrete adjoint GARK method. In matrix notation the coefficients read
\begin{equation}
\label{GARK.symplectic-conjugate}
   \Ahat^{\{q,m\}} = \one^{\{q\}}\, \b^{\{m\}}\tr
- \B^{\{q\}-1}\, \mathbf{A}^{\{m,q\}}\tr\, \B^{\{m\}},
\quad \bhat^{\{q\}} = \b^{\{q\}},
\end{equation}
and is equivalent to the symplecticness condition~\eqref{part.GARK.sympl.cond2} for partitioned GARK schemes. 
The matrix $\Ahat^{\{q,m\}}$ given by \eqref{GARK.symplectic-conjugate}
is called the symplectic conjugate of ${\A}^{\{q,m\}}$. If ${\A}^{\{q,m\}}=\Ahat^{\{q,m\}}$ holds for all $q,m$ the GARK method is called {\it self-adjoint.} We have the following result.
\begin{lemma}
\label{lemma.adj.sym}
Symplecticity and self-adjointness of a  GARK scheme $(\b^{\{m\}}, \A^{\{m,n\}})$ with $b^{\{m\}}_i \ne 0$ are equivalent properties. 
\end{lemma}
\begin{proof}
For nonzero weights the symplecticness condition \eqref{eqn:GARK-symplectic-condition} is equivalent to
\begin{equation*}
\mathbf{P}^{\{q,m\}}  
= \Zero_{s^{\{q\}} \times s^{\{m\}}}
\quad \Leftrightarrow \quad
\A^{\{q,m\}}  =  \one^{\{q\}}\, \b^{\{m\}}\tr
-\mathbf{\B}^{\{q\}-1}\, \mathbf{A}^{\{m,q\}}\tr\, \B^{\{m\}},
\end{equation*}
and therefore is equivalent to self-adjointness using \eqref{GARK.symplectic-conjugate}.
\end{proof}

It was shown in \cite{Sandu_2021_GARK-adjoint} that the order of the discrete adjoint method \eqref{eq:GARK_costate}, \eqref{eq:GARK_costate.fwd} coincides with the order of the base GARK scheme when computing solution derivatives, i.e., $\lambda_n = (d\Phi/d \y_n)\tr$ for some functional $\Phi$ defined on the solution $\{\y_i\}$. This is true if in \eqref{eq:costate_stages_ell} the Jacobians $\fun^{\{q\}}_{\yy}\tr(Y_{n,i}^{\{q\}})$ are evaluated at the forward GARK stages. However, this is not true for the formal discrete adjoint GARK \eqref{eq:GARK_costate.fwd} regarded as a general GARK integration scheme.

\begin{remark}
\ifreport
Sanz-Serna \cite{SanzSerna_2016_adjointSymplectic} considers Lagrangian problems:
\begin{equation*}
\int \mathcal{L}(q,\dot{q},t) \to \min, \quad \mathcal{L}(q,\dot{q},t) = T(q,\dot{q},t) - U(q, t),
\end{equation*}
where the control co-states are equivalent to momenta:
\begin{equation*}
p = \nabla_{\dot{q}}\,\mathcal{L}(q,\dot{q},t).
\end{equation*}
\fi
Zanna~\cite{zanna20dvm} has shown that integrating the states $Q^{\{m\}}$ with a GARK scheme, and the co-states $P^{\{m\}}$ with the discrete adjoint of the GARK method~\eqref{eq:GARK_costate} results in a symplectic numerical method. This result generalizes the results of Sanz-Serna \cite{SanzSerna_2016_adjointSymplectic} for RK methods to GARK methods.
\end{remark}

%%%%%%%%%%%%%%%%%%%%%%%%%%%%
\subsection{Order conditions}
%%%%%%%%%%%%%%%%%%%%%%%%%%%%
In this section we consider a partitioned GARK scheme \eqref{eq:sympl.part.gark}--\eqref{eqn:Partitioned-Butcher-tableau} satisfying the symplecticness condition \eqref{eqn:symplecticness}, and study the order conditions when applied to solve a partitioned system of the form \eqref{eqn:split-TV}.

\begin{remark}
The choice of setting the coefficients $\A^{\{m,n\}}$ and $\A^{\{\nparts+m,\nparts+n\}}$ for $1 \le m,n \le \nparts$ to zero (or any other value) in \eqref{eqn:Partitioned-Butcher-tableau} follows the fact that they do not contribute to the order conditions since the associated elementary differentials  are identically equal to zero.
To see this, consider general N-trees of the form 
\[
v = \Big[[t_1]_{\{m_1\}} \dots [t_M]_{\{m_M\}},[r_1]_{\{\ell_1 + \nparts\}} \dots [r_L]_{\{\ell_L + \nparts\}} \Big]_{\{n\}}, \quad 1 \le m_i,\ell_i \le \nparts.
\]
The associated order condition involves $\A^{\{n,m_i\}}$ and $\A^{\{n,\ell_i + \nparts\}}$, and 
the corresponding elementary differential is:
\[
\begin{split}
F(v) = \begin{cases}
\frac{\partial^{M+L} \mathbf{f}^{\{n\}}(\p)}{\partial \q^M\, \partial \p^L} \big( \cdots \big) 
   %= \frac{\partial^{L+M} \begin{bmatrix} T_{\p}^{\{m\}}(\p) \\ 0 \end{bmatrix}}{\partial \q^L\, \partial \p^M} \big( \dots \big) 
   \equiv \scriptstyle 0 & \textnormal{for} ~~ 1 \le n \le \nparts ~~\textnormal{and} ~~ M \ge 1, \\[3pt]
\frac{\partial^{M+L} \mathbf{f}^{\{n\}}(\q)}{\partial \q^M\, \partial \p^L} \big( \cdots \big) 
    \equiv \scriptstyle 0 & \textnormal{for} ~~ \nparts+ 1 \le n \le 2\,\nparts ~~\textnormal{and} ~~ L \ge 1.
\end{cases}
\end{split}
\]
The only non-zero elementary differentials correspond to N-trees where any $\q$-node has only $\p$-children, and vice-versa.
\end{remark}

The order conditions (up to order four) for a symplectic partitioned GARK scheme \eqref{eq:sympl.part.gark}--\eqref{eqn:Partitioned-Butcher-tableau} are obtained directly from Theorem~\ref{th.sympl.order.cond}:
\begin{subequations}
\label{eq.ordercond.3.2}
\begin{align}
\label{eq.ordercond.3.2.1a}
\b ^{\{m\}}\tr \cdot \one^{{\{m\}}} = 1,  & \quad \forall\; m, &  (\textnormal{order}~ 1) \\
\label{eq.ordercond.3.2.1b}
\bhat^{\{m\}}\tr \cdot \one^{{\{m\}}} = 1,  & \quad \forall\; m, &  (\textnormal{order}~ 1) \\
\label{eq.ordercond.3.2.2}
 \bhat^{\{m\}}\tr \cdot \c^{\{m,\ell\}}  =  \sfrac{1}{2}, &\quad  \forall\; m,\ell,  &  (\textnormal{order}~ 2) \\
 \label{eq.ordercond.3.2.3a}
\b^{\{m\}}\tr \cdot \left(  \chat^{\{m,\ell\}} \times \chat^{\{m,s\}} \right)
 =  \sfrac{1}{3}, & \quad \forall\; m,\; \forall\; \ell \leq s, & (\textnormal{order}~ 3) \\
 \label{eq.ordercond.3.2.3b}
\bhat^{\{m\}}\tr \cdot \left(  \mathbf{c}^{\{m,\ell\}} \times \mathbf{c}^{\{m,s\}} \right)
 =  \sfrac{1}{3}, & \quad \forall\; m,\; \forall\; \ell \leq s, & (\textnormal{order}~ 3) \\
 \label{eq.ordercond.3.2.4a}
\b^{\{m\}}\tr \!\cdot\! \left( \chat^{\{m,\ell\}} \times \chat^{\{m,s\}} \times \chat^{\{m,t\}} \right)
 =  \sfrac{1}{4}, & \quad \forall\; m,\; \forall\; \ell \le s \le t, & (\textnormal{order}~ 4) \\
 \label{eq.ordercond.3.2.4b}
\bhat^{\{m\}}\tr \!\cdot\! \left( \mathbf{c}^{\{m,\ell\}} \times \mathbf{c}^{\{m,s\}} \times \mathbf{c}^{\{m,t\}} \right)
 =  \sfrac{1}{4}, & \quad \forall\; m,\; \forall\; \ell \le s \le t, & (\textnormal{order}~ 4) \\
% \label{eq.ordercond.3.2.4c}
%%
%\b^{\{m\}}\tr \cdot \left(\chat^{\{m,t\}} \times \Ahat^{\{m,s\}}  \cdot \mathbf{c}^{\{s,\ell\}} \right) 
% =  \sfrac{1}{8}, & \quad \forall\; m,t,s,\ell\ge 1, & (\textnormal{order}~ 4) \\
 \label{eq.ordercond.3.2.4d}
 \bhat^{\{m\}}\tr \!\cdot\! \left(\mathbf{c}^{\{m,\ell\}} \times \mathbf{A}^{\{m,s\}}  \cdot \chat^{\{s,t\}} \right) 
 =  \sfrac{1}{8}, & \quad \forall\; m,\ell,s,t. & (\textnormal{order}~ 4) 
\end{align}
\end{subequations}

\begin{remark}
Using \eqref{part.GARK.sympl.cond} these order conditions can be rewritten in terms of
$\mathbf{d}^{\{m,\ell\}}$  by substituting $\chat^{\{m,\ell\}} = \one^{\{m\}} - \B^{\{m\}-1}\,\mathbf{d}^{\{m,\ell\}}$:
\begin{align*}
    \eqref{eq.ordercond.3.2.3a} \quad \Leftrightarrow \quad &
    \b ^{\{\ell\}}\tr \cdot \left(  \B^{\{\ell\}-1}\mathbf{d}^{\{\ell,s\}} \times \B^{\{\ell\}-1} \mathbf{d}^{\{\ell,t\}} \right)
 =  \sfrac{1}{3}, \\
   \eqref{eq.ordercond.3.2.4a} \quad \Leftrightarrow \quad &
\b ^{\{\ell\}}\tr \cdot \left(  \B^{\{\ell\}-1}\mathbf{d}^{\{\ell,s\}} \times \B^{\{\ell\}-1} \mathbf{d}^{\{\ell,t\}} \times \B^{\{\ell\}-1} \mathbf{d}^{\{\ell,s\}} \right)
 =  \sfrac{1}{4}, \\  
 \eqref{eq.ordercond.3.2.4d} \quad \Leftrightarrow \quad & 
  \bhat^{\{m\}}\tr \cdot \left(  
  \c^{\{m,s\}} \times \A^{\{m,t\}} \B^{\{t\}-1}
  \mathbf{d}^{\{t,\ell\}}  \right) = \sfrac{5}{24}.
\end{align*}
For $\nparts=1$ and $\b^{\{1\}}=\bhat^{\{1\}}$, these order conditions coincide with the order conditions given by Hager ~\cite[Table 1]{Hager_2000_RKadj} (neglecting the superscripts for simplicity):
\begin{align*}
    \sum_i b_i & = 1, &
    \sum_i d_i & = \sfrac{1}{2}, &
    \sum_i \sfrac{d_i^2}{b_i} & = \sfrac{1}{3}, &
    \sum_i b_i c_i^2 & = \sfrac{1}{3}, \\
    \sum_i \sfrac{d_i^3}{b_i^2} & = \sfrac{1}{4}, &
    \sum_i b_i c_i^3 & = \sfrac{1}{4}, &
%    \sum_{i,j} \frac{d_i a_{j,i}b_j c_j }{b_i} & = \frac{5}{24}, &
    \sum_{i,j} \frac{b_i c_i a_{i,j} d_j }{b_j} & = \sfrac{5}{24}.
\end{align*}
Note that the remaining order conditions in Hager (the first order three condition and the second, third, fourth, fifth and eighth order four conditions in ~\cite[Table 1]{Hager_2000_RKadj}) are redundant due to the symplecticity of the partitioned scheme.
\end{remark}

\begin{lemma}[The same weights case]
\label{lem:order4.partGARK}
Let us now assume that $s^{\{\ell\}}={\hat s}^{\{\ell\}}$ and  $\b^{\{\ell\}}={\bhat}^{\{\ell\}}$ holds for all $\ell=1,2,\ldots,\nparts$. 
If the base GARK scheme $(\b^{\{m\}},\A^{\{m,\ell\}})$ has order four, then the symplectic partitioned scheme \eqref{eq:sympl.part.gark}--\eqref{eqn:Partitioned-Butcher-tableau}  has at least order two. In addition, we have:
\begin{itemize}
    \item The partitioned GARK scheme has order three if condition~\eqref{eq.ordercond.3.2.3a} holds.
    \item In addition, the partitioned GARK scheme has order four, if conditions~\eqref{eq.ordercond.3.2.4a} and~\eqref{eq.ordercond.3.2.4d} holds. These are automatically fulfilled for symmetric schemes.
\end{itemize}
\end{lemma}

\begin{example}
\label{ex:Verlet}
The Verlet algorithm (the two-stage Lobatto IIIA--IIIB pair of order 2) is of the form \eqref{eqn:Partitioned-Butcher-tableau}  with $\nparts=1$:
\begin{equation}
\label{eq:Verlet}
\renewcommand{\arraystretch}{1.5}
\begin{array}{c|c} \chat & \Ahat \\ \hline & \bhat\tr \end{array}  
~=~
\raisebox{9pt}{$
\begin{array}{c|cc}
\scriptstyle \scriptstyle 0 &\scriptstyle  \scriptstyle 0 &\scriptstyle  0 \\ \scriptstyle 1 & \frac{1}{2} & \frac{1}{2} \\ \hline & \frac{1}{2} & \frac{1}{2}
\end{array}$},
\qquad
\begin{array}{c|c} \c & \A \\ \hline & \b\tr \end{array} 
~=~
\raisebox{9pt}{$
\begin{array}{c|cc}
\frac{1}{2} & \frac{1}{2} &\scriptstyle  0 \\ \frac{1}{2} & \frac{1}{2} &\scriptstyle  0 \\ \hline & \frac{1}{2} & \frac{1}{2}
\end{array}$}.
\end{equation}
Similarly, the three-stage Lobatto IIIA--IIIB pair of order 4 has the form \eqref{eqn:Partitioned-Butcher-tableau}
\begin{equation}
\label{eq:Lobatto}
\renewcommand{\arraystretch}{1.5}
\begin{array}{c|c} \chat & \Ahat \\ \hline & \bhat\tr \end{array}  
~=~
\raisebox{18pt}{$
\begin{array}{c|rrr}
\scriptstyle \scriptstyle 0 &\scriptstyle  \scriptstyle 0 &\scriptstyle   0  &\scriptstyle  0  \\
\frac{1}{2} &  \frac{5}{24}  & \frac{1}{3}   & -\frac{1}{24}\\
\scriptstyle 1 & \frac{1}{6}  & \frac{2}{3}  & \frac{1}{6} \\ \hline
   &  \frac{1}{6}  & \frac{2}{3}  & \frac{1}{6}
 \end{array}$},
\qquad
\begin{array}{c|c} \c & \A \\ \hline & \b\tr \end{array} 
~=~
\raisebox{18pt}{$
\begin{array}{c|rrr}  
  \scriptstyle 0  & \frac{1}{6}  & -\frac{1}{6}  & \scriptstyle 0 \\
   \frac{1}{2}  & \frac{1}{6} &  \frac{1}{3}  & \scriptstyle 0 \\
  \scriptstyle  1  & \frac{1}{6}  & \frac{5}{6} &\scriptstyle   0\\ \hline
       &  \frac{1}{6}  & \frac{2}{3}  & \frac{1}{6}
\end{array}$}.
\end{equation}
%\end{remark}
\end{example}

%%%%%%%%%%%%%%%%%%%%%%%%%%%%
\subsection{Symmetry and time-reversibility}
\label{sec:symmetry-and-reversibility}
%%%%%%%%%%%%%%%%%%%%%%%%%%%%
The symmetry and time-reversible results for GARK schemes can be adapted to partitioned GARK schemes as follows.

\begin{theorem}[Symmetric partitioned GARK schemes]
The partitioned GARK scheme \eqref{eq:sympl.part.gark} is symmetric \eqref{cond.symm.comp} if the following conditions hold:
%
%\begin{subequations}
%\label{symmetry.part}
%\begin{eqnarray*}
%b_j^{\{q\}} & = & b_{s+1-j}^{\{q\}}, \qquad j=1,\ldots,\hat s^{\{q\}}, \\
%a_{i,j}^{\{m,q\}} & = &  b_{j}^{\{q\}} - a_{s^{\{m\}}+1-i,\hat s^{\{q\}}+1-j}^{\{m,q\}}, \qquad
%i=1,\ldots, s^{\{m\}}, \, j=1,\ldots,\hat s^{\{q\}}, \\
%\widehat{b}_j^{\{q\}} & = & \widehat{b}_{s+1-j}^{\{q\}}, \qquad j=1,\ldots,s^{\{q\}}, \\
%\ahat_{i,j}^{\{m,q\}} & = &  \widehat{b}_{j}^{\{q\}} - \ahat_{\hat s^{\{m\}}+1-i,s^{\{q\}}+1-j}^{\{m,q\}}, \qquad
%i=1,\ldots, \widehat s^{\{m\}}, \, j=1,\ldots,s^{\{q\}}. 
%\end{eqnarray*}
%\end{subequations}
%
\begin{subequations}
\label{eqn:GARK-symmetry-matrix}
\begin{eqnarray}
\label{eqn:GARK-symmetry-matrix.a}
\b^{\{m\}} & = & \underline{\b}^{\{m\}},  
%= \mathcal{P}\,\b^{\{m\}}, 
\qquad \bhat^{\{m\}}  = \underline{\bhat}^{\{m\}} , 
%= \mathcal{P}\,\hat{\b}^{\{m\}}, 
\\
\label{eqn:GARK-symmetry-matrix.b}
\A^{\{\ell,m\}} & = & \underline{\A}^{\{\ell,m\}}, %= \one\,\b^{\{m\}}\tr - \mathcal{P}^{\{\ell\}}\,\A^{\{\ell,m\}}\,\mathcal{P}^{\{m\}}, 
\quad
\Ahat^{\{\ell,m\}} = \underline{\Ahat}^{\{\ell,m\}} 
%= \one\,\hat{\b}^{\{m\}}\tr - \mathcal{P}\,\hat{\A}^{\{\ell,m\}}\,\mathcal{P}\,. 
\end{eqnarray}
\end{subequations}
In other words, a partitioned GARK scheme is symmetric, iff both the base GARK scheme and its discrete adjoint scheme are symmetric.
\end{theorem}

\begin{proof}
The time-reversed tableau \eqref{eqn:GARK-time-reversed-tableau} 
of the partitioned method \eqref{eqn:Partitioned-Butcher-tableau} is:
\begin{equation}
\label{eqn:Partitioned-Butcher-tableau-reversed}
\renewcommand{\arraystretch}{1.5}
\begin{array}{c}
\underline{\A}_{\textsc{gark}}
\\
\Xhline{2\arrayrulewidth}
\underline{\b}_{\textsc{gark}}\tran
\end{array}
~=~
\raisebox{10pt}{$
\begin{array}{c|c}
%\one\,\b \tr  & \one\,\bhat\tr - \mathcal{P} \mathbf{\hat A } \mathcal{P} \\ \hline
%\one\,\b \tr  - \mathcal{P} \mathbf{A} \mathcal{P}  & \one\,\bhat\tr
%\\
\Zero_{s \times s}  & \one\,\bhat\tr - \mathcal{P}\, \Ahat \,\mathcal{P} \\ \hline
\one\,\b \tr  - \mathcal{P}\, \mathbf{A}\, \mathcal{P}  & \Zero_{s \times s}
\\
\Xhline{2\arrayrulewidth}
(\mathcal{P}\, \b )\tr & (\mathcal{P}\,  \bhat)\tr
\end{array}
$}
~=~
\raisebox{10pt}{$
\begin{array}{c|c}
\Zero_{s \times s}  & \underline{\Ahat} \\ \hline
\underline{\mathbf{A}} & \Zero_{s \times s}
\\
\Xhline{2\arrayrulewidth}
\underline{\b}\tr & \underline{\bhat}\tr
\end{array}
$}.
\end{equation}
The diagonal blocks in the time-reversed scheme \eqref{eqn:Partitioned-Butcher-tableau-reversed} are set to $0$ as they do not have any effect in the computation.
\end{proof}

\begin{lemma}[Symmetry of partitioned GARK schemes]
\label{lem:symmetry-PGARK}
Consider a GARK scheme $(\b^{\{m\}},\A^{\{m,\ell\}})$ that is symmetric \eqref{eqn:GARK-symmetry}. Then its discrete adjoint scheme $(\b^{\{m\}},\Ahat^{\{m,\ell\}})$ given by \eqref{GARK.symplectic-conjugate} is also symmetric, and so is the partitioned GARK scheme $(\b^{\{m\}},\A^{\{m,\ell\}},\Ahat^{\{m,\ell\}})$ given by \eqref{eq:sympl.part.gark}--\eqref{eqn:Partitioned-Butcher-tableau} with $\bhat^{\{m\}} = \b^{\{m\}}$.
\end{lemma}

\begin{proof}
We only have to show that $\A^{\{\ell,m\}}  =  \underline{\A}^{\{\ell,m\}}$ implies
$\Ahat^{\{\ell,m\}} =  \underline{\Ahat}^{\{\ell,m\}}$. From \eqref{eqn:GARK-time-reversed-1.b} the symmetry of the GARK scheme is:
\[
\A^{\{\ell,m\}} = \underline{\A}^{\{\ell,m\}} = \one^{\{\ell\}}\,\b^{\{m\}}\tr - \mathcal{P}^{\{\ell\}}\,\A^{\{\ell,m\}}\,\mathcal{P}^{\{m\}}.
\]
From \eqref{eqn:GARK-time-reversed-1.b}, \eqref{GARK.symplectic-conjugate}, and the symmetry of the GARK scheme assumption, the time reversed discrete adjoint matrix is:
\begin{equation*}
\begin{split}
\underline{\Ahat}^{\{\ell,m\}} &= \one^{\{\ell\}}\,\b^{\{m\}}\tr - \mathcal{P}^{\{\ell\}}\,\Ahat^{\{\ell,m\}}\,\mathcal{P}^{\{m\}} \\
&= \one^{\{\ell\}}\,\b^{\{m\}}\tr - \mathcal{P}^{\{\ell\}}\,\left( \one^{\{\ell\}}\, \b^{\{m\}}\tr
- \B^{\{\ell\}-1}\, \mathbf{A}^{\{m,\ell\}}\tr\, \B^{\{m\}}\right)\,\mathcal{P}^{\{m\}} \\
%&= \one^{\{\ell\}}\,\b^{\{m\}}\tr - \one^{\{\ell\}}\, \underline{\b}^{\{m\}}\tr
%+ \underline{\B}^{\{\ell\}-1}\, \mathbf{A}^{\{m,\ell\}}\tr\, \underline{\B}^{\{m\}} \\
&= \B^{\{\ell\}-1}\, \mathbf{A}^{\{m,\ell\}}\tr\, \B^{\{m\}} \\
%&= \B^{\{\ell\}-1}\, \left(  \one^{\{m\}}\,\b^{\{\ell\}}\tr - \mathcal{P}^{\{m\}}\,\A^{\{m,\ell\}}\,\mathcal{P}^{\{\ell\}} \right)\tr\, \B^{\{m\}} \\
&= \B^{\{\ell\}-1}\, \left(  \b^{\{\ell\}}\,\one^{\{m\}}\tr - \mathcal{P}^{\{\ell\}}\,\A^{\{m,\ell\}}\tr\,\mathcal{P}^{\{m\}} \right)\, \B^{\{m\}} \\
&= \one^{\{\ell\}}\,\b^{\{m\}}\tr - \B^{\{\ell\}-1}\,\A^{\{m,\ell\}}\tr\,\B^{\{m\}} \\
&= \Ahat^{\{\ell,m\}}.
\end{split}
 \end{equation*} 
%
%\begin{equation*}
%  \Ahat^{\{\ell,m\}} =  \underline{\Ahat}^{\{\ell,m\}}  \Leftrightarrow
% \Ahat^{\{\ell,m\}} = \one^{\{\ell\}}\, \b^{\{m\}}\tr - 
%\mathcal{P}^{\{\ell\}}\, \Ahat^{\{\ell,m\}}\, \mathcal{P}^{\{m\}},
% \end{equation*} 
% %
% which is equivalent to
% %
%  \begin{equation*}
%  \one \b^{\{m\}}\tr - (\B^{\{\ell\}})^{-1} \A^{\{m,\ell\}}\tr \B^{\{m\}} = 
% \one \b^{\{m\}}\tr - 
%  \mathcal{P} \left( \one \b^{\{m\}}\tr - (\B^{\{\ell\}})^{-1} \A^{\{m,\ell\}}\tr \B^{\{m\}} \right) \mathcal{P}
%\end{equation*}
%%
%as the partitioned GARK scheme is symplectic.
%Multiplying with $\B^{\{\ell\}}$ and $(\B^{\{m\}})^{-1}$
%from left and right, respectively, yields
%%
%\begin{equation*}
% \A^{\{m,\ell\}\tr}  = 
% \B^{\{\ell\}}  \mathcal{P} \left( \one \b^{\{m\}}\tr - (\B^{\{\ell\}})^{-1} \A^{\{m,\ell\}}\tr \B^{\{m\}} \right) \mathcal{P}  (\B^{\{m\}})^{-1}
%\end{equation*}
%With $\B^{\{\ell\}}  \mathcal{P} = 
%\mathcal{P} \B^{\{\ell\}}$ and 
%$\mathcal{P}  (\B^{\{m\}})^{-1} = (\B^{\{m\}})^{-1} \mathcal{P} $ we get
%\begin{equation*}
%    \A^{\{m,\ell\}}\tr   = 
% \mathcal{P}  \B^{\{\ell\}}  \left( \one \b^{\{m\}}\tr - (\B^{\{\ell\}})^{-1} \A^{\{m,\ell\}}\tr \B^{\{m\}} \right)  (\B^{\{m\}})^{-1} \mathcal{P}  \Leftrightarrow
% \Ahat^{\{m,\ell\}} =  \underline{\Ahat}^{\{m,\ell\}}
%\end{equation*}
\end{proof}

%%%%%%%%%%%%%%%%%%%%%%%%%%%%%%%%%%%%%%%%%%%%%%%%%%%%%%%
\subsubsection{Construction of symmetric and symplectic methods starting from a symmetric scheme}
\label{sec.3.4.1}
%%%%%%%%%%%%%%%%%%%%%%%%%%%%%%%%%%%%%%%%%%%%%%%%%%%%%%%

Lemma \ref{lem:symmetry-PGARK} provides an easy way to construct partitioned GARK schemes of order four, which are both symmetric and symplectic. One starts with a symmetric GARK scheme $(\b^{\{m\}},\A^{\{m,\ell\}})$ and constructs the symplectic partitioned GARK scheme  $(\b^{\{m\}},\A^{\{m,\ell\}},\Ahat^{\{m,\ell\}})$. This scheme is symmetric by Lemma \eqref{lem:symmetry-PGARK}. If condition~\eqref{eq.ordercond.3.2.3a} holds the partitioned scheme has order three, and therefore order four is ensured by symmetry.

\begin{remark}
Consider a partitioned GARK scheme that is both symplectic and symmetric. The symmetry condition \eqref{eqn:GARK-symmetry-matrix.b} together with the symplecticness condition \eqref{eqn:symplecticness} give
\begin{equation}
\label{eq.cond.symmetry}
\begin{split}
\A^{\{\ell,m\}}  &=  \one^{\{\ell\}}\,\b^{\{m\}}\tr - \mathcal{P}^{\{\ell\}}\,\A^{\{\ell,m\}}\,\mathcal{P}^{\{m\}} =  \one^{\{\ell\}}\,\bhat^{\{m\}}\tr - \B^{\{\ell\}-1}\, \Ahat^{\{m,\ell\}}\tr\,\Bhat^{\{m\}}, \quad \\
\Ahat^{\{\ell,m\}} &= \one^{\{\ell\}}\,\bhat^{\{m\}}\tr - \mathcal{P}^{\{\ell\}}\,\Ahat^{\{\ell,m\}}\,\mathcal{P}^{\{m\}} = \one^{\{\ell\}}\, \b^{\{m\}}\tr
-\Bhat^{\{\ell\}-1}\, \mathbf{A}^{\{m,\ell\}}\tr\, \B^{\{m\}}\,. 
\end{split}
\end{equation}
and condition~\eqref{cond.merge} reads 
\begin{equation*}
\begin{split}
 \Ahat^{\{m,\ell\}}\tr\,\Bhat^{\{m\}} &=  \B^{\{\ell\}}\,\mathcal{P}^{\{\ell\}}\,\A^{\{\ell,m\}}\,\mathcal{P}^{\{m\}} + \b ^{\{\ell\}}\,(\bhat^{\{m\}}\tr -\b^{\{m\}}\tr ), \quad \\
 \mathbf{A}^{\{m,\ell\}}\tr\, \B^{\{m\}} &= \Bhat^{\{\ell\}}\,  \mathcal{P}^{\{\ell\}}\,\Ahat^{\{\ell,m\}}\,\mathcal{P}^{\{m\}} + \bhat^{\{\ell\}}\, (\b^{\{m\}}\tr - \bhat^{\{m\}}\tr).
\end{split}
\end{equation*}
The last terms vanish when $\b^{\{m\}}=\bhat^{\{m\}}$.
\end{remark}

%%%%%%%%%%%%%%%%%%%%%%%%%%%%%%%%%%%%%%%%%%%%%%%%%%%%%%%
\subsubsection{Construction of symmetric and symplectic methods starting from a symplectic scheme}
%%%%%%%%%%%%%%%%%%%%%%%%%%%%%%%%%%%%%%%%%%%%%%%%%%%%%%%
\label{sec3.4.2}
Remark~\ref{remark.sympl.tr} can also be applied to the time-reversed scheme~\eqref{eqn:GARK-time-reversed-tableau}, i.e., it is symplectic, iff the underlying partitioned GARK scheme is symplectic and all weights are palyndromic:
$\b^{\{m\}}={\underline{\b}}^{\{m\}}$ and $\bhat^{\{m\}}={\underline{\bhat}}^{\{m\}}$ for all $m$. Similar to Theorem~\ref{theorem.construction.symm.sympl.} we have the following result.

\begin{theorem}
\label{theorem.construction.symm.sympl2.}
Consider a partitioned GARK scheme \eqref{eqn:Partitioned-Butcher-tableau}
%
%\begin{align*}
%\begin{array}{c}
%\underline{\A}_{\textsc{gark}}
%\\
%\Xhline{2\arrayrulewidth}
%\underline{\b}_{\textsc{gark}}\tran
%\end{array}
% & =
%\begin{array}{c|c}
%\Zero_{s \times s}  & \Ahat  \\ \hline
%\mathbf{A}   & \Zero_{s \times s}
%\\
%\Xhline{2\arrayrulewidth}
%\b \tr & \bhat\tr
%\end{array}
%\end{align*}
that is symplectic and has palindromic weights,
$\b^{\{m\}}={\underline{\b}}^{\{m\}}$ and $\bhat^{\{m\}}={\underline{\bhat}}^{\{m\}}$ for all $m$. 
Applying one step with the partitioned GARK scheme, followed by one step with its time-reversed partitioned GARK scheme \eqref{eqn:Partitioned-Butcher-tableau-reversed}, 
%
%\begin{align*}
%\renewcommand{\arraystretch}{1.3}
%\underline{\A}_{\textsc{gark}}  =
%\begin{array}{c|c}
%\Zero_{s \times s}  & \one\,\bhat\tr - \mathcal{P} \Ahat \mathcal{P} \\ \hline
%\one\,\b \tr  - \mathcal{P} \mathbf{A} \mathcal{P}  & \Zero_{s \times s}
%\\
%\Xhline{2\arrayrulewidth}
%\underline{\b}\tr & \underline{\bhat}\tr
%\end{array},
%\end{align*}
%
defines a new GARK scheme with the Butcher tableau 
\begin{equation}
\label{eqn:symmetric-composition}
\renewcommand{\arraystretch}{1.3}
\begin{array}{c}
\widetilde{\A}_{\textsc{gark}}
\\
\Xhline{2\arrayrulewidth}
\widetilde{\b}_{\textsc{gark}}\tran
\end{array}
~~=~~
\raisebox{24pt}{$
\begin{array}{cc|cc}
\Zero_{s \times s}  & \frac{1}{2} \Ahat  & \Zero_{s \times s} & \Zero_{s \times s} \\
\frac{1}{2} \mathbf{A}   & \Zero_{s \times s} & \Zero_{s \times s} & \Zero_{s \times s} \\
\hline
\frac{1}{2} \one\b\tr & \frac{1}{2} \one\,\bhat\tr & \Zero_{s \times s}  & \frac{1}{2} \underline{\Ahat} \\
\frac{1}{2} \one\b\tr & \frac{1}{2} \one\,\bhat\tr & \frac{1}{2} \underline{\mathbf{A}} & \Zero_{s \times s} \\
\Xhline{2\arrayrulewidth}
\frac{1}{2} \b\tr  & \frac{1}{2} \bhat\tr & \frac{1}{2} \underline{\b}\tr  & \frac{1}{2} \underline{\bhat}\tr
\end{array}
$},
\end{equation}
which is both symmetric and symplectic. 
\end{theorem}

\begin{proof}
See proof of Theorem~\ref{theorem.construction.symm.sympl.}. Note that the coefficients are divided by two in order to recover the standard form over one step. 
\end{proof}
\begin{remark}
    Also note that \eqref{eqn:symmetric-composition} is not a partitioned GARK scheme of the form \eqref{eqn:Partitioned-Butcher-tableau}.
\end{remark}
\begin{leaveout}
\sandu{Add comment: Note that \eqref{eqn:symmetric-composition} is NOT a partitioned GARK scheme of the form \eqref{eqn:Partitioned-Butcher-tableau}.}
The composition of the schemes:
\begin{align*}
{\A}_{\textsc{gark}} & =
\begin{array}{c|c}
\Zero_{s \times s}  & \Ahat  \\ \hline
\mathbf{A}   & \Zero_{s \times s}
\\
\Xhline{2\arrayrulewidth}
\b \tr & \bhat\tr
\end{array}, \qquad
\underline{\A}_{\textsc{gark}}  =
\begin{array}{c|c}
\Zero_{s \times s}  & \one\,\bhat\tr - \mathcal{P} \Ahat \mathcal{P} \\ \hline
\one\,\b \tr  - \mathcal{P} \mathbf{A} \mathcal{P}  & \Zero_{s \times s}
\\
\Xhline{2\arrayrulewidth}
\underline{\b}\tr & \underline{\bhat}\tr
\end{array}
\end{align*}
leads to:
\[
\begin{array}{cc|cc}
\Zero_{s \times s}  & \Ahat  & \Zero_{s \times s} & \Zero_{s \times s} \\
\mathbf{A}   & \Zero_{s \times s} & \Zero_{s \times s} & \Zero_{s \times s} \\
\hline
\one\b\tr & \one\,\bhat\tr & \Zero_{s \times s}  & \underline{\Ahat} \\
\one\b\tr & \one\,\bhat\tr & \underline{\mathbf{A}} & \Zero_{s \times s} \\
\Xhline{2\arrayrulewidth}
\b\tr  & \bhat\tr & \underline{\b}\tr  & \underline{\bhat}\tr
\end{array}
\]
which is NOT a partitioned split GARK scheme. This is the form we actually compute when we compose steps. Note that the tableau has three nonzero entries in some row; even if we permute rows and columns we do not get the form below. This is the closest permutation:
\[
\begin{array}{cc|cc}
 \Zero_{s \times s}  & \Zero_{s \times s} & \Ahat  &\Zero_{s \times s} \\
\one\b\tr  & \Zero_{s \times s} & \one\,\bhat\tr  & \underline{\Ahat} \\
\hline
\mathbf{A}   &  \Zero_{s \times s} &  \Zero_{s \times s} & \Zero_{s \times s} \\
 \one\b\tr   &\underline{\mathbf{A}} & \one\,\bhat\tr&  \Zero_{s \times s}  \\
\Xhline{2\arrayrulewidth}
 \b\tr   & \underline{\b}\tr & \bhat\tr  & \underline{\bhat}\tr 
\end{array}
\]
\[
\begin{array}{cc|cc}
  \Zero_{s \times s} & \Zero_{s \times s} & \mathbf{A}   &  \Zero_{s \times s}\\
  \one\,\bhat\tr&  \Zero_{s \times s} & \one\b\tr   &\underline{\mathbf{A}}  \\
 \hline
  \Ahat  &\Zero_{s \times s} &  \Zero_{s \times s}  & \Zero_{s \times s} \\
 \one\,\bhat\tr  & \underline{\Ahat} & \one\b\tr  & \Zero_{s \times s} \\
\Xhline{2\arrayrulewidth}
 \b\tr   & \underline{\b}\tr & \bhat\tr  & \underline{\bhat}\tr \\
\end{array}
\]
Below we set the diagonal blocks to zero. Are we sure this leads to the same scheme as the composition form above? The positions in the second step start from the positions at the end of the first step, which depend on the momenta in the first step.

This last tableau is the form which we need to check for symmetry (automatic, I believe) and symplecticness  \eqref{eqn:P-matrix-gark} . 
\begin{equation*}
\b_{\textsc{gark}} \coloneqq \mbox{diag}\left(\b_{\textsc{gark}}\right), \quad 
\mathbf{P} = \mathbf{A}_{\textsc{gark}}\tran\,  \b_{\textsc{gark}} + \b_{\textsc{gark}}\, \mathbf{A}_{\textsc{gark}} 
-  \b_{\textsc{gark}} \,\underline{\b}_{\textsc{gark}}\tran.
\end{equation*}

The composition of a symplectic partitioned GARK scheme \eqref{eqn:Partitioned-Butcher-tableau} with its its time-reversed scheme \eqref{eqn:Partitioned-Butcher-tableau-reversed}
defines the scheme:
\begin{align*}
\raisebox{6pt}{$
\begin{array}{c|c}
\Zero_{2s \times 2s}  & \Ahat_{\textsc{comp}}
\\ \hline
\A_{\textsc{comp}}
  & \Zero_{2s \times 2s}
\\
\Xhline{2\arrayrulewidth}
\b_{\textsc{comp}} & \bhat_{\textsc{comp}}
\end{array}
$}
& ~\coloneqq ~
\raisebox{19pt}{$
\begin{array}{c|c}
\Zero_{2s \times 2s}  & 
\begin{array}{cc} 
\Ahat & \Zero_{s \times s} \\ 
\one\,\bhat\tr & \underline{\Ahat} % \one\, \bhat\tr - \mathcal{P} \Ahat {\mathcal{P}}
\end{array}
\\ \hline
\begin{array}{cc} 
\mathbf{A} & \Zero_{s \times s} \\ 
\one\,\b\tr & \underline{\A} % \one\, \b\tr - \mathcal{P} \mathbf{A } {\mathcal{P}}
\end{array}
  & \Zero_{2s \times 2s}
\\
\Xhline{2\arrayrulewidth}
\begin{array}{l} \b\tr  \hspace*{22pt} \underline{\b}\tr \end{array} & \begin{array}{l}   \bhat\tr  \hspace*{22pt} \underline{\bhat}\tr  \end{array}
\end{array}$},
\end{align*}
which, after a permutation of rows and columns, can be seen to be a partitioned scheme defined by
\begin{align}
\label{eqn:partitioned-composition}
& \left(\A_{\textsc{comp}}^{\{\ell,m\}}, \b_{\textsc{comp}}^{\{m\}}, \Ahat_{\textsc{comp}}^{\{\ell,m\}}, \bhat_{\textsc{comp}}^{\{m\}}\right) \\
\notag &\qquad =
\left(
    \begin{bmatrix}
        \A^{\{\ell,m\}} & \scriptstyle 0  \\ \one^{\{\ell\}} \b^{\{m\}}\tr &  \underline {\A}^{\{\ell,m\}}
    \end{bmatrix},
    \begin{bmatrix}
        \b^{\{m\}} \\ \underline{\b}^{\{m\}} 
    \end{bmatrix},
    \begin{bmatrix}
        \Ahat^{\{\ell,m\}} & \scriptstyle 0  \\ \one^{\{\ell\}} \bhat^{\{m\}}\tr & \underline{\Ahat}^{\{\ell,m\}}
    \end{bmatrix},
    \begin{bmatrix}
        \bhat^{\{m\}} \\ \underline{\bhat}^{\{m\}} 
    \end{bmatrix}
    \right).
\end{align}
Need to better explain how after permutations is equivalent to the form with blocks \eqref{eqn:partitioned-composition}. What are the ranges of $l,m$ in this context.

\begin{theorem}
\label{theorem.construction.part.symm.sympl.}
Start with a symplectic partitioned GARK scheme $(\A^{\{\ell,m\}},\b^{\{m\}},\Ahat^{\{\ell,m\}},\bhat^{\{m\}})$ and construct the partitioned composition method \eqref{eqn:partitioned-composition} defined by applying one step with the GARK scheme, followed by one step with its time-reversed version.
If $\bhat^{\{m\}} = \mathcal{P}\,\b^{\{m\}} = \underline{\b}^{\{m\}}$ and ${\b}^{\{m\}} = \bhat^{\{m\}}$ 
then the resulting composition scheme is both symmetric and symplectic.
\end{theorem}
\begin{proof}
Symmetry follows from the construction of the composition scheme. 
The composition scheme should be symmetric by construction. Note that the composition with the time-reversed scheme leads to some non-zero off-diagonal blocks. Why do we need additional conditions  $\bhat^{\{m\}} = \mathcal{P}\,\b^{\{m\}} = \underline{\b}^{\{m\}}$ and ${\b}^{\{m\}} = \bhat^{\{m\}}$ to prove it? The fact that we artificially set the diagonal blocks to zero in the composition scheme may be the cause?

We check the symmetry relations \eqref{eqn:GARK-symmetry-matrix.a}:
\begin{align*}
 \underline{\b}^{\{m\}}_{\textsc{comp}} = \mathcal{P}_{2s} \,\b_{\textsc{comp}}^{\{m\}} & =
 \begin{bmatrix}
      \scriptstyle 0 & \mathcal{P} \\ \mathcal{P} & \scriptstyle 0 
  \end{bmatrix}
  \begin{bmatrix}
  \b^{\{m\}} \\ \underline{\b}^{\{m\}}
  \end{bmatrix} = 
  \begin{bmatrix}
  \mathcal{P}\,\underline{\b}^{\{m\}}\\ \mathcal{P}\,\b^{\{m\}} 
  \end{bmatrix} =
  \b_{\textsc{comp}},
\end{align*}
and similarly we check that $\bhat_{\textsc{comp}} = \mathcal{P}_{2s}\,  \bhat_{\textsc{comp}}$.

We next check the symmetry relations \eqref{eqn:GARK-symmetry-matrix.b}:
\[
\underline{\A}^{\{\ell,m\}}  \coloneqq  \one^{\{\ell\}}\,\b^{\{m\}}\tr - \mathcal{P}^{\{\ell\}}\,\A^{\{\ell,m\}}\,\mathcal{P}^{\{m\}}
\]
\begin{align*}
\underline{\A}_{\textsc{comp}}^{\{\ell,m\}} & = 
    \begin{bmatrix}
        \one^{\{\ell\}} \b^{\{m\}}\tr & \one^{\{\ell\}} \underline{\b}^{\{m\}}\tr \\
        \one^{\{\ell\}} \b^{\{m\}}\tr & \one^{\{\ell\}} \underline{\b}^{\{m\}}\tr
    \end{bmatrix}
    -
    \begin{bmatrix}
        \scriptstyle 0 & \mathcal P \\ \mathcal P & \scriptstyle 0  
    \end{bmatrix}
    \begin{bmatrix}
        \A^{\{\ell,m\}} & \scriptstyle 0  \\ \one^{\{\ell\}} \b^{\{m\}}\tr &  \underline {\A}^{\{\ell,m\}}
    \end{bmatrix}
    \begin{bmatrix}
        \scriptstyle 0 & \mathcal P \\ \mathcal P & \scriptstyle 0  
    \end{bmatrix} \\
& = \begin{bmatrix}
    \one^{\{\ell\}} \b^{\{m\}}\tr - \mathcal{P}\, \underline{\A}^{\{\ell,m\}}\,  \mathcal{P} 
    & \scriptstyle 0  \\
    \one^{\{\ell\}} \b^{\{m\}}\tr & \one^{\{\ell\}} \underline{\b}^{\{m\}}\tr - \mathcal{P}\, \A^{\{\ell,m\}}\,  \mathcal{P} \\
\end{bmatrix} \\
& = \begin{bmatrix}
    \one^{\{\ell\}} (\b^{\{m\}}-\underline{\b}^{\{m\}})\tr + \A^{\{\ell,m\}}
    & \scriptstyle 0  \\
    \one^{\{\ell\}} \b^{\{m\}}\tr & \one^{\{\ell\}} (\underline{\b}^{\{m\}}-\b^{\{m\}})\tr + \underline{\A}^{\{\ell,m\}}
\end{bmatrix}
\end{align*}

If $\underline{\b}^{\{m\}} = \b^{\{m\}}$ then

The condition 
$\Ahat_{\textsc{comp}}=\underline{\Ahat}_{\textsc{comp}} $ reads
\begin{align*}
    \begin{bmatrix}
\Ahat & \Zero_{s \times s} \\ 
\one\,\bhat\tr & \one\, \bhat\tr - \mathcal{P} \Ahat {\mathcal{P}}
\end{bmatrix} & =
\begin{bmatrix}
    \one \bhat\tr & \one \underline{\bhat}\tr \\
    \one \bhat\tr & \one \underline{\bhat}\tr
\end{bmatrix} -
\mathcal{P}_{2s} 
\begin{bmatrix}
  \Ahat & \Zero_{s \times s} \\ 
\one\,\bhat\tr & \one\, \bhat\tr - \mathcal{P} \Ahat {\mathcal{P}}  
\end{bmatrix}
\mathcal{P}_{2s} \\
& = \begin{bmatrix}
    \one \bhat\tr - \left( \mathcal{P} \one \bhat\tr \mathcal{P} - \mathcal{P}^2 \Ahat \mathcal{P}^2 \right) & \one \underline{\bhat}\tr - \mathcal{P} \one \bhat\tr \mathcal{P} \\
    \one \bhat\tr & \one \underline{\bhat}\tr - \mathcal{P} \Ahat \mathcal{P}
\end{bmatrix}
\end{align*}
which is equivalent to  
$\bhat = \underline{\bhat}$. 
%$\mathcal P \b^{\{m\}} = \b^{\{m\}}=\underline{\b}^{\{m\}}$. 
The condition 
$\A_{\textsc{comp}}=\underline{\A}_{\textsc{comp}} $
holds, if $\b =  \underline{\b}$.

Symplecticness requires
\begin{align*}
 \Ahat_{\textsc{comp}} & = \ones \b_{\textsc{comp}}\tr - \Bhat_{\textsc{comp}}^{-1}
 \A_{\textsc{comp}}\tr \B_{\textsc{comp}},
\end{align*}
%\begin{align*}
%    \begin{bmatrix}
%        \Ahat^{\{\ell,m\}} & \scriptstyle 0  \\ \one \bhat^{\{m\}} & \underline{\Ahat}^{\{\ell,m\}}
%    \end{bmatrix}  = &
%    \begin{bmatrix}
%        \one \b^{\{m\}} & \one \b^{\{m\}} \\
%        \one \underline{\b}^{\{m\}} & \one \underline{\b}^{\{m\}}
%    \end{bmatrix} - \\
%    &
%    \begin{bmatrix}
%        \B^{\{\ell\}-1} & \scriptstyle 0  \\ \scriptstyle 0 & \Bhat^{\{\ell\}-1}
%    \end{bmatrix}%
%    \begin{bmatrix}
%        \A^{\{m,\ell\}}\tr & \b^{\{\ell\}} \one\tr \\ \scriptstyle 0 &  \underline {\A}^{\{m,\ell\}}\tr
%    \end{bmatrix}
%     \begin{bmatrix}
%        \B^{\{m\}} & \scriptstyle 0  \\ \scriptstyle 0 & \Bhat^{\{m\}}
%    \end{bmatrix},
%\end{align*}
which holds for a symplectic partitioned scheme $(\A^{\{\ell,m\}},\b^{\{m\}},\Ahat^{\{\ell,m\}},\bhat^{\{m\}})$, if in addition $\b=\bhat$ holds. 

\end{proof}
\begin{remark}
\label{rem.constr.expl}
One drawback of this construction is the low number of degrees of freedom $\lceil s^{\{m\}}/2 \rceil$ in selecting the weights $\b^{\{m\}}$.
\end{remark}
\color{black}

If we aim at explicit schemes with an even number of stages $s^{\{\ell\}}$, a necessary property is that the matrices \guenther{$\A_{\textsc{comp}}^{\{l,l\}}$} %$\A^{\{l,l\}}$ are at least (after some rearrangement of stages)  lower triangular matrices. Symmetry then requires \guenther{we write better suggests} the general form \sandu{Why is this the most general form? Clearly this leads to composition schemes, but not clear we cannot do otherwise.}
\begin{equation*}
%\A^{\{l,l\}} 
\A_{\textsc{comp}}^{\{l,l\}}=
    \begin{bmatrix}
        \left( \one^{\{\ell\}} \one^{\{\ell\}}\tr - X_{l,l}\right) \widetilde  \B^{\{\ell\}} & \scriptstyle 0  \\
        \one^{\{\ell\}} \widetilde \b^{\{\ell\}}\tr & P_l \left( X_{l,l} \widetilde \B^{\{\ell\}} \right) P_l
    \end{bmatrix}
\end{equation*}
with weights $\b_{\textsc{comp}}^{\{\ell\}}\tr=(\widetilde \b^{\{\ell\}}\tr,(\mathcal{P} \widetilde \b^{\{\ell\}})\tr)\tr$, $\mathcal{P}_l$ being the permutation matrix of dimension $s^{\{\ell\}}/2$, $\widetilde B^{\{\ell\}}\coloneqq\diag(\widetilde b^{\{\ell\}})$ and arbitrary matrices $X_{l,l} \in \Re^{{s^{\{\ell\}}/2 \times{s^{\{\ell\}}/2}}}$.
The conjugate scheme is given by
\begin{equation*}
  % \Ahat^{\{l,l\}} 
   \Ahat_{\textsc{comp}}^{\{l,l\}}=
    \begin{bmatrix}
        X_{l,l}\tr \widetilde \B^{\{\ell\}}  & \scriptstyle 0  \\
        \one^{\{\ell\}} \widetilde \b^{\{\ell\}}\tr & 
        P_l  (\one^{\{\ell\}} \one^{\{\ell\}}\tr - X_{l,l}\tr) \widetilde \B^{\{\ell\}}  P_l
    \end{bmatrix}, \quad
    \bhat^{\{\ell\}} = \b^{\{\ell\}}.
\end{equation*}
To get an overall symmetric and symplectic scheme, the coupling matrices might be set by
\begin{align*}
%\A^{\{\ell,m\}} 
\A_{\textsc{comp}}^{\{\ell,m\}} &=
    \begin{bmatrix}
        \left( \one^{\{\ell\}} \one_m\tr - X_{l,m}\right) \widetilde  \B^{\{m\}} & \scriptstyle 0  \\
        \one^{\{\ell\}} \widetilde \b^{\{m\}}\tr & P_l \left( X_{l,m} \widetilde \B^{\{m\}} \right) P_m
    \end{bmatrix}, \\
%    \Ahat^{\{\ell,m\}} 
 \Ahat_{\textsc{comp}}^{\{\ell,m\}}   
    &=
    \begin{bmatrix}
        X_{m,l}\tr \widetilde \B^{\{m\}}  & \scriptstyle 0  \\
        \one^{\{\ell\}} \widetilde \b^{\{m\}}\tr & 
        P_l  (\one^{\{\ell\}} \one_m\tr - X_{m,l}\tr) \widetilde \B^{\{m\}}  P_m
    \end{bmatrix}.
\end{align*}

\begin{remark}
This construction is equivalent to the construction of Theorem~\ref{theorem.construction.symm.sympl.} with choosing symmetric matrices
$\A^{\{\ell,m\}}\coloneqq\left( \one^{\{\ell\}} \one_m\tr - X_{l,m}\right) \widetilde  \B^{\{m\}}$ and $\Ahat^{\{\ell,m\}}\coloneqq X_{m,l}\tr \widetilde \B^{\{m\}}$
based on time-reversed schemes, if $\mathcal P \widetilde \b^{\{m\}}= \widetilde \b^{\{m\}}$ holds. However, this construction allows for arbitrary weights $\b^{\{m\}}$ to keep the degree of freedom $s^{\{m\}}$ in these weigths.
\end{remark}
\end{leaveout}

%%%%%%%%%%%%%%%%%%%%%%%%%
\subsection{Partitioned GARK schemes for potential splitting}
%%%%%%%%%%%%%%%%%%%%%%%%%

We now consider the often encountered case where only the potential is split: 
    \begin{subequations}
    \label{eqn:split-V}
    \begin{equation}
    \label{eqn:split-V-Hamiltonian-sum}
    H(\p,\q) = T(\p) + V(\q) \quad \mbox{with} \quad  V(\q) = \sum_{m=2}^{\nparts} V^{\{m\}}(\q),
    \end{equation}
    and where we have the following partitioned Hamiltonian  \eqref{eqn:partitioned-Hamiltonians}
    \begin{align}
    \label{eqn:split-V-Hamiltonian}
        H(\p,\q)\!=\! &\sum_{m=1}^{\nparts}\! H^{\{m\}}(\p,\q) \;\; \mbox{with} \;\; 
        H^{\{m\}}(\p,\q)  = \begin{cases} 
        T(\p), & m = 1, \\
        V^{\{m\}}(\q), & m=2,\ldots,\nparts.
        \end{cases}
    \end{align}
    \end{subequations}
We note that the potential split system \eqref{eqn:split-V} is a special case of a partitioned system \eqref{eqn:split-TV}  with $V^{\{1\}}(\q)= 0$ and $T^{\{m\}}(\p) = 0$ for $m = 2,\dots,\nparts$.

The partitioned GARK scheme \eqref{eq:sympl.part.gark} applied to the potential splitting \eqref{eqn:split-V} reads:
\begin{subequations}
\label{eqn:GARK-symplectic.partitioned.separable}
\begin{align}
\label{eqn:GARK-symplectic.stage-solutions.partitioned.separable}
P_i^{\{1\}} & =  \p_0 + h \sum_{m=2}^{\nparts} \sum_{j=1}^{\hat s^{\{m\}}} \ahat_{i,j}^{\{1,m\}} k_j^{\{m\}},
%, \qquad i=1,\ldots,s^{\{q\}},\\
& Q_i^{\{q\}} & =  \q_0 + h \sum_{j=1}^{s^{\{1\}}}  a_{i,j}^{\{q,1\}} \ell_j^{\{1\}}, \\
%\qquad i=1,\ldots,\hat s^{\{q\}},\\
\label{eqn:GARK-symplectic.solutions.partitioned.separable}
\p_1 & =  \p_0 + h \sum_{q=2}^{\nparts} \sum_{i=1}^{\hat s^{\{q\}}} \widehat{b}_{i}^{\{q\}} k_i^{\{q\}}, &
\q_1 & =  \q_0 + h \sum_{i=1}^{s^{\{1\}}} {b}_{i}^{\{1\}} \ell_i^{\{1\}}, \\
\label{eqn:GARK-symplectic.stage-functions.partitioned.separable}
k_i^{\{m\}} & =  -V_{\q}^{\{m\}}(Q_i^{\{m\}}), &
\ell_i^{\{m\}} & =  T_{\p}^{\{m\}}(P_i^{\{m\}}).
\end{align}
\end{subequations}
One notes that the stage vectors $P_i^{\{q\}}$ for $q=2,\ldots,\nparts$ and $Q_i^{\{1\}}$ are not needed for this type of splitting. 
The corresponding generalized Butcher tableau \eqref{eqn:general-Butcher-tableau} is:
%
%
%%\label{eqn:Partitioned-Butcher-tableau}
%
%The GARK scheme \eqref{eqn:GARK-symplectic} applied to the potential splitting \eqref{eqn:split-V} reads:
%%
%\begin{subequations}
%\label{eqn:GARK-symplectic.partitioned.separable}
%\begin{align}
%\label{eqn:GARK-symplectic.stage-solutions.partitioned.separable}
%P_i^{\{1\}} & = \p_0 + h \sum_{m=2}^{\nparts} \sum_{j=1}^{s^{\{m\}}} a_{i,j}^{\{1,m\}} k_j^{\{m\}}, & 
%Q_i^{\{q\}} & = \q_0 + h \sum_{j=1}^{s^{\{m\}}} a_{i,j}^{\{q,1\}} \ell_j^{\{1\}},\\
%\label{eqn:GARK-symplectic.stage-functions.partitioned.separable}
%k_i^{\{m\}} & = - V^{\{m\}}_{\q}\left(Q_i^{\{m\}}\right),~~ m \ge 2, &
%\ell_i^{\{1\}} & = T_{\p}(P_i^{\{1\}}), \\ 
%\label{eqn:GARK-symplectic.solutions.partitioned.separable}
%\p_1 & =  \p_0 + h \sum_{q=2}^{\nparts} \sum_{i=1}^{s^{\{q\}}} b_{i}^{\{q\}} k_i^{\{q\}}, &
%\q_1 & =  \q_0 + h \sum_{i=1}^{s^{\{1\}}} b_{i}^{\{1\}} \ell_i^{\{1\}}.
%\end{align}
%\end{subequations}
%%
%\begin{subequations}
%\label{eqn:GARK-symplectic.partitioned.separable}
%\begin{eqnarray}
%P_i^{\{1\}} & = & \p_0 - h  \sum_{m=2}^{\nparts} \sum_{j=1}^{s^{\{m\}}} a_{i,j}^{\{1,m\}} \cdot 
%V^{\{m\}}_{\q}(Q_i^{\{m\}}), \\
%Q_i^{\{m\}} & = & \q_0 + h  \sum_{j=1}^{s^{\{1\}}} a_{i,j}^{\{m,1\}} \cdot T_{\p}(P_i^{\{1\}}), \quad q=2,\ldots,\nparts, \\
%\p_1 & = & \p_0 - h  \sum_{m=2}^{\nparts} \sum_{i=1}^{s^{\{q\}}} b_{i}^{\{m\}} \cdot V^{\{m\}}_{\q}(Q_i^{\{m\}}), \\
%\q_1 & = & \q_0 + h  \sum_{i=1}^{s^{\{1\}}} b_{i}^{\{1\}} \cdot T_{\p}(P_i^{\{1\}}).
%\end{eqnarray}
%\end{subequations}
%
\begin{equation}
\renewcommand{\arraystretch}{1.5}
\label{eqn:Partitioned-Butcher-tableau-PS}
\begin{array}{c|ccc}
\Zero  & \Ahat^{\{1,2\}} & \cdots & \Ahat^{\{1,\nparts\}} \\
\hline
 \mathbf{A}^{\{2,1\}} &\Zero&&\Zero\\
 \vdots &&\ddots&\\
 \mathbf{A}^{\{\nparts, 1\}} &\Zero&&\Zero\\ 
 \Xhline{2\arrayrulewidth}
\b ^{\{1\}}\tr & \bhat^{\{2\}}\tr & \cdots & \bhat^{\{\nparts\}}\tr
\end{array}\,\raisebox{-26pt}{$.$}
\end{equation}
The generalized momenta are obtained by integrating each potential $V^{\{m\}}$ with a Runge-Kutta scheme $(\bhat^{\{m\}},\Ahat^{\{1,m\}})$ for $m = 2,\ldots,\nparts$. The generalized positions are obtained by integrating the kinetic energy with  a Runge-Kutta scheme $(\b ^{\{1\}},\mathbf{A}^{\{m,1\}})$ for $m = 2,\ldots,\nparts$. All other coupling coefficients are zero.

\begin{leaveout}
\sandu{This probably not needed here, proof to be moved earlier:
\begin{corollary}[Symplecticity]
\label{cor.part.ham.2}
The necessary and sufficient conditions for the symplecticity of the GARK scheme~\eqref{eqn:GARK-symplectic.partitioned.separable}
are $\mathbf{P}^{\{1,m\}}=0$ for $m \in \{2,\ldots,\nparts\}$. %The other conditions \eqref{eqn:GARK-symplectic} are not needed.
\end{corollary}
\begin{proof}
From \eqref{eqn:GARK-symplectic.stage-functions} applied to \eqref{eqn:split-V}  we see that $k_i^{\{1\}} = \Zero$ and $\ell_i^{\{m\}} = \Zero$ for $m=2,\dots,\nparts$. Consequently, using the partition
\[
\mathfrak{g}_i^{\{m\}} = \begin{bmatrix}
\mathfrak{g}_{k,i}^{\{m\}} \\
\mathfrak{g}_{l,i}^{\{m\}}
\end{bmatrix}
\]
in the proof of Theorem~\ref{theorem.symplectic} gives:
\[
\begin{bmatrix}
k_i^{\{1\}} \\ \ell_i^{\{1\}}
\end{bmatrix} 
=
\begin{bmatrix}
0 \\ \textnormal{NB}(\mathfrak{g}_{l,i}^{\{m\}},[\q_0, \p_0])
\end{bmatrix};
\qquad 
%\]
%and 
%\[
\begin{bmatrix}
k_i^{\{m\}} \\ \ell_i^{\{m\}}
\end{bmatrix} 
=
\begin{bmatrix}
\textnormal{NB}(\mathfrak{g}_{k,i}^{\{m\}},[\q_0, \p_0]) \\
0
\end{bmatrix}, ~~ m \ge 2. 
\]
Hence we have for
\begin{align*}
    m=1: \qquad & U \coloneqq  
    \begin{bmatrix}
    0 \\ \mathfrak{g}_\ell^{\{1\}} (u_1)
    \end{bmatrix} 
    \times \cdots \times
    \begin{bmatrix}
    0 \\ \mathfrak{g}_\ell^{\{1\}} (u_r)
    \end{bmatrix},  \quad
    V \coloneqq  
    \begin{bmatrix}
    0 \\ \mathfrak{g}_\ell^{\{1\}} (v_1)
    \end{bmatrix} 
    \times \cdots \times
    \begin{bmatrix}
    0 \\ \mathfrak{g}_\ell^{\{1\}} (v_r)
    \end{bmatrix}, \\
    m >1: \qquad & 
   U \coloneqq  \begin{bmatrix}
    \mathfrak{g}_k^{\{m\}} (u_1) \\ 0
    \end{bmatrix} \times \cdots \times
    \begin{bmatrix}
    \mathfrak{g}_k^{\{m\}} (u_r) \\ 0
    \end{bmatrix},  \quad
    V \coloneqq  \begin{bmatrix}
    \mathfrak{g}_k^{\{m\}} (v_1) \\ 0
    \end{bmatrix} \times \cdots \times
    \begin{bmatrix}
    \mathfrak{g}_k^{\{m\}} (v_r) \\ 0
    \end{bmatrix},   
\end{align*}
for $m>1$.
Thus $U_i V_j \neq 0$ for $m=1, \ell >1$ or  $m>1, \ell =1$ and
\eqref{eq.derivation.sympl} implies
$\mathbf{P}^{\{1,m\}}=0$ for $m=2,\dots,\nparts$ due to the symmetry of $\mathbf{P}$.
\end{proof}
}
\end{leaveout}
\begin{leaveout}
\begin{proof}
From \eqref{eqn:GARK-symplectic.stage-functions} applied to \eqref{eqn:split-V}  we see that $k_i^{\{1\}} = \Zero$ and $\ell_i^{\{m\}} = \Zero$ for $m=2,\dots,\nparts$. Consequently (see \Cref{app:proof-of-symplecticness}):
\begin{eqnarray*}
dk_i^{\{1\},J} \wedge d\ell_j^{\{m\},J}=0 \quad \mbox{for } m=1,\ldots,\nparts, \\
dk_i^{\{\ell\},J} \wedge d\ell_j^{\{m\},J}=0 \quad \mbox{for } \ell,m =2,\ldots,\nparts.
\end{eqnarray*}
The conditions $\mathbf{P}^{\{1,1\}}=0$  and $\mathbf{P}^{\{\ell,m\}}=0$ for $\ell,m \in \{2,\ldots,\nparts\}$
are not necessary for the symplecticity of the GARK scheme. Conditions $\mathbf{P}^{\{1,m\}}=0$ 
imply $\mathbf{P}^{\{m,1\}}=0$ for $m=2,\dots,\nparts$ due to the symmetry of $\mathbf{P}$.
\end{proof}
\end{leaveout}

As the partitioned GARK scheme~\eqref{eqn:GARK-symplectic.partitioned.separable} is a special case of~\eqref{eq:sympl.part.gark}, Theorem~\ref{cor.part.ham} yields 
%The necessary and sufficient conditions for the symplecticity of the GARK scheme~\eqref{eqn:GARK-symplectic.partitioned.separable}
%are $\mathbf{P}^{\{1,m\}}=0$ for $m \in \{2,\ldots,\nparts\}$.
%
the symplecticity conditions \eqref{eqn:symplecticness} %read:
\begin{equation*}
    \mathbf{P}^{\{1,m\}} = \mathbf{A}^{\{m,1\}}\tr\,  \mathbf{\hat B}^{\{m\}} + \B^{\{1\}}\, \Ahat^{\{1,m\}} -  \b ^{\{1\}} \,\bhat^{\{m\}}\tr  = 0,  \qquad m=2,\ldots,\nparts,
\end{equation*}
and can be solved for $\Ahat^{\{1,m\}}$ when all entries of $\b^{\{1\}}$ are nonzero \eqref{part.GARK.sympl.cond}:
\begin{align}
\label{def.conjugate}
  \Ahat^{\{1,m\}} & %= \mathbf{\hat A}^{\{1,m\}} 
  \coloneqq \one^{\{1\}}  \,\bhat^{\{m\}}\tr - \B^{\{1\}-1}\, \mathbf{A}^{\{m,1\}}\tr\, \Bhat^{\{m\}},
  \qquad m = 2,\ldots,\nparts.
\end{align}
In the symplectic case, each $\Ahat^{\{1,m\}}$ is uniquely defined in terms of
the Runge Kutta scheme  $(\b ^{\{1\}},\mathbf{A}^{\{m,1\}})$ and the weight vector $\bhat^{\{m\}}$.

%\begin{remark}
%In \eqref{eqn:Partitioned-Butcher-tableau-PS} we formally set the coefficients $\A^{\{1,1\}}$ and $\A^{\{m,n\}}$ for $m,n>1$ to zero. This follows the fact that all elementary differentials associated to trees involving these matrices are zero, and thus the corresponding order conditions are not required:
%%
%\begin{align*}
%    F(\tau_1) & = \begin{bmatrix}
%    T_{\p}(\p) \\ 0
%    \end{bmatrix}, \quad
%    F(\tau_{\{m\}}) = \begin{bmatrix}
%    0 \\ - V_{\q}^{\{m\}}(\q)
%    \end{bmatrix}  \Rightarrow \\
%    F([\tau_1]_1) & = \frac{\partial f^{\{1\}}(\p)}{\partial (\q,\p)} F(\tau_1) = \begin{bmatrix}
%    \scriptstyle 0 & T_{\p,\p}(\p) \\ \scriptstyle 0 & \scriptstyle 0 
%    \end{bmatrix}
%    \cdot \begin{bmatrix}
%    T_{\p}(\p) \\ 0
%    \end{bmatrix} = 0, \\
%    F([\tau_{\{m\}}]_n) & = \frac{\partial f^{\{n\}}(\q)}{\partial (\q,\p)} F(\tau_{\{m\}}) = \begin{bmatrix}
%    \scriptstyle 0 & \scriptstyle 0  \\ -V_{\q,\q}^{\{n\}}(\q)  & \scriptstyle 0 
%    \end{bmatrix}
%    \cdot  \begin{bmatrix}
%    0 \\ - V_{\q}^{\{m\}}(\q)
%    \end{bmatrix}   =0.
%\end{align*}
%\end{remark}

The order conditions up to order four are obtained directly from \eqref{eq.ordercond.3.2}
(refer to Theorem~\ref{th.sympl.order.cond}):
\begin{align*}
\b ^{\{1\}}\tr \cdot \one^{{\{1\}}} = 1,  &  &  (\textnormal{order}~ 1) \\
\bhat^{\{m\}}\tr \cdot \one^{{\{m\}}} = 1,  & \quad \forall\; m>1, &  (\textnormal{order}~ 1) \\
 \bhat^{\{m\}}\tr \cdot \mathbf{c}^{\{m,1\}}  =  \sfrac{1}{2}, &\quad  \forall\; m>1,  &  (\textnormal{order}~ 2) \\
\b ^{\{1\}}\tr \cdot \left(  \chat^{\{1,\ell\}} \times \chat^{\{1,s\}} \right)
 =  \sfrac{1}{3}, & \quad \forall\; \ell, s>1,\; \ell \leq s, & (\textnormal{order}~ 3) \\
\bhat^{\{m\}}\tr \cdot \left(  \mathbf{c}^{\{m,1\}} \times \mathbf{c}^{\{m,1\}} \right)
 =  \sfrac{1}{3}, & \quad \forall\; m>1, & (\textnormal{order}~ 3) \\
\b ^{\{1\}}\tr \cdot \left( \chat^{\{1,\ell\}} \times \chat^{\{1,s\}} \times \chat^{\{1,t\}} \right)
 =  \sfrac{1}{4}, & \quad \forall\; \ell,s,t>1,\; \ell \leq s \leq t, & (\textnormal{order}~ 4) \\
\bhat^{\{m\}}\tr \cdot \left( \mathbf{c}^{\{m,1\}} \times \mathbf{c}^{\{m,1\}} \times \mathbf{c}^{\{m,1\}} \right)
 =  \sfrac{1}{4}, & \quad \forall\; m>1, & (\textnormal{order}~ 4) \\
\b ^{\{1\}}\tr  \cdot  \left(\chat^{\{1,\ell\}} \times \Ahat^{\{1,s\}}  \cdot \mathbf{c}^{\{s,1\}} \right) 
 =  \sfrac{1}{8}, & \quad \forall\; \ell,s>1, & (\textnormal{order}~ 4) \\
 \bhat^{\{m\}}\tr  \cdot  \left(\mathbf{c}^{\{m,1\}} \times \mathbf{A}^{\{m,1\}}  \cdot \chat^{\{1,t\}} \right) 
 =  \sfrac{1}{8}, & \quad \forall\; m,t>1. & (\textnormal{order}~ 4) 
\end{align*}

%
%The GARK scheme reads \sandu{Why do we repeat \eqref{eqn:GARK-symplectic.partitioned.separable}?}
%%
%\begin{subequations}
%\label{eqn:GARK-symplectic.partitioned.separable2}
%\begin{eqnarray}
%P_i^{\{1\}} & = & \p_0 - h  \sum_{m=2}^{\nparts} \sum_{j=1}^{s^{\{m\}}} {\hat a}_{i,j}^{\{1,m\}} 
%\cdot 
%V^{\{m\}}_{\q}(Q_i^{\{m\}}), \\
%Q_i^{\{q\}} & = & \q_0 + h  \sum_{j=1}^{s^{\{m\}}} a_{i,j}^{\{q,1\}} \cdot T_{\p}(P_i^{\{1\}}), \quad q=2,\ldots,\nparts\\
%\p_1 & = & \p_0 - h  \sum_{m=2}^{\nparts} \sum_{i=1}^{s^{\{q\}}} b_{i}^{\{m\}} \cdot V^{\{m\}}_{\q}(Q_i^{\{m\}}), \\
%\q_1 & = & \q_0 + h  \sum_{i=1}^{s^{\{q\}}}  b_{i}^{\{1\}} \cdot T_{\p}(P_i^{\{1\}}).
%\end{eqnarray}
%\end{subequations}
\begin{remark}
For $\nparts=2$ and $s^{\{1\}}=s^{\{2\}}$ the scheme \eqref{eqn:Partitioned-Butcher-tableau-PS} is equivalent to a traditional Partitioned RK scheme with Butcher tableau: 
\begin{equation}
\renewcommand{\arraystretch}{1.5}
\label{eqn:Partitioned-Butcher-tableau2}
\begin{array}{c|c}
\mathbf{A} & \Ahat \\
 \Xhline{2\arrayrulewidth}
\b \tr & \bhat\tr
\end{array}
\quad \textnormal{with} \quad
\left\{
\begin{array}{ll}
   \mathbf{A} = \mathbf{A}^{\{2,1\}}, &
   \Ahat = \Ahat^{\{1,2\}}=\ones  \,\bhat\tr - \B^{-1}\, \mathbf{A}^{\{2,1\}}\tr\, \Bhat, \\
   \b  = \b ^{\{1\}}, &
   \bhat = \b ^{\{2\}}. 
\end{array} 
\right.
\end{equation}
%
%where 
%%
%\begin{align*}
%   \mathbf{A}&= \mathbf{A}^{\{2,1\}}, \quad
%   \Ahat = \Ahat^{\{1,2\}}=\ones  \,\bhat\tr - \B^{-1}\, \mathbf{A}^{\{2,1\}}\tr\, \Bhat, \quad
%   \b  = \b ^{\{1\}}, \quad
%   \bhat = \b ^{\{2\}}. 
%\end{align*} 
%The choice of $N-1$ basic RK schemes $(\b ^{\{2\}},\mathbf{A}^{\{2,1\}}), \ldots, (\b ^{\{\nparts\}},\mathbf{A}^{\{\nparts,1\}})$, together with the choice of $\b ^{\{1\}}$ defines the partitioned symplectic GARK scheme~\eqref{eqn:GARK-symplectic.partitioned.separable}. It can be described by its Butcher tableau
%\begin{equation*}%begin{equation}
%\renewcommand{\arraystretch}{1.5}
%%\label{eqn:Partitioned-Butcher-tableau2}
%\begin{array}{cccc}
 % & \mathbf{\hat A}^{\{2,1\}} & \cdots & \mathbf{\hat A}^{\{2,1\}} \\
 %\mathbf{A}^{\{1,2\}} &&&\\
 %\vdots &&&\\
 %\mathbf{A}^{\{1,\nparts\}} &&&\\ \hline
%\b ^{\{1\}}\tr & \bhat^{\{2\}}\tr & \cdots & \bhat^{\{\nparts\}}\tr
%\end{array}
%\end{equation*}%end{equation}
%with \sandu{Why repeat?}
%\begin{equation*}
 %\mathbf{\hat A}^{\{1,m\}} & \coloneqq  \ones  \,\bhat^{\{m\}}\tr - \B^{\{1\}-1} \mathbf{A}^{\{m,1\}}\tr \mathbf{\hat B}^{\{m\}},
  %\quad m=2,\ldots,\nparts\\
  %\bhat^{\{m\}} &\coloneqq\b^{\{m\}}, \quad  m=2,\ldots,\nparts.
%\end{equation*}
\end{remark}

\subsection{Explicit partitioned GARK schemes}
One idea for constructing explicit symmetric and symplectic GARK schemes is  to define a two-step scheme: take a first step with an explicit symplectic partitioned GARK scheme, and then a second step with its time-reversed scheme (see Section~\ref{sec3.4.2}). However, the resulting symplectic and symmetric scheme might not fall into the class of partitioned GARK schemes, as discussed in the proof of Theorem \ref{theorem.construction.symm.sympl2.}.

In the following we will consider ideas how to construct explicit schemes that define partitioned GARK schemes.

\subsubsection*{Explicit symplectic partitioned GARK schemes}
%\sandu{A composition scheme is not necessarily a partitioned GARK. Taking an explicit symplectic GARK followed by its time reversed scheme is probably a case worth pursuing, but it does not immediately fit in the discussion below.}

%
Partitioned GARK schemes \eqref{eqn:Partitioned-Butcher-tableau} are explicit iff they fulfill the condition~\cite{Sandu_2015_GARK} %\sandu{Why do we take the absolute values below?}
\begin{align}
\label{cond.explicit}
    \mathbf{S}^{\{\ell,m\}}:=  \Ahat^{\{\ell,m\}}  \times  \A^{\{m,\ell\}}\tr = 0 \quad \forall\, \ell,m =1,\ldots,\nparts,
\end{align}
where $|\cdots|$ takes  element-wise  absolute  values,  and $\times$ is  the  element-wise  product. 

For symplectic partitioned GARK schemes \eqref{eqn:Partitioned-Butcher-tableau} schemes with non-vanishing weights $\b^{\{\ell\}}$
and $\bhat^{\{m\}}$, equation~\eqref{part.GARK.sympl.cond} 
\[
\ahat^{\{\ell,m\}}_{i,j} =  \widehat{b}^{\{m\}}_{j} \, \left( 1
- a^{\{m,\ell\}}_{j,i}/b^{\{\ell\}}_i \right), \\
\]
together with condition \eqref{cond.explicit} lead to:
\begin{equation}
\label{eq:explicit-conditions}
\ahat^{\{\ell,m\}}_{i,j} = \widehat{b}^{\{m\}}_j\,(1-x_{j,i}^{\{m,\ell\}}) \mbox{ and } a^{\{m,\ell\}}_{j,i} = b^{\{\ell\}}_i\, x_{j,i}^{\{m,\ell\}},
\quad
x_{j,i}^{\{m,\ell\}} \in \{0,1\},
\end{equation}
%
%\begin{cases}
%\textnormal{either} &  \ahat^{\{\ell,m\}}_{i,j} = \widehat{b}^{\{m\}}_j \mbox{ and } a^{\{m,\ell\}}_{j,i}=0, \\[3pt]
%\textnormal{or} & 
%\ahat^{\{\ell,m\}}_{i,j}=0  \mbox{ and } a^{\{m,\ell\}}_{j,i} = b^{\{\ell\}}_i,
%\end{cases}
%\quad
%\forall\; i = 1,\dots,s ^{\{\ell\}}; \, j = 1,\dots,s ^{\{m\}}.
%
for all $i = 1,\dots,s^{\{\ell\}}; \, j = 1,\dots,s^{\{m\}}$,
or in compact notation with $\mathbf{x} \coloneqq (x_{i,j})_{i,j}$ and $\mathbf{\widehat{x}} \coloneqq (1-x_{i,j})_{i,j}$:
\begin{subequations}
\begin{align*}
    \A^{\{m,\ell\}} & = \mathbf{x}^{\{m,\ell\}} \cdot \B^{\{\ell\}}, \\
    \Ahat^{\{\ell,m\}} & = \mathbf{\widehat{x}}^{\{\ell,m\}}  \cdot \Bhat^{\{m\}}.
\end{align*}
\end{subequations}

The scheme \eqref{eq:sympl.part.gark} then reads:
\begin{align*}
P_i^{\{m\}} & =  \p_0 + h \sum_{\ell=1}^{\nparts} \sum_{j=1}^{\hat s^{\{\ell\}}} \ahat_{i,j}^{\{m,\ell\}} k_j^{\{\ell\}}
=  \p_0 + h \sum_{\ell=1}^{\nparts} \sum_{j=1}^{\hat s^{\{\ell\}}} \widehat{x}_{j,i}^{\{\ell,m\}}\, \widehat{b}^{\{\ell\}}_j\, k_j^{\{\ell\}} , \\
%, \qquad i=1,\ldots,s^{\{q\}},\\
 Q_i^{\{m\}} & =  \q_0 + h \sum_{\ell=1}^{\nparts} \sum_{j=1}^{s^{\{\ell\}}}  a_{i,j}^{\{m,\ell\}} \ell_j^{\{\ell\}}
 =  \q_0 + h \sum_{\ell=1}^{\nparts} \sum_{j=1}^{s^{\{\ell\}}} x_{i,j}^{\{m,\ell\}}\, b^{\{\ell\}}_j\,  \ell_j^{\{\ell\}} , \\
%\qquad i=1,\ldots,\hat s^{\{q\}},\\
\p_1 & =  \p_0 + h \sum_{\ell=1}^{\nparts} \sum_{i=1}^{\hat s^{\{\ell\}}} \widehat{b}_{i}^{\{\ell\}} k_i^{\{\ell\}}, \\
\q_1 & =  \q_0 + h \sum_{\ell=1}^{\nparts} \sum_{i=1}^{s^{\{\ell\}}} {b}_{i}^{\{\ell\}} \ell_i^{\{\ell\}}.
\end{align*}

{%\color{blue} 
After reordering the rows and columns such as to reflect the order in which stages are computed, we have the following cases.
\begin{enumerate}
\item If we compute all the stages $i=1,\dots,s^{\{m\}}$ for partition $m$ before moving on to partition $m+1$, then:
\begin{itemize}
\item both $\A^{\{m,\ell\}}$ and $\Ahat^{\{m,\ell\}}$ are zero for $m < \ell$; 
\item they are full matrices for  $m > \ell$, with $\A^{\{m,\ell\}} = \one^{\{m\}} \b^{\{\ell\}T}$ and $\Ahat^{\{m,\ell\}} = \one^{\{m\}} \bhat^{\{\ell\}T}$; and 
\item are lower triangular for $m = \ell$ with $a_{i,i}^{\{m,m\}} \cdot \ahat_{i,i}^{\{m,m\}} = 0$ for explicitness.
\end{itemize}
\item If we compute stage $i$ for each partition $m = 1,\dots,\nparts$ before moving on to stage $i+1$, then each of the coefficient matrices has to be lower triangular such as to preserve explicitness. This implies $x_{i,j}^{\{m,\ell\}} = 0$ and $1-x_{j,i}^{\{m,\ell\}} = 0$ for $i<j$, therefore $x_{i,j}^{\{m,\ell\}} = 1$ for $i>j$. We have the following structures when $s^{\{m\}} \ge s^{\{\ell\}}$:
\begin{align*}
   \scalebox{0.55}{$
    \A^{\{m,\ell\}}  = \begin{bmatrix}
    b^{\{\ell\}}_1\,x_{1,1}^{\{m,\ell\}} & \scriptstyle 0 & \cdots  & \scriptstyle 0  \\
    b^{\{\ell\}}_1 & b^{\{m\}}_2\,x_{2,2}^{\{m,\ell\}}   & &  \vdots \\
        \vdots &  \vdots  &\ddots  &  \vdots \\
        b^{\{\ell\}}_1 & b^{\{m\}}_2 &  \cdots &  b^{\{\ell\}}_{s^{\{\ell\}}} \,x_{s^{\{\ell\}},s^{\{\ell\}}}^{\{m,\ell\}} \\
        \vdots &  \vdots  &\ddots  &  \vdots \\
        b^{\{\ell\}}_1 & b^{\{\ell\}}_2 &  \cdots &  b^{\{\ell\}}_{s^{\{\ell\}}}
    \end{bmatrix}
    $}, ~~
    \scalebox{0.55}{$
    \Ahat^{\{\ell,m\}} = \begin{bmatrix}
    \widehat{b}^{\{m\}}_1\,\widehat{x}_{1,1}^{\{m,\ell\}} & \scriptstyle 0 & \cdots  & \scriptstyle 0 & \scriptstyle 0 & \dots & \scriptstyle 0  \\
    \widehat{b}^{\{m\}}_1 & \widehat{b}^{\{m\}}_2\,\widehat{x}_{2,2}^{\{m,\ell\}}   & &  \vdots & \vdots  && \vdots \\
        \vdots &  \vdots  &\ddots  &  \vdots & \vdots & & \vdots \\
        \widehat{b}^{\{m\}}_1 & \widehat{b}^{\{m\}}_2 &  \cdots &  \widehat{b}^{\{m\}}_{s^{\{\ell\}}} \,\widehat{x}_{s^{\{\ell\}},s^{\{\ell\}}}^{\{m,\ell\}} & \scriptstyle 0 & \dots & \scriptstyle 0  
    \end{bmatrix} $}.
\end{align*}
The last $s^{\{m\}}-s^{\{\ell\}}-1$ stages of $\A^{\{m,\ell\}} $ are redundant (equal to each other) and the last $s^{\{m\}}-s^{\{\ell\}}-1$ columns of $\Ahat^{\{m,\ell\}} $ are zero. To avoid redundancy, it makes sense to only consider $s^{\{m\}} = s^{\{\ell\}} + 1$ if $x_{s^{\{\ell\}},s^{\{\ell\}}}^{\{m,\ell\}}=0$, and $s^{\{m\}} = s^{\{\ell\}}$ if $x_{s^{\{\ell\}},s^{\{\ell\}}}^{\{m,\ell\}}=1$.
For $s^{\{m\}} \le s^{\{\ell\}}$ we get
\begin{align*}
   \scalebox{0.55}{$
    \A^{\{m,\ell\}}  = \begin{bmatrix}
  {b}^{\{m\}}_1\,{x}_{1,1}^{\{m,\ell\}} & \scriptstyle 0 & \cdots  & \scriptstyle 0 & \scriptstyle 0 & \dots & \scriptstyle 0  \\
    {b}^{\{m\}}_1 & {b}^{\{m\}}_2\,{x}_{2,2}^{\{m,\ell\}}   & &  \vdots & \vdots  && \vdots \\
        \vdots &  \vdots  &\ddots  &  \vdots & \vdots & & \vdots \\
        {b}^{\{m\}}_1 & {b}^{\{m\}}_2 &  \cdots &  {b}^{\{m\}}_{s^{\{\ell\}}} \,{x}_{s^{\{\ell\}},s^{\{\ell\}}}^{\{m,\ell\}} & \scriptstyle 0 & \dots & \scriptstyle 0  
    \end{bmatrix} 
    $}, ~~
    \scalebox{0.55}{$
    \Ahat^{\{\ell,m\}} = \begin{bmatrix}
    \widehat{b}^{\{\ell\}}_1\,\widehat{x}_{1,1}^{\{m,\ell\}} & \scriptstyle 0 & \cdots  & \scriptstyle 0  \\
    \widehat{b}^{\{\ell\}}_1 & \widehat{b}^{\{m\}}_2\,\widehat{x}_{2,2}^{\{m,\ell\}}   & &  \vdots \\
        \vdots &  \vdots  &\ddots  &  \vdots \\
        b^{\{\ell\}}_1 & b^{\{m\}}_2 &  \cdots &  b^{\{\ell\}}_{s^{\{\ell\}}} \,\widehat{x}_{s^{\{\ell\}},s^{\{\ell\}}}^{\{m,\ell\}} \\
        \vdots &  \vdots  &\ddots  &  \vdots \\
        \widehat{b}^{\{\ell\}}_1 & \widehat{b}^{\{\ell\}}_2 &  \cdots &  \widehat{b}^{\{\ell\}}_{s^{\{\ell\}}}
    \end{bmatrix}
    $},
    \end{align*}
    Similar to the first case, the last $s^{\{\ell\}}-s^{\{m\}}-1$ stages of $\Ahat^{\{m,\ell\}} $ are redundant (equal to each other) and the last $s^{\{\ell\}}-s^{\{m\}}-1$ columns of $\A^{\{m,\ell\}} $ are zero. To avoid redundancy, it makes again sense to only consider $s^{\{\ell\}} = s^{\{m\}} + 1$ if $x_{s^{\{\ell\}},s^{\{\ell\}}}^{\{m,\ell\}}=1$, and $s^{\{m\}} = s^{\{\ell\}}$ if $x_{s^{\{\ell\}},s^{\{\ell\}}}^{\{m,\ell\}}=0$.

\end{enumerate}

\begin{remark}
For $|s^{\{m\}}-s^{\{\ell\}}| \ge 2$ the redundant stages may arise not only at the end for the last $|s^{\{m\}}-s^{\{\ell\}}| -1$ stages, but also before. The corresponding matrices are then no longer tridiagonal, but have a step form, i.e., if one element in a row is zero, all elements above are zero, too. 

One example for such a setting will be given in Example~\ref{ex.extended.yoshida} for the matrices $\Ahat^{\{1,3\}}$ and $\A^{\{3,1\}}$.
 \end{remark}
}

\begin{leaveout}
%
%MG: I do not understand this example
\begin{remark}
Consider the scheme \eqref{eq:sympl.part.gark} and restrict coefficients such as to make it explicit:
\begin{subequations}
\label{eq:sympl.part.gark.explicit}
\begin{align*}
%P_i^{\{m\}} & =  \p_0 + h \sum_{q=1}^{m-1} \sum_{j=1}^{\hat s^{\{q\}}} \ahat_{i,j}^{\{m,q\}} k_j^{\{q\}}
%+ h  \sum_{j=1}^{i-1} \ahat_{i,j}^{\{m,m\}} k_j^{\{m\}} \\
%& \quad + h \,\ahat_{i,i}^{\{m,m\}} k_j^{\{m\}} + {\color{blue} h \sum_{q=m+1}^{\nparts} \sum_{j=1}^{i-1} \ahat_{i,j}^{\{m,q\}} k_j^{\{q\}}}, \dots \\
%
P_i^{\{q\}} & =  \p_0 + h \sum_{m=1}^{q-1} \sum_{j=1}^{\hat s^{\{m\}}} \ahat_{i,j}^{\{q,m\}} k_j^{\{m\}}
+ h  \sum_{j=1}^{i-1} \ahat_{i,j}^{\{q,q\}} k_j^{\{q\}} \\
& \quad + h \,\ahat_{i,i}^{\{q,q\}} k_j^{\{q\}} + h \sum_{m=q+1}^{\nparts} \sum_{j=1}^{i-1} \ahat_{i,j}^{\{q,m\}} k_j^{\{m\}} \\
& =  \p_0 + h \sum_{m=1}^{q-1} \sum_{j=1}^{\hat s^{\{m\}}} 
\widehat{b}^{\{m\}}_j\,(1-x_{j,i}^{\{m,q\}}) k_j^{\{m\}}
+ h  \sum_{j=1}^{i-1} 
\widehat{b}^{\{q\}}_j\,(1-x_{j,i}^{\{q,q\}}) k_j^{\{q\}}
\\
& \quad + h \,\widehat{b}^{\{q\}}_i\,(1-x_{i,i}^{\{q,q\}}) 
k_j^{\{q\}} + h \sum_{m=q+1}^{\nparts} \sum_{j=1}^{i-1} 
\widehat{b}^{\{m\}}_j\,(1-x_{j,i}^{\{m,q\}}) 
k_j^{\{m\}} \\
%P_i^{\{1\}} & =  \p_0 
%+ h  \sum_{j=1}^{i-1} \ahat_{i,j}^{\{1,1\}} k_j^{\{1\}} + h \,\ahat_{i,i}^{\{1,1\}} k_i^{\{1\}} + {\color{blue} h \sum_{j=1}^{i-1} \ahat_{i,j}^{\{1,2\}} k_j^{\{2\}}}, \dots \\
%%
%Q_i^{\{1\}} & =  \q_0 
%+ h  \sum_{j=1}^{i-1}  a_{i,j}^{\{1,1\}} \ell_j^{\{1\}} + h \, a_{i,i}^{\{1,1\}} \ell_i^{\{1\}}+ h  \sum_{j=1}^{i-1}  a_{i,j}^{\{1,2\}} \ell_j^{\{2\}},\\
%%
%P_i^{\{2\}} & =  \p_0 + h \sum_{j=1}^{i} \ahat_{i,j}^{\{2,1\}} k_j^{\{1\}}
%+ h  \sum_{j=1}^{i-1} \ahat_{i,j}^{\{2,2\}} k_j^{\{2\}} + h \,\ahat_{i,i}^{\{2,2\}} k_i^{\{2\}}, \\
%%
%Q_i^{\{2\}} & =  \q_0 +  {\color{blue} h \sum_{j=1}^{i}  a_{i,j}^{\{2,1\}} \ell_j^{\{1\}} }
%+ h  \sum_{j=1}^{i-1}  a_{i,j}^{\{2,2\}} \ell_j^{\{2\}} + h \, a_{i,i}^{\{2,2\}} \ell_i^{\{2\}}.
%
Q_i^{\{q\}} & =  \q_0 + {%\color{purple} 
h \sum_{m=1}^{q-1} \sum_{j=1}^{s^{\{m\}}}  a_{i,j}^{\{q,m\}} \ell_j^{\{m\}} }
+ h  \sum_{j=1}^{i-1}  a_{i,j}^{\{q,q\}} \ell_j^{\{q\}} \\
& \quad + h \, a_{i,i}^{\{q,q\}} \ell_j^{\{q\}}+ h \sum_{m=q+1}^{\nparts} \sum_{j=1}^{i-1}  a_{i,j}^{\{q,m\}} \ell_j^{\{m\}} \\
& =  \q_0 + {%\color{purple} 
h \sum_{m=1}^{q-1} \sum_{j=1}^{s^{\{m\}}}   
b^{\{m\}}_j\,x_{i,j}^{\{q,m\}}
\ell_j^{\{m\}} }
+ h  \sum_{j=1}^{i-1}  
b^{\{q\}}_j\,x_{i,j}^{\{q,q\}}
\ell_j^{\{q\}} \\
& \quad + h \, 
b^{\{q\}}_j\,x_{i,i}^{\{q,q\}}
\ell_j^{\{q\}}+ h \sum_{m=q+1}^{\nparts} \sum_{j=1}^{i-1}  
b^{\{m\}}_j\,x_{i,j}^{\{q,m\}}
\ell_j^{\{m\}}.
\end{align*}
\end{subequations}
Here $\A^{\{q,q\}}$ and $\Ahat^{\{q,q\}}$ are lower triangular, with $a_{i,i}^{\{q,q\}} \cdot \ahat_{i,i}^{\{q,q\}} = 0$ for explicitness. {%\color{purple} 
$\A^{\{q,m\}}$ for $m < q$ is full}, and $\Ahat^{\{q,m\}}$ is lower triangular, with
\begin{subequations}
\begin{align*}
    \A^{\{q,m\}} & = \bf{x}^{\{q,m\}} \cdot \diag{\b^{\{m\}}}, \\
    \Ahat^{\{q,m\}} & = \left( \ones \ones\tr - \bf{x}^{\{m,q\},\top}\right) \cdot \diag{\b^{\{m\}}}.
\end{align*}
\end{subequations}
\end{remark}
\end{leaveout}

\subsubsection*{Explicit symplectic and symmetric partitioned GARK schemes} If, in addition, the scheme is also symmetric then condition~ \eqref{eqn:GARK-symmetry}  leads to
\begin{align*}
a_{j,i}^{\{m,\ell\}} + a_{s^{\{m\}}+1-j,s^{\{\ell\}}+1-i}^{\{m,\ell\}} &=  b_{i}^{\{\ell\}}   =  b_{s^{\{\ell\}}+1-i}^{\{\ell\}}, \\
\widehat{a}_{i,j}^{\{\ell,m\}} + \widehat{a}_{s^{\{\ell\}}+1-i,s^{\{m\}}+1-j}^{\{\ell,m\}} &=  \widehat{b}_{j}^{\{m\}}   =  \widehat{b}_{s^{\{m\}}+1-j}^{\{m\}}.
\end{align*}
Since $a_{j,i}^{\{m,\ell\}}$ can take either values $0$ or $b_{i}^{\{\ell\}}$, $a_{s^{\{m\}}+1-j,s^{\{\ell\}}+1-i}^{\{m,\ell\}}$ takes the complementary values $b_{s^{\{\ell\}}+1-i}^{\{\ell\}}$ or $0$, respectively.
All possible solutions of equation \eqref{eq:explicit-conditions} have the form:
\begin{equation}
\label{eq:explicit-conditions.symmetric}
\begin{split}
 \ahat^{\{\ell,m\}}_{i,j} &= \widehat{b}^{\{m\}}_j \,(1-x_{j,i}^{\{m,\ell\}}), \\
 \widehat{a}_{s^{\{\ell\}}+1-i,s^{\{m\}}+1-j}^{\{\ell,m\}} &= \widehat{b}_{s^{\{m\}}+1-j}^{\{m\}} \, x_{j,i}^{\{m,\ell\}}, \\
a^{\{m,\ell\}}_{j,i} &=  b^{\{\ell\}}_i\, x_{j,i}^{\{m,\ell\}}, \\
a_{s^{\{m\}}+1-j,s^{\{\ell\}}+1-i}^{\{m,\ell\}} &=  b_{s^{\{\ell\}}+1-i}^{\{\ell\}} (1-x_{j,i}^{\{m,\ell\}}),
\end{split}
\qquad
x_{j,i}^{\{m,\ell\}} \in \{0,1\}.
%\begin{cases}
%\textnormal{either} &  \ahat^{\{\ell,m\}}_{i,j} = \widehat{b}^{\{m\}}_j   \mbox{ and }  \widehat{a}_{s^{\{\ell\}}+1-i,s^{\{m\}}+1-j}^{\{\ell,m\}} = 0  \\
%& \mbox{ and } a^{\{m,\ell\}}_{j,i}=0
% \mbox{ and } a_{s^{\{m\}}+1-j,s^{\{\ell\}}+1-i}^{\{m,\ell\}}=  b_{s^{\{\ell\}}+1-i}^{\{\ell\}}, \\[3pt]
%\textnormal{or} & 
%\ahat^{\{\ell,m\}}_{i,j}=0  \mbox{ and } \widehat{a}_{s^{\{\ell\}}+1-i,s^{\{m\}}+1-j}^{\{\ell,m\}}  =  \widehat{b}_{s^{\{m\}}+1-j}^{\{m\}} \\
%&  \mbox{ and }  a^{\{m,\ell\}}_{j,i} = b^{\{\ell\}}_i
%  \mbox{ and }  a_{s^{\{m\}}+1-j,s^{\{\ell\}}+1-i}^{\{m,\ell\}} = 0.
%\end{cases}
\end{equation}

We finish with two examples for symplectic and symmetric GARK schemes.
%related to the second case above. \sandu{Which is the first and which is the second case?}

\begin{example}[Yoshida~\cite{Yoshida_1990_splitting}]
\label{ex.classical.Yoshida}
The classical fourth order symplectic and symmetric scheme of Yoshida can be written as a partitioned GARK scheme with $\nparts=1$ and
\begin{align*}
    \A^{\{1,1\}} = \begin{bmatrix}
        \frac{d_1}{2} & \scriptstyle 0 & \scriptstyle 0 & \scriptstyle 0  \\
        \frac{d_1}{2} & \frac{d_1+d_2}{2} & \scriptstyle 0 & \scriptstyle 0  \\
        \frac{d_1}{2} & \frac{d_1+d_2}{2} & \frac{d_1+d_2}{2} & \scriptstyle 0 
        \end{bmatrix},
        \quad &
        \Ahat^{\{1,1\}}= \begin{bmatrix}
        \scriptstyle 0 & \scriptstyle 0 & \scriptstyle 0  \\
        \scriptstyle d_1 & \scriptstyle 0 & \scriptstyle 0   \\
        \scriptstyle d_1 & \scriptstyle d_2 & \scriptstyle 0  \\
        \scriptstyle d_1 & \scriptstyle d_2 & \scriptstyle d_1
        \end{bmatrix},  \\
        \b=\left[ \sfrac{d_1}{2}, \sfrac{d_1+d_2}{2}, \sfrac{d_1+d_2}{2},\sfrac{d_1}{2}\right]\tr, \quad &
        \bhat=\left[d_1,d_2,d_1\right]\tr, \quad 
        d_1=\sfrac{1}{2-2^{\frac{1}{3}}},~~ d_2=-2^{\frac{1}{3}} \cdot d_1.
\end{align*}
\end{example}

\begin{example}[Extension of Yoshida's scheme~\cite{Yoshida_1990_splitting}]
\label{ex.extended.yoshida}
An explicit partitioned symmetric and symplectic scheme of type~\eqref{eqn:GARK-symplectic.partitioned.separable} for  potential splitting with $\nparts=3$
is given by the following extension of Yoshida's fourth order scheme~\cite{Yoshida_1990_splitting}:
\begin{equation*}
\renewcommand{\arraystretch}{1.5}
%\label{eqn:Partitioned-Butcher-tableau2}
\begin{array}{c|c|c}
\Zero_{4 \times 4}   & 
         \begin{array}{ccc}
         \scriptstyle 0 & \scriptstyle 0 & \scriptstyle 0  \\
        \scriptstyle d_1 & \scriptstyle 0 & \scriptstyle 0   \\
        \scriptstyle d_1 & \scriptstyle d_2 & \scriptstyle 0  \\
        \scriptstyle d_1 & \scriptstyle d_2 & \scriptstyle d_1
        \end{array} &
        \begin{array}{cc}
        \scriptstyle 0 & \scriptstyle 0  \\
        \frac{1}{2} & \scriptstyle 0  \\
        \frac{1}{2} & \scriptstyle 0  \\
        \frac{1}{2} & \frac{1}{2}
        \end{array}  \\
\hline
\begin{array}{cccc}
\frac{d_1}{2} & \scriptstyle 0 & \scriptstyle 0 & \scriptstyle 0  \\
        \frac{d_1}{2} & \frac{d_1+d_2}{2} & \scriptstyle 0 & \scriptstyle 0  \\
        \frac{d_1}{2} & \frac{d_1+d_2}{2} & \frac{d_1+d_2}{2} & \scriptstyle 0 
        \end{array} 
        & \Zero_{3 \times 3} & \Zero_{3 \times 2}  \\ \hline 
   \begin{array}{cccc}
        \frac{d_1}{2} & \scriptstyle 0 & \scriptstyle 0 & \scriptstyle 0  \\
        \frac{d_1}{2} & \frac{d_1+d_2}{2} & \frac{d_1+d_2}{2} & \scriptstyle 0 
        \end{array} 
        & \Zero_{2 \times 3}  & \Zero_{2 \times 2}  \\
 \Xhline{2\arrayrulewidth}
        \begin{array}{cccc}
        \frac{d_1}{2} & \frac{d_1+d_2}{2} & \frac{d_1+d_2}{2} & \frac{d_1}{2} 
        \end{array} &
        \begin{array}{ccc}
        \scriptstyle d_1 & \scriptstyle d_2 & \scriptstyle d_1
        \end{array} &
        \begin{array}{cc}
        \frac{1}{2} & \frac{1}{2}
        \end{array} 
        \end{array}
\end{equation*}
Note that this scheme has order four for $H_1(\p,\q)=T(\p)+V_1(\q)$ and order two for $H_2(\p,\q)=V_2(\q)$. This multi-order character of the scheme is tailored to exploiting a fast dynamics and cheap evaluation costs in $H_1$, and  a slow dynamics and expensive evaluation costs in $H_2$.
\end{example}

\section{Numerical examples}
We finish with two examples for symplectic and time-reversible GARK schemes: the kdV equation as an example for a general skew-symmetric matrix $\mathbf{J}$ discussed in Section~\ref{sec:symplectic-GARK}, and a mathematical pendulum as an example for multirate potential in a potential splitting discussed in Section~\ref{sec-partitionedsymplecticgark}.

%%%%%%%%%%%%%%%%%%%%%%%%%%
\subsection{Symplectic integration with non-symplectic partitions} 
\label{sec-non-symplectic-splitting}
%%%%%%%%%%%%%%%%%%%%%%%%%%

\begin{figure}[ht!]
	\centering
	        \subfloat[Solution at different times]{
		\centering
		\includegraphics[width=0.45\textwidth]{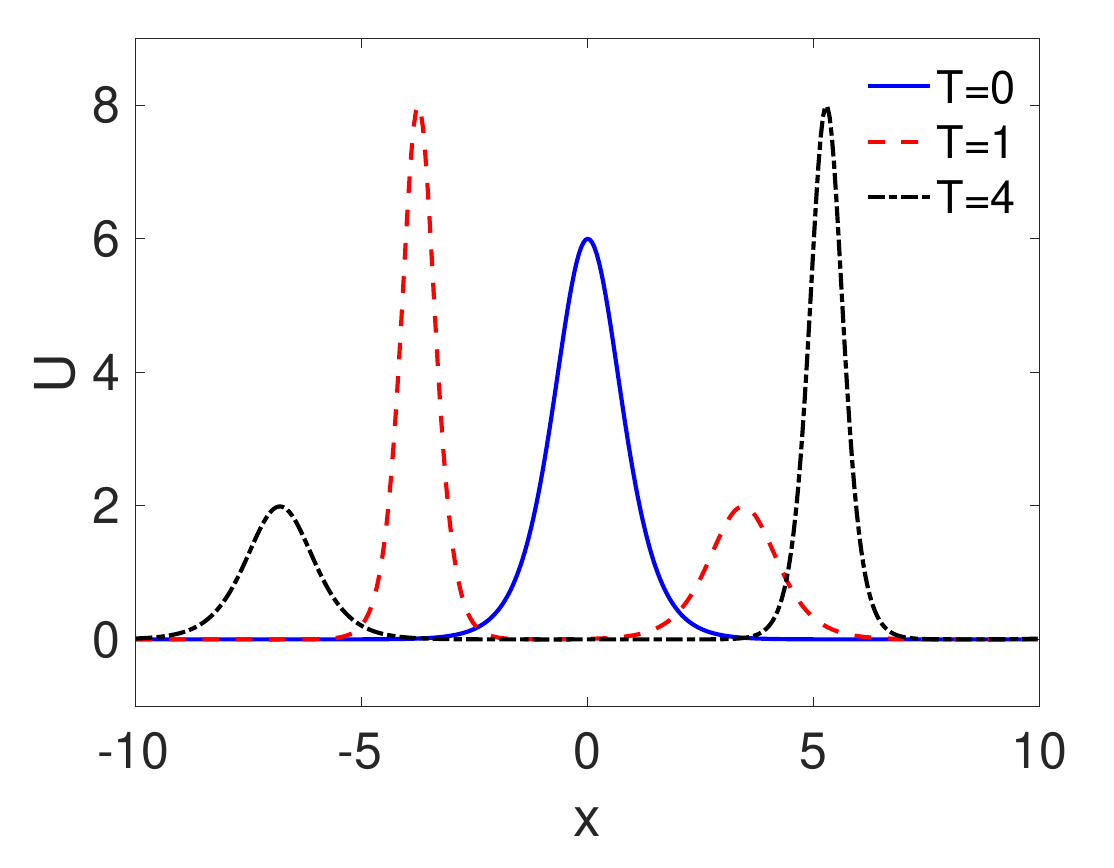}
                 \label{fig:KdV-solution}
	}
	\hfill
		\subfloat[Hamiltonian error evolution]{
		\centering
		\includegraphics[width=0.45\textwidth]{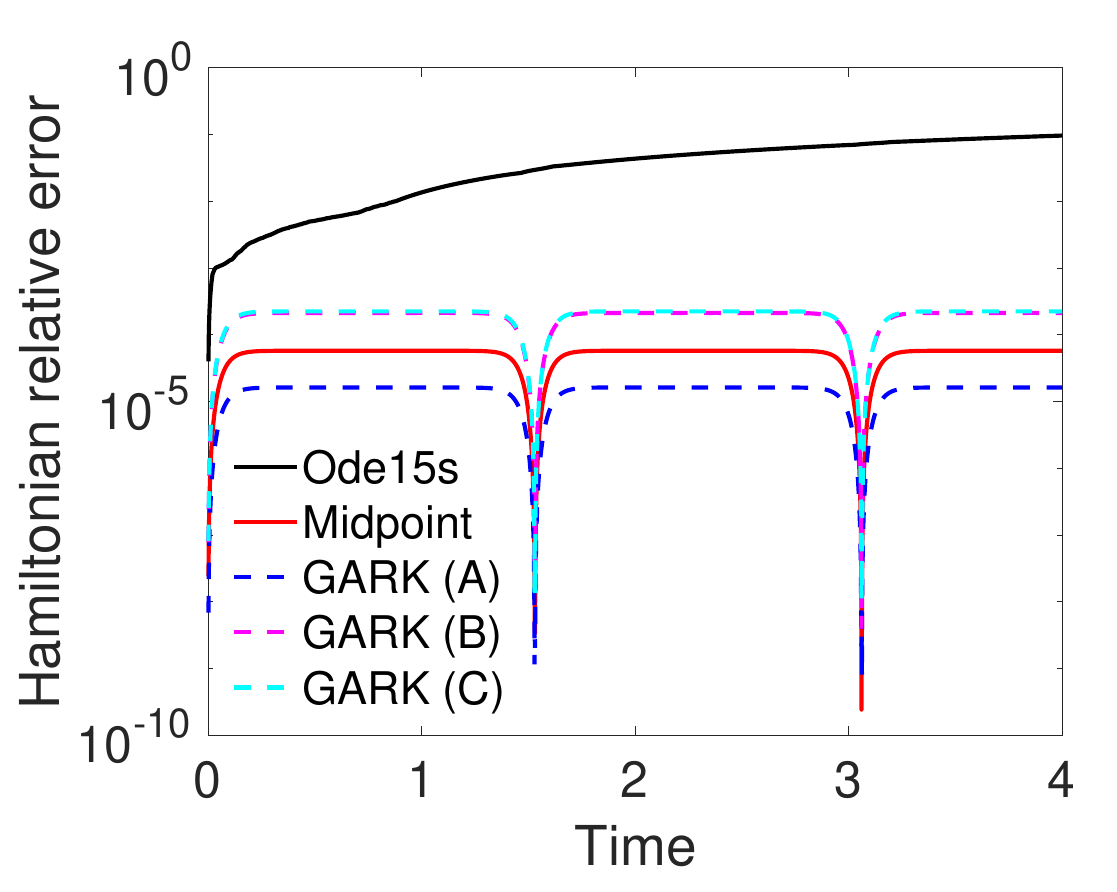}
                 \label{fig:KdV-hamiltonian}
	}
    \caption{Numerical results for the KdV system \eqref{eqn:kdv-numerical} solved with different time integration methods.}
    \label{fig:KdV}
\end{figure}

We consider symplectic time integration for the Korteweg-de Vries (KdV) equation \cite{Ascher_2005_KdV,Dutykh_2013_KdV}, a non-dissipative nonlinear hyperbolic equation with smooth solutions:
\begin{equation}
\label{eqn:kdv}
\begin{split}
u' &= \alpha\,(u^2)_x + \rho\,u_x + \nu\,u_{xxx} = V'(u)_x +  \nu\,u_{xxx}, \\
V(u) &= \sfrac{\alpha}{3}\,u^3 + \sfrac{\rho}{2}\,u^2, \quad u(t=0,x) = 6\,\operatorname{sech}(x)^2,
\quad \alpha = -3, ~ \rho = 1, ~\nu = -1,
\end{split}
\end{equation}
and periodic boundary conditions $u(0,t)=u(10,t)$. The initial condition leads to the formation of two solitons traveling at different speeds \cite{Dutykh_2013_KdV}, as seen in Figure \ref{fig:KdV-solution}.
The discrete Hamiltonian leads to a symplectic semi-discretization in space:
\begin{equation}
\label{eqn:kdv-numerical}
\begin{split}
H(\mathbf{u}) &= \Delta x \, \sum_{i} \left( V(u_i) - \sfrac{\nu}{2}\, \left(\sfrac{u_{i+1}-u_{i}}{\Delta x}\right)^2 \right), \quad
u_i' =  \sfrac{1}{2\, \Delta x}\, \left( \sfrac{\partial H}{\partial u_{i+1}} - \sfrac{\partial H}{\partial u_{i-1}} \right), \\
%    \Leftrightarrow ~\mathbf{u}' &= \fun(\mathbf{u}) =\alpha\, D_1\,\mathbf{u}^2  + \rho\,D_1\,\mathbf{u} + \nu\,D_3\,\mathbf{u},\\
%& \qquad D_1 = \sfrac{[-1, 0,1]}{2\Delta x}\,, \quad D_3 = \sfrac{[-1, 2, 0, -2, 1]}{\Delta x^3}.
%%[-2, -1, 1, 2]
\end{split}
\end{equation}
{which can be written as a generalized Hamiltonian system with 
\begin{equation}
\label{kdV2}
    \mathbf{u}' = \mathbf{J} \cdot \nabla H(\mathbf{u})
\end{equation}
and the skew-symmetric matrix $J$ given by
\[
\mathbf{J}= \begin{bmatrix}
0 & 1      &       &        \\
-1 & \ddots & \ddots &        \\
   & \ddots      & \ddots & 1  & \\
   & & -1 & 0 
\end{bmatrix} - e_1 e_n^\top + e_n e_1^\top, \qquad n \coloneqq  10/\Delta x.
\]
}
%}
%
\begin{leaveout}
\begin{remark}
Note that system~\eqref{eqn:kdv-numerical} has the structure of a generalized Hamiltonian system with 
\begin{equation}
    \mathbf{u}' = J \cdot \nabla H(\mathbf{u})
\end{equation}
with the skew-symmetric matrix $J$ given by
\[
J= \begin{bmatrix}
0 & 1      &       &        \\
-1 & \ddots & \ddots &        \\
   & \ddots      & \ddots & 1  & \\
   & & -1 & 0 
\end{bmatrix} - e_1 e_n^\top + e_n e_1^\top, \qquad n \coloneqq  10/\Delta x.
\]
\end{remark}
\end{leaveout}
We integrate the system \eqref{kdV2} %\eqref{eqn:kdv-numerical} 
with ode15s in Matlab, the symplectic implicit midpoint scheme, and with the symplectic GARK-IMIM scheme \eqref{eqn:GARK-IMIM2} using different partitions:
{
\begin{equation}
\label{eqn:kdv-partitions}
\fun^{\{1\}}(\mathbf{u})  =
\begin{cases}
\mathbf{J} \nabla H_1(\mathbf{u}), & (A) \\
\mathbf{J} \nabla H_2(\mathbf{u}), & (B) \\
\mathbf{J} \nabla (H_1(\mathbf{u})+H_2(\mathbf{u})),  & (C)
\end{cases}; \quad
\fun^{\{2\}}(\mathbf{u}) = \mathbf{J} \nabla H(\mathbf{u}) %\fun(\mathbf{u}) 
- \fun^{\{1\}}(\mathbf{u}).
\end{equation}
with
\begin{equation}
    H_1(\mathbf{u})= \Delta x \, \sum_{i} \frac{\sigma}{2} u_i^2, \quad
    H_2(\mathbf{u})= \Delta x \, \sum_{i} \frac{\alpha}{3} u_i^3.
\end{equation}
%
%
%\begin{equation}
%\label{eqn:kdv-partitions}
%\fun^{\{1\}}(\mathbf{u})  =
%\begin{cases}
%\rho\,D_1\,\mathbf{u}, & (A) \\
%\alpha\,D_1\,\mathbf{u}^2, & (B) \\
%\alpha\,D_1\,\mathbf{u}^2 + \rho\,D_1\,\mathbf{u},  & (C)
%\end{cases}; \quad
%\fun^{\{2\}}(\mathbf{u}) = \fun(\mathbf{u}) - \fun^{\{1\}}(\mathbf{u}).
%\end{equation}
%%
}
A fixed time step $\Delta t = 10^{-3}$ is used. Results are shown in Figure \ref{fig:KdV}. For all partitions \eqref{eqn:kdv-partitions} the GARK scheme is symplectic and thus preserves a nearby shadow Hamiltonian. Consequently, the error in the Hamiltonian oscillates around the true value as it can be seen in Figure \ref{fig:KdV-hamiltonian}.

%\textcolor{red}{Are partitions   \eqref{eqn:kdv-partitions} Hamiltonian? H is not conserved when I add and subtract a random term.}

%\sandu{
%Looks like it works (only) for partitions that are skew symmetric, even if not Hamiltonian:
%\begin{equation}
%   \mathbf{u}' = J \cdot \sum_m g^{\{m\}}(\mathbf{u}) = \sum_m \fun^{\{m\}}(\mathbf{u}),
%    \quad \fun^{\{m\}} = J \cdot g^{\{m\}}.
%\end{equation}
%Can we explain this?
%}

%%%%%%%%%%%%%%%%%%%%%%%%%%
\subsection{Symplectic and time-reversible GARK schemes for Hamiltonians with multirate potential} 
\label{sec-multiratepotential}
%%%%%%%%%%%%%%%%%%%%%%%%%%

Consider a Hamiltonian $H(\p,\q)=T(\p)+V(\q)$, where the potential can be split into two parts $V(\q)=V_1(\q)+V_2(\q)$. Assuming that $V_1$ is characterized by a fast dynamics and cheap evaluation costs, and  
$V_2$ by a slow dynamics and expensive evaluation costs, respectively. Then 
the Hamiltonian can be partitioned into two parts $H_1(\p,\q)+H_2(\p,\q)$ (with $H_1(\p,\q) \coloneqq T(\p)$, $H_2(\p,\q) \coloneqq V_1(\q)$)
%$H_1(\p,\q)=T(\p)+V_1(\q)$ 
and $H_3(\p,\q)=V_2(\q)$ with fast/slow dynamics and cheap/expensive evaluation costs, respectively.

\begin{figure}
\centering
%\centerline{\includegraphics[width=0.7\textwidth]{Figures/pendulum_adnotated.pdf}}
\setlength{\unitlength}{1cm}
\scalebox{0.7}{
\begin{picture}(20,5)
\put(4.5,3.2){\Large{$\alpha$}}
\put(4.5,4.8){\Large{$x_2$}}
\put(9.2,3.6){\Large{$x_1$}}
\put(3.,2.8){\Large{$g$}}
\put(6.5,0.8){\Large{$m_{\rm{\rm pend}}$}}
\put(9.2,0.8){\Large{$m_{\rm{\rm osc}}$}}
\put(6.,2.7){\Large{$\ell$}}
\put(8.15,2.0){\Large{$k$}} %{\Large{$k,\,d$}}
\put(2.5,0){\includegraphics[width=2.9in]{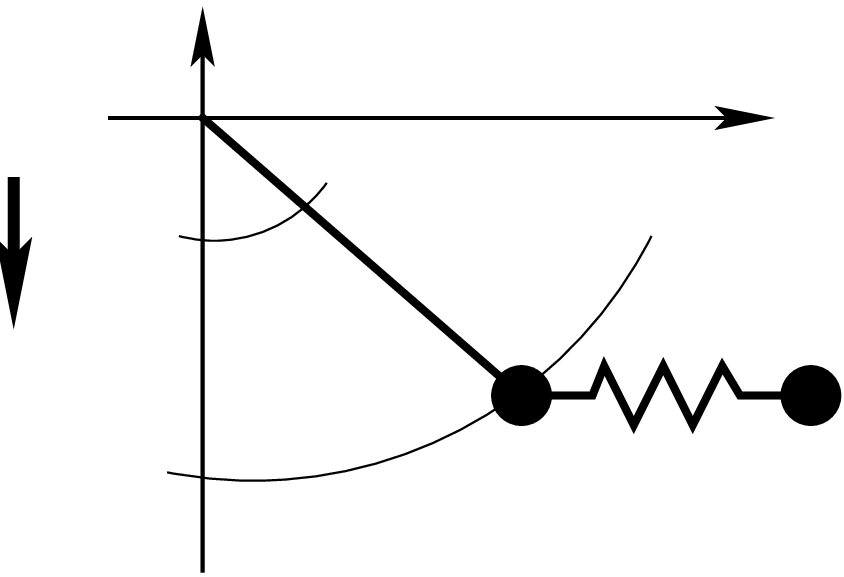}}
\end{picture}
}
\caption{Mathematical pendulum coupled to an oscillator (taken after~\cite{Arnold_2007_MR-mutibody}).
}
\label{fig:coupled-pendulum} 
\end{figure}

As an example of such a system with multiscale behaviour  we consider a mathematical pendulum of constant length $\ell$ that is coupled to a damped  oscillator with a horizontal degree of freedom, as illustrated in Figure~\ref{fig:coupled-pendulum}. 
The system consists of two rigid bodies: the first mass $m_{\rm{\rm pend}}$ is connected to a second mass $m_{\rm{\rm osc}}$ by a soft spring with stiffness $k$. Neglecting the friction of the spring, the system is Hamiltonian.

The minimal set of coordinates $\q^\top=(q_1,q_2) \coloneqq (\alpha,x_1)$  and generalized momenta $\p^\top=(p_1,p_2)$ uniquely describe the position and momenta of both bodies. The Hamiltonian of the system is given by;
\begin{align*}
    H(\p,\q)=H_1(\p,\q)+H_2(\p,\q)
\end{align*}
with the fast Hamiltonian
\begin{align*}
    H_1(\p,\q) & = T(\p) + V_1(\q), \\
    T(\p) &= \sfrac{1}{2\, m_{\rm osc}} p_2^2 + \sfrac{1}{2\, m_{\rm pend}} \left( \sfrac{p_1}{\ell} \right)^2, \\
    V_1(\q) & = - m_{\rm pend} \, g \, \ell \, \cos (q_1),
\end{align*}
and the slow Hamiltonian
\begin{align*}
    H_2(\p,\q) & = V_2(\q) = \frac{1}{2} k \big( q_2- \ell \sin (q_1) \big)^2.
\end{align*}
The equations of motion are then given by the second-order ODE system
\begin{eqnarray*}
\begin{pmatrix}
m_{\rm{\rm pend}}\, \ell & 0 \\ 0 & m_{\rm{\rm osc}}
\end{pmatrix} \,
\ddot{q} 
= 
\begin{pmatrix}
- m_{\rm{\rm pend}}\, g\, \sin(\alpha) + \cos(\alpha)\, F \\
-F
\end{pmatrix}
\eqqcolon f(\q),
\end{eqnarray*}
where the following abbreviation stands for the spring force:

\[
F= k \, \big(x_1 - \ell\, \sin(\alpha) \big). % + d \, \big(\dot x_1 - \ell\, \dot \alpha\, \cos (\alpha) \big),
\]
%

\begin{figure}
    %\centering
    \includegraphics{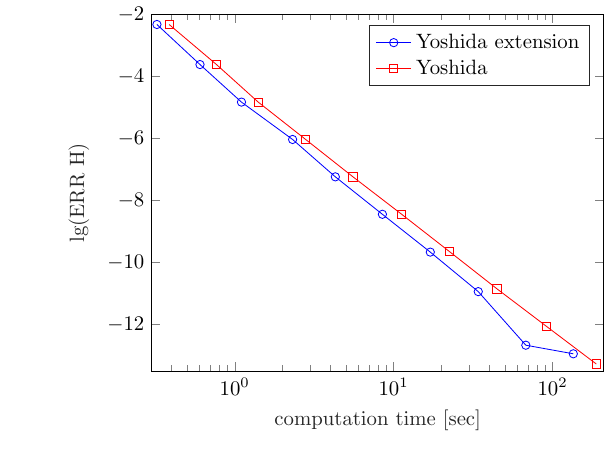}
    \caption{Numerical results for parameters $m_{\rm{\rm pend}} = m_{\rm{\rm osc}} = \ell = 1$, $k=5 \cdot 10^{-6}$, and $10000$ times higher evaluation costs for $V_2$: absolute error in the Hamiltonian $H$ for Yoshida and the Yoshida extension vs. computation time.}
    \label{fig.results.mbs2}
\end{figure}

%Figures~\ref{fig.results.mbs} and~
Figure~\ref{fig.results.mbs2} shows the numerical results obtained for this benchmark for the GARK extension of Yoshida's fourth order scheme derived in Example \ref{ex.extended.yoshida} and, for comparison, Yoshida's fourth order method from Example \ref{ex.classical.Yoshida}. 
Note that per integration step Yoshida's scheme needs three function evaluations of both $V_1$ and $V_2$, whereas the extension needs three for $V_1$, but only two for $V_2$. Figure~\ref{fig.results.mbs2} shows the achieved accuracy compared to the number of $V_2$ evaluations
%. A comparison of the computation times is given in~\ref{fig.results.mbs2}
assuming that the evaluation costs of $V_2$ are 10,000 times higher than the ones of $V_1$. In this case, the extension clearly outperforms the basic scheme of Yoshida. This  situation in typical for many problems with a fast but cheap and slow but expensive force as in Lattice Quantum Chromodynamics, for example, with a cheap gauge field with fast dynamics and an expensive fermionic force with slow dynamics~\cite{BookLQCD}.

%\begin{figure}
%    %\centering
%    \includegraphics{ctc_H_k=5e-6.pdf}
%    \caption{Numerical results for parameters $m_{\rm{\rm pend}} = m_{\rm{\rm osc}} = \ell = 1$ and $k=5 \cdot 10^{-6}$: absolute error in the Hamiltonian $H$ for Yoshida and the Yoshida extension vs. computation time.}
%    \label{fig.results.mbs}
%\end{figure}

\begin{comment}
\input{results/err_H_k-small.tex}

\input{results/rtc_H_k-small.tex}

The results can be explained as follows: the symplectic GARK scheme yields an error expansion 
with an order two term due to $V_2$ and an order four term due to $V_1$,
as the Hamiltonian $H_1$ and $H_2$ are solved with fourth and second order accuracy, resp.
For large step sizes %$h>\sqrt{C_2/C_1}$, 
the fourth error term dominates, for smaller step sizes the second order term. This explains the slopes 4 and 2 visible in Fig.~\ref{fig.results.mbs} for the symplectic GARK scheme for large and small step sizes, respectively.
As one can see in Fig.~\ref{fig.results.mbs2}, the symplectic GARK scheme becomes advantageous for higher accuracy demands in terms of computation time. Compared to 
the fourth-order scheme of Yoshida~\cite{Yoshida_1990_splitting}, which has to use three function evaluations per step, the symplectic GARK scheme only needs two function evaluations for the slow, but expensive part $V_2$ per step.
\end{comment}

%%%%%%%%%%%%%%%%%%%%%%%%%%
\section{Conclusions}
\label{sec-conclusions}
%%%%%%%%%%%%%%%%%%%%%%%%%%

This paper derives partitioned symplectic schemes in the GARK framework, which allows for arbitrary splittings of the Hamiltonian into different Hamiltonian subsystems,  which works also in the case of a more general Hamiltonian flow $\fun(\y)=\mathbf J \nabla H(\y)$ with an arbitrary, but skew-symmetric matrix $\mathbf J=- \mathbf J^\top$. The derived symplecticity conditions reduce drastically the number of GARK order conditions. We show that symmetric GARK schemes are time-reversible and construct symmetric and time-reversible GARK schemes based on composing a symplectic GARK scheme and its time-reversed scheme. A special attention is given to partitioned symplectic GARK schemes, which can be tailored to a specific splitting w.r.t. potentials or potentials and kinetic parts, resp. We show that symplecticity and self-adjointness are equivalent, and  show how the coupling matrices $\A^{\{\ell,m\}}$ and ${\hat \A}^{\{\ell,m\}}$ can  be chosen such as to construct explicit schemes. Using different discretization orders for different parts of the splitting defines one way to exploit the multiscale behavior of different potentials $V_1$ and $V_2$ of a Hamiltonian, where  $V_1$ is  characterized  by  a  fast  dynamics  and  cheap  evaluation  costs,  and $V_2$ by  a  slow dynamics and expensive evaluation costs, respectively. Numerical tests for a coupled oscillator confirm the theoretical results.

Future work will be to derive efficient symplectic GARK schemes tailored for couplings arising in port-Hamiltonian modeling on the one hand, and to generalize symplectic GARK schemes to multirate symplectic GARK schemes, which use different step sizes for different partitions to exploit the multirate potential. Another task will be to generalize this Abelian setting to a Non-Abelian setting used in lattice QCD, for example, where the equations of motion are defined on Lie groups and their associated Lie algebras.


\begin{thebibliography}{10}

\bibitem{armusa97}
{\sc A.~L. Ara{\'u}jo, A.~Murua, and J.~M. Sanz-Serna}, {\em Symplectic methods
  based on decompositions}, SIAM Journal on Numerical Analysis, 34 (1997),
  pp.~1926--1947.

\bibitem{Arnold_2007_MR-mutibody}
{\sc M.~Arnold}, {\em Multi-rate time integration for large scale multibody
  system models}, in IUTAM Symposium on Multiscale Problems in Multibody System
  Contacts, Eberhard P., ed., IUTAM Bookseries, vol.1, Springer, Dordrecht,
  2007, pp.~1--10.

\bibitem{Ascher_2005_KdV}
{\sc U.~M. Ascher and R.~I. McLachlan}, {\em On symplectic and multisymplectic
  schemes for the kdv equation}, Journal of Scientific Computing, 25 (2005),
  pp.~83--104.

\bibitem{paper_BGJR}
{\sc A.~Bartel, M.~G{\"u}nther, B.~Jacob, and T.~Reis}, {\em Operator
  {Splitting} {Based} {Dynamic} {Iteration} for {Linear Port-Hamiltonian
  Systems}}, Numerische Mathematik, 155 (2023), pp.~1--34.

\bibitem{Blanes_2003_composition}
{\sc S.~Blanes and P.C. Moan}, {\em Practical symplectic partitioned
  runge--kutta and runge--kutta--nystr{\"o}m methods}, Journal of Computational
  and Applied Mathematics, 142 (2002), pp.~313--330.

\bibitem{Dutykh_2013_KdV}
{\sc D.~Dutykh, M.~Chhay, and F.~Fedele}, {\em Geometric numerical schemes for
  the {KdV} equation}, Computational Mathematics and Mathematical Physics, 53
  (2013), pp.~221--236.

\bibitem{BookLQCD}
{\sc M.~G\"unther F.~Knechtli and M.~Peardon}, {\em Lattice Quantum
  Chromodynamics: Practical Essentials (SpringerBriefs in Physics)},
  Springer-Verlag, 2017.

\bibitem{Gonzalez1996TimeIA}
{\sc Oscar Gonzalez}, {\em Time integration and discrete hamiltonian systems},
  Journal of Nonlinear Science, 6 (1996), pp.~449--467.

\bibitem{Sandu_2022_GARK_splittingSchemes}
{\sc Severiano Gonz{\'a}lez-Pinto, Domingo Hern{\'a}ndez-Abreu, Maria~S.
  P{\'e}rez-Rodr{\'\i}guez, Arash Sarshar, Steven Roberts, and Adrian Sandu},
  {\em A unified formulation of splitting-based implicit time integration
  schemes}, Journal of Computational Physics, 448 (2022), p.~110766.

\bibitem{Sandu_2016_MR-GARK}
{\sc M.~G\"{u}nther and A.~Sandu}, {\em Multirate generalized additive
  {Runge-Kutta} methods}, Numerische Mathematik, 133 (2016), pp.~497--524.

\bibitem{Hager_2000_RKadj}
{\sc W.~Hager}, {\em {Runge-Kutta methods in optimal control and the
  transformed adjoint system}}, Numerische Mathematik, 87 (2000), pp.~247--282.

\bibitem{Hairer_2006_geometric-book}
{\sc Ernst Hairer, Christian Lubich, and Gerhard Wanner}, {\em Geometric
  numerical integration: structure-preserving algorithms for ordinary
  differential equations}, vol.~31, Springer Science \& Business Media, 2006.

\bibitem{Hairer_book_I}
{\sc E.~Hairer, S.P. Norsett, and G.~Wanner}, {\em Solving ordinary
  differential equations {I}: {N}onstiff problems}, no.~8 in Springer Series in
  Computational Mathematics, Springer-Verlag Berlin Heidelberg, 1993.

\bibitem{Kennedy_2003}
{\sc Christopher~A. Kennedy and Mark~H. Carpenter}, {\em Additive
  {Runge--Kutta} schemes for convection--diffusion--reaction equations},
  Applied Numerical Mathematics, 44 (2003), pp.~139--181.

\bibitem{McLachlan}
{\sc R.~I. McLachlan, G.~R.~W. Quispel, and N.~Robidoux}, {\em Geometric
  integration using discrete gradients}, Phil. Trans. R. Soc., Serie A 357
  (19969), pp.~1021--1045.

\bibitem{Sandu_2021_MR-GARK_Implicit}
{\sc S.~Roberts, J.~Loffeld, A.~Sarshar, C.S. Woodward, and A.~Sandu}, {\em
  Implicit multirate {GARK} methods}, Journal of Scientific Computing, 87
  (2021), p.~4.

\bibitem{Sandu_2020_MRI-GARK_Coupled}
{\sc S.~Roberts, A.~Sarshar, and A.~Sandu}, {\em Coupled multirate
  infinitesimal {GARK} methods for stiff differential equations with multiple
  time scales}, SIAM Journal on Scientific Computing, 42 (2020),
  pp.~A1609--A1638.

\bibitem{Sandu_2021_GARK-adjoint}
{\sc U.~Romer, M.~Narayanamurthi, and A.~Sandu}, {\em Goal-oriented a
  posteriori estimation of numerical errors in the solution of multiphysics
  systems}.
\newblock Submitted, 2021.

\bibitem{Sandu_2006_dadjRK}
{\sc A.~Sandu}, {\em On the properties of {Runge-Kutta} discrete adjoints}, in
  Lecture Notes in Computer Science, vol.~LNCS 3994, Part IV, International
  Conference on Computational Science, 2006, pp.~550--557.

\bibitem{Sandu_2019_MRI-GARK}
\leavevmode\vrule height 2pt depth -1.6pt width 23pt, {\em A class of multirate
  infinitesimal {GARK} methods}, SIAM Journal on Numerical Analysis, 57 (2019),
  pp.~2300--2327.

\bibitem{Sandu_2021_GARK-ROS}
{\sc A.~Sandu, M.~Guenther, and S.B. Roberts}, {\em Linearly implicit {GARK}
  schemes}, Applied Numerical Mathematics, 161 (2021), pp.~286--310.

\bibitem{Sandu_2015_GARK}
{\sc A.~Sandu and M.~G\"{u}nther}, {\em A generalized-structure approach to
  additive {Runge-Kutta} methods}, SIAM Journal on Numerical Analysis, 53
  (2015), pp.~17--42.

\bibitem{SanzSerna_2016_adjointSymplectic}
{\sc J.M. Sanz-Serna}, {\em Symplectic runge--kutta schemes for adjoint
  equations, automatic differentiation, optimal control, and more}, SIAM
  Review, 58 (2016), pp.~3--33.

\bibitem{SanzSerna_1992_symplectic-review}
{\sc J.~M. Sanz-Serna}, {\em Symplectic integrators for hamiltonian problems:
  an overview}, Acta Numerica, 1 (1992), pp.~243--286.

\bibitem{SanzSerna_1992_symplectic}
{\sc J.~M. Sanz-Serna}, {\em Symplectic {Runge--Kutta} and related methods:
  recent results}, Physica D, 60 (1992), pp.~293--302.

\bibitem{saca93}
{\sc J.~M. Sanz-Serna and M.~P. Calvo}, {\em Numerical Hamiltonian Problems},
  Chapman and Hall, 1993.

\bibitem{Sandu_2019_MR-GARK_High-Order}
{\sc A.~Sarshar, S.~Roberts, and A.~Sandu}, {\em Design of high-order decoupled
  multirate {GARK} schemes}, SIAM Journal on Scientific Computing, 41 (2019),
  pp.~A816--A847.

\bibitem{Tanner_2018_PhD}
{\sc G.M. Tanner}, {\em Generalized additive {Runge-Kutta} methods for stiff
  {ODE}s}, PhD thesis, University of Iowa, 1988.

\bibitem{Yoshida_1990_splitting}
{\sc H.~Yoshida}, {\em Construction of higher order symplectic integrators},
  Physics Letters, 150 (1990), pp.~262--268.

\bibitem{zanna20dvm}
{\sc A.~Zanna}, {\em Discrete variational methods and symplectic generalized
  additive {R}unge---{K}utta methods}.
\newblock https://arxiv.org/abs/2001.07185, January 2020.

\end{thebibliography}
\end{document}